%% file: icat-a1.tex
\documentclass[11pt,reqno]{amsart}

\input{macros}

\title[Nil-Brauer category and $\mathrm{i}$quantum groups]{Nil-Brauer categorifies the split $\mathrm{i}$quantum group of rank one}

\author{Jonathan Brundan}
\address[J.B.]{
  Department of Mathematics,
  University of Oregon,
  Eugene, OR, USA
}
\urladdr{\href{https://pages.uoregon.edu/brundan}{https://pages.uoregon.edu/brundan}, \textrm{\textit{ORCiD}:} \href{https://orcid.org/0009-0009-2793-216X}{https://orcid.org/0009-0009-2793-216X}}
\email{brundan@uoregon.edu}

\author{Weiqiang Wang}
\address[W.W.]{Department of Mathematics, University of Virginia, Charlottesville, VA, USA}
\urladdr{\href{
https://uva.theopenscholar.com/weiqiang-wang}{
https://uva.theopenscholar.com/weiqiang-wang}}
\email{ww9c@virginia.edu}

\author{Ben Webster}
\address[B.W.]{Department of Pure Mathematics, University of Waterloo \& Perimeter Institute for Theoretical Physics,
Waterloo, ON, Canada}
\urladdr{\href{https://uwaterloo.ca/scholar/b2webste}{https://uwaterloo.ca/scholar/b2webste}, \textrm{\textit{ORCiD}:} \href{https://orcid.org/0000-0003-1896-5540}{https://orcid.org/0000-0003-1896-5540}}
\email{ben.webster@uwaterloo.ca}

\thanks{J.B. is supported in part by NSF grant DMS-2101783. W.W. is supported in part by the NSF grant DMS-2001351. B.W. is supported by Discovery Grant RGPIN-2018-03974 from the Natural Sciences and Engineering Research Council of Canada. This research was also supported by Perimeter Institute for Theoretical Physics. Research at Perimeter Institute is supported by the Government of Canada through the Department of Innovation, Science and Economic Development and by the Province of Ontario through the Ministry of Research and Innovation.
Part of this research was carried out at the program ``Representation theory, combinatorics and geometry'' at the Institute for Mathematical Sciences, National University of Singapore, in 2022.
}


\begin{document}

\begin{abstract}
We prove that the Grothendieck ring of the monoidal category of finitely generated graded projective modules for the nil-Brauer category is isomorphic to an integral form of the split iquantum group of rank one. Under this isomorphism, the indecomposable graded projective modules correspond to the icanonical basis. We also derive character formulae for irreducible graded modules and deduce various branching rules.
\end{abstract}

\maketitle
\tableofcontents
\input{s1-introduction}
\input{s2-numerology}
\input{s3-nilbrauer}
\input{s4-idempotents}
\input{s5-representations}
\bibliographystyle{alphaurl}
\bibliography{icat-a1}
\end{document}

%% file: macros.tex
\usepackage{
  mathabx,
  amsmath,
  amsfonts,
  amssymb,
  amsthm,
  amscd,
  graphicx,
  mathtools,
  enumitem,
  mathdots,
  txfonts,
  stackrel
}
\usepackage[all]{xy}
\usepackage[dvipsnames]{xcolor}
\makeatletter
\@namedef{subjclassname@2020}{%
  \textup{2020} Mathematics Subject Classification}
\makeatother

\usepackage{dsfont}
\usepackage[colorlinks=true, linkcolor=black, citecolor=black, urlcolor=purple, breaklinks=true]{hyperref}

\usepackage{tikz}
\usetikzlibrary{decorations.markings}
\usetikzlibrary{calc}
\usetikzlibrary{cd}
\usetikzlibrary{patterns}
\usetikzlibrary{shapes}

\newcommand*\getscale[1]{%
  \begingroup
    \pgfgettransformentries{\scaleA}{\scaleB}{\scaleC}{\scaleD}{\whatevs}{\whatevs}%
    \pgfmathsetmacro{#1}{sqrt(abs(\scaleA*\scaleD-\scaleB*\scaleC))}%
    \expandafter
  \endgroup
  \expandafter\def\expandafter#1\expandafter{#1}%
}

\tikzset{anchorbase/.style={outer sep=auto,baseline={([yshift=-0.5ex]current bounding box.center)}}}
\tikzset{wipe/.style={white,line width=4pt}}
\colorlet{objcolor}{OliveGreen}
\newcommand\botlabel[1]{node[inner sep=0.5pt,anchor=north] {$\color{objcolor}\scriptstyle{#1}$}}

\newcommand{\txtdot}{\begin{tikzpicture}[anchorbase,scale=1.2]
\draw[-] (0,-.2) to (0,.2);
\closeddot{0,0};
\end{tikzpicture}\!\!}

\newcommand{\txtcrossing}{\begin{tikzpicture}[anchorbase,scale=1.2]
\draw[-] (-.15,-.2) to (.15,.2);
\draw[-] (.15,-.2) to (-.15,.2);
\end{tikzpicture}}

\newcommand{\txtcap}{\begin{tikzpicture}[anchorbase,scale=1.2]
\draw[-] (-.15,-.15) to [out=90,in=180] (0,.15) to [out=0,in=90] (.15,-.15);
\end{tikzpicture}}

\newcommand{\txtcup}{\begin{tikzpicture}[anchorbase,scale=1.2]
\draw[-] (-.15,.15) to [out=-90,in=180] (0,-.15) to [out=0,in=-90] (.15,.15);
\end{tikzpicture}}

\newcommand{\smalltxtdot}{\begin{tikzpicture}[anchorbase,scale=.8]
\draw[-] (0,-.2) to (0,.2);
\closeddot{0.01,0};
\end{tikzpicture}\!\!}

\newcommand{\smalltxtcrossing}{\begin{tikzpicture}[anchorbase,scale=.8]
\draw[-] (-.15,-.2) to (.15,.2);
\draw[-] (.15,-.2) to (-.15,.2);
\end{tikzpicture}}

\newcommand{\smalltxtcap}{\begin{tikzpicture}[anchorbase,scale=.8]
\draw[-] (-.15,-.15) to [out=90,in=180] (0,.15) to [out=0,in=90] (.15,-.15);
\end{tikzpicture}}

\newcommand{\smalltxtcup}{\begin{tikzpicture}[anchorbase,scale=.8]
\draw[-] (-.15,.15) to [out=-90,in=180] (0,-.15) to [out=0,in=-90] (.15,.15);
\end{tikzpicture}}

\newcommand{\stringlabel}[1]{\scriptstyle\color{blue}#1}

\newcommand{\cross}[2][0]{
\draw[line width=1.3pt,rotate=#1] ($#2+(-.077,-.077)$) -- ($#2+(.077,.077)$);
\draw[line width=1.3pt,rotate=#1] ($#2+(-.077,.077)$) -- ($#2+(.077,-.077)$);
}
\newcommand{\smallercross}[2][0]{
\draw[line width=1.3pt,rotate=#1] ($#2+(-.066,-.066)$) -- ($#2+(.066,.066)$);
\draw[line width=1.3pt,rotate=#1] ($#2+(-.066,.066)$) -- ($#2+(.066,-.066)$);
}
\newcommand{\smallestcross}[2][0]{
\draw[line width=1.3pt,rotate=#1] ($#2+(-.058,-.058)$) -- ($#2+(.058,.058)$);
\draw[line width=1.3pt,rotate=#1] ($#2+(-.058,.058)$) -- ($#2+(.058,-.058)$);
}

\tikzset{bulb/.style={thin,black,fill=black}}

\renewcommand{\circledbar}[2]{
\getscale{\scalefactor};
\node[circle,draw,bulb,inner sep=.85pt] at (#1) {$\hspace{0.06pt}\color{white}\scriptscriptstyle \pmb#2$};\draw[line width=1.1pt,white] ($(#1)+({-.06/ \scalefactor},{.08 / \scalefactor})$) to ($(#1)+({.06/ \scalefactor},{.08 / \scalefactor})$)}

\newcommand{\circled}[2]{\node[circle,draw,bulb,inner sep=.85pt] at (#1) {$\hspace{0.06pt}\color{white}\scriptscriptstyle\pmb #2$}}

\newcommand{\pinX}[5][black]{
\path #2 node[inner sep=0pt] (x) {$\color{#1}\bullet$};
\path #4 node[rectangle,rounded corners,draw,#1!60!white,fill=#1!10!white,inner sep=2.5pt](y) {$\color{#1}\scriptstyle#5$};
\draw[Triangle Cap-,thick,#1!60!white] (x)-- node[above,inner sep=1pt]{$\color{#1!60!white}\scriptstyle#3$}(y);
}

\newcommand\closeddot[1]{\node at (#1) {$\bullet$}}

\usepackage{zref-clever}

\zcsetup{cap,nameinlink=true,
countertype = {enumi=item, enumii=item, enumiii=item, enumiv=item},
counterresetby = {enumii=enumi, enumiii=enumii, enumiv=enumiii}}
    
\zcRefTypeSetup{item}{
    Name-sg = {} ,
    name-sg = {} ,
    Name-pl = {} ,
    name-pl = {} ,
}
\zcRefTypeSetup{corollary}{
    Name-sg = Corollary ,
    name-sg = corollary ,
    Name-pl = Corollaries ,
    name-pl = corollaries ,
}
\zcRefTypeSetup{definition}{
    Name-sg = Definition ,
    name-sg = definition ,
    Name-pl = Definitions ,
    name-pl = definitions ,
}
\zcRefTypeSetup{example}{
    Name-sg = Example ,
    name-sg = example ,
    Name-pl = Examples ,
    name-pl = examples ,
}
\zcRefTypeSetup{equation}{
    Name-sg = {} ,
    name-sg = {} ,
    Name-pl = {} ,
    name-pl = {} ,
}
\zcRefTypeSetup{lemma}{
    Name-sg = Lemma ,
    name-sg = lemma ,
    Name-pl = Lemmas ,
    name-pl = lemmas ,
}
\zcRefTypeSetup{remark}{
    Name-sg = Remark ,
    name-sg = remark ,
    Name-pl = Remarks ,
    name-pl = remarks ,
}
\zcRefTypeSetup{conjecture}{
    Name-sg = Conjecture ,
    name-sg = conjecture ,
    Name-pl = Conjectures ,
    name-pl = conjectures ,
}
\zcsetup{countertype={subsection=subsection}}
\zcRefTypeSetup{subsection}{
    Name-sg = Subsection ,
    name-sg = subsection ,
    Name-pl = Subsections ,
    name-pl = subsections ,
}
\zcRefTypeSetup{table}{
    Name-sg = Table ,
    name-sg = table ,
    Name-pl = Tables ,
    name-pl = tables ,
}
\zcRefTypeSetup{theorem}{
    Name-sg = Theorem ,
    name-sg = theorem ,
    Name-pl = Theorems ,
    name-pl = theorems ,
}
\zcsetup{countertype={alphatheorem=alphatheorem}}
\zcRefTypeSetup{alphatheorem}{
    Name-sg = Theorem ,
    name-sg = theorem ,
    Name-pl = Theorems ,
    name-pl = theorems ,
}

\newcommand{\cref}[1]{\zcref{#1}}
\newcommand{\Cref}[1]{\zcref[S]{#1}}

\AddToHook{env/conjecture/begin}{\zcsetup{countertype={theorem=conjecture}}}
\AddToHook{env/corollary/begin}{\zcsetup{countertype={theorem=corollary}}}
\AddToHook{env/definition/begin}{\zcsetup{countertype={theorem=definition}}}
\AddToHook{env/example/begin}{\zcsetup{countertype={theorem=example}}}
\AddToHook{env/lemma/begin}{\zcsetup{countertype={theorem=lemma}}}
\AddToHook{env/remark/begin}{\zcsetup{countertype={theorem=remark}}}
\AddToHook{env/table/begin}{\zcsetup{countertype={theorem=table}}}

\newcommand\Z{\mathbb{Z}}
\newcommand\Q{\mathbb{Q}}

\newcommand\N{\mathbb{N}}
\def\K{w}
\def\B{\mathcal{B}}
\def\I{\mathcal{I}}
\def\J{\mathcal{J}}
\newcommand\kk{\Bbbk}
\newcommand{\sqbinom}[2]{\genfrac{[}{]}{0pt}{}{#1}{#2}}
\def\O{\mathbb{O}}
\def\op{{\operatorname{op}}}
\def\rev{{\operatorname{rev}}}
\newcommand\one{\mathds{1}}

\def\sl{\mathfrak{sl}}
\def\RSD{\mathrm{D}}
\def\Par{\mathcal{P}}
\def\e{\mathbf{e}}
\def\f{\mathbf{f}}

\def\u{\mathbf{u}}
\def\vv{\mathbf{v}}
\def\w{\mathbf{w}}
\def\x{\mathbf{x}}
\def\y{\mathbf{y}}

\newcommand\NB{\mathrm{NB}}
\newcommand\ZZ{\mathrm{Z}}

\newcommand\cA{\mathbf{A}}
\newcommand\cB{\mathbf{B}}
\newcommand\cNB{\mathbf{NB}}
\newcommand\CNB{\mathbf{cNB}}

\newcommand\fU{\mathfrak{U}}
\newcommand\cV{\mathbf{V}}
\newcommand\tT{{\scriptstyle\mathrm{T}}}
\newcommand\tR{{\scriptstyle\mathrm{R}}}

\newcommand\tD{{\scriptstyle\mathrm{D}}}
\newcommand{\ndots}{\bullet}
\newcommand\gVec{\mathbf{gVec}}
\def\gMOD{\operatorname{\!-gMod}}
\def\gMod{\operatorname{\!-gmod}}
\def\pgMod{\operatorname{\!-pgmod}}
\def\END{\mathbf{End}}
\def\gEND{\mathbf{gEnd}}
\def\rank{\operatorname{rank}}
\def\lround{(\!(}
\def\rround{)\!)}

\DeclareMathOperator{\End}{End}

\DeclareMathOperator{\Hom}{Hom}
\def\soc{\operatorname{soc}}
\def\head{\operatorname{hd}}
\def\rad{\operatorname{rad}}
\DeclareMathOperator{\id}{id}
\DeclareMathOperator{\im}{im}
\DeclareMathOperator{\coker}{coker}

\DeclareMathOperator{\Res}{Res}
\DeclareMathOperator{\Ind}{Ind}

\def\ch{\operatorname{ch}}
\def\ostar{{\medcirc\hspace{-2.22mm}\star\hspace{.5mm}}}

\def\atled{\text{\rotatebox[origin=c]{180}{$\delta$}}}

\def\X{\mathrm{X}}
\def\Y{\mathrm{Y}}
\def\H{\mathrm{H}}
\def\HH{\mathring{\mathrm{H}}}
\newcommand{\A}{{\Z[q,q^{-1}]}}
\newcommand{\U}{{\mathrm U}}

\newcommand{\Ui}{{\mathrm U}_t^{\imath}}
\newcommand{\Uinone}{{\mathrm U}^{\imath}}
\newcommand{\Uizero}{{\mathrm U}_0^{\imath}}
\newcommand{\Uione}{{\mathrm U}_1^{\imath}}
\newcommand{\UAi}{{_\Z}{\mathrm U}_t^{\imath}}

\def\pmod#1{{\:\left(\operatorname{mod} #1\right)}}

\def\el{\ell}

\newtheorem{alphatheorem}{Theorem}

\newtheorem{theorem}{Theorem}[section]
\newtheorem{lemma}[theorem]{Lemma}
\newtheorem*{lemma*}{Lemma}
\newtheorem{corollary}[theorem]{Corollary}

\theoremstyle{definition}
\newtheorem{definition}[theorem]{Definition}
\newtheorem{remark}[theorem]{Remark}
\newtheorem{example}[theorem]{Example}

\numberwithin{equation}{section}
\allowdisplaybreaks

%% file: s1-introduction.tex
\setcounter{section}{0}

\section{Introduction}\label{intro}

In \cite{letzter1},
Letzter introduced what we now call 
the {\em iquantum groups} associated to symmetric pairs.
These can be viewed as a
generalization of Drinfeld-Jimbo quantum 
groups---the latter are the
iquantum groups 
arising from diagonal symmetric pairs.
Lusztig's canonical bases for quantum groups, with their favorable positivity properties, 
provided one source of motivation for the categorification of quantum groups via the 
{\em Kac-Moody 2-categories}
of Khovanov, Lauda and Rouquier \cite{KL3, Rou}. A theory of {\em icanonical bases} for iquantum groups was  developed in \cite{BW18QSP, BW18KL}. In special cases,
these again have positive structure constants; see \cite{LiW18} which treats the quasi-split types AIII.
Therefore, it is reasonable to hope that there should be a categorification of iquantum groups. 

In rank one, there are three quasi-split
iquantum groups of finite type. First, there is the usual $\mathrm U_q(\mathfrak{sl}_2)$, which was
categorified by Lauda and Rouquier in \cite{Lauda,Rou}.
The second, arising from the Satake diagram of $A_2$ with 
non-trivial diagram involution, was categorified in \cite{BSWW}. 
In this article, we explain how to categorify the remaining case, the split iquantum group 
$\mathrm{U}_q^\imath(\sl_2)$ 
corresponding to the symmetric pair $(\mathrm{SL}_2,\mathrm{SO}_2)$.
This is a basic building block for general iquantum groups,
and it is expected to play a key role in the categorification of quasi-split iquantum groups of higher rank.

Our
categorification of $\mathrm{U}_q^\imath(\sl_2)$  arises from the
{\em nil-Brauer category} $\cNB_t$ introduced in \cite{BWWbasis}.
This is a strict graded
$\kk$-linear monoidal category
defined over a field $\kk$ of characteristic different from 2. 
It has one self-dual generating object $B$ and four generating morphisms
represented diagrammatically by
$\txtdot\;$ (degree 2), $\;\txtcrossing\;$ (degree $-2$), $\;\txtcap\;$ (degree 0), and $\;\txtcup\;$ (degree 0), subject to some natural relations recorded
in \cref{NBdef}.
The parameter $t$ gives the value of
$\;\begin{tikzpicture}[anchorbase]
\draw (0,0) circle (.2);
\end{tikzpicture}\,:\one\rightarrow\one$. The defining relations imply that $t^2 1_\one = t 1_\one$, hence, 
we must have that $t=0$ or $t=1$ in order for the category to be non-trivial;
see \cite[(2.9)]{BWWbasis}. We assume that this is the case from now on.

To formulate the main results precisely, 
rather than working in terms of idempotents, as is often done in the categorification literature, 
we use the language of modules.
By a {\em graded $\cNB_t$-module},
we mean a graded $\kk$-linear functor
from $\cNB_t$ to graded vector spaces.
The endofunctor of $\cNB_t$
defined by tensoring with its generating object extends
to an exact endofunctor, also denoted $B$, of the category of graded $\cNB_t$-modules.
Let $[n] := q^{n-1}+q^{n-3}+\cdots+q^{1-n}$
be the quantum integer, 
and $V^{\oplus [n]}$ denote the corresponding
direct sum of degree-shifted copies of a graded module $V$.

\begin{alphatheorem}
\label{introthm:PIM}
There are unique
(up to isomorphism) indecomposable 
projective graded $\cNB_t$-modules $P(n)\:(n \geq 0)$
such that $P(0)$ is the projective graded module associated to the
idempotent that is the identity endomorphism of the unit object, and for $n \geq 0$ we have that
$$
B P(n) \cong 
\begin{dcases}
P(n+1)^{\oplus [n+1]} \oplus P(n-1)^{\oplus [n]}
&\text{if $n \equiv t \pmod{2}$}\\
P(n+1)^{\oplus [n+1]}
&\text{if $n \not\equiv t \pmod{2}$,}
\end{dcases}
$$
interpreting $P(-1)$ as 0.
These modules give a full set of indecomposable
projective graded $\cNB_t$-modules (up to isomorphism and grading shift). 
\end{alphatheorem}

The proof of \cref{introthm:PIM} 
is similar in spirit to 
Lauda's proof of the analogous result for the 2-category $\fU(\sl_2)$ obtained in \cite{Lauda}.
It involves the explicit construction of appropriate homogeneous primitive idempotents. These resemble primitive idempotents in the nil-Hecke algebra familiar from Schubert calculus, but they are considerably more subtle; see \cref{id3,idempotenttheorem}.
Another important ingredient 
needed to establish the indecomposability of $P(n)$ is the identification of the Cartan form 
on the Grothendieck ring of $\cNB_t$ 
with an explicitly defined sesquilinear form on the iquantum group. This is discussed further after the statement of the next theorem, which is our main categorification result.

Let $\U^\imath := \mathrm{U}_q^\imath(\sl_2)$ be the split iquantum group of rank 1.
As a $\Q(q)$-algebra, this is simply a polynomial algebra on one generator $b$,
but it has a 
non-trivial $\Z[q,q^{-1}]$-form
$\UAi$ associated to the parameter $t \in \{0,1\}$.
As a $\Z[q,q^{-1}]$-module, $\UAi$ is free 
with a distinguished basis
given by the {\em idivided powers} $b^{(n)}$. These arise from the {\em icanonical basis} of $\UAi$ constructed in \cite{BW18KL}
in terms of the
finite-dimensional irreducible $\sl_2$-modules of highest weight
$\lambda\equiv t\pmod{2}$, and computed explicitly in the split
rank one case in \cite{BeW18}.
The recursion for the 
indecomposable projective graded modules 
in \cref{introthm:PIM} exactly matches the recurrence relation for
idivided powers $b^{(n)}\:(n \geq 0)$ from \cite{BeW18}.
This coincidence is the essence of our next main theorem; see \cref{therealthing}.
For the statement, let $K_0(\cNB_t)$ be the split Grothendieck ring of the monoidal
category of finitely generated projective graded $\cNB_t$-modules. This is a $\Z[q,q^{-1}]$-algebra,
with the action of $q$ arising from the grading shift functor.

 \begin{alphatheorem}
\label{introthm:iso}
There is a unique $\Z[q,q^{-1}]$-algebra isomorphism
$$
\kappa_t:
K_0(\cNB_t)
\stackrel{\sim}{\rightarrow} \UAi
$$
intertwining the endomorphism of $K_0(\cNB_t)$
induced by the endofunctor $B$ with the
endomorphism of $\UAi$ defined by multiplication by the generator $b$ of the iquantum group.
For any $n \geq 0$, $\kappa_t$ maps
the isomorphism class of the indecomposable projective module $P(n)$
to the icanonical basis element $b^{(n)}$.
\end{alphatheorem}

Under the isomorphism of \cref{introthm:iso},
the non-degenerate symmetric 
bilinear form $(\cdot,\cdot)^\imath$
on $\UAi$ constructed in \cite{BW18QSP}
is equal (after twisting with the bar involution to make it 
sesquilinear in the appropriate sense, and some rescaling) to the Cartan form on 
$K_0(\cNB_t)$.
The proof of this depends ultimately on the basis theorem for $\cNB_t$
from \cite{BWWbasis} together with some combinatorics of chord diagrams which is of independent interest; see \cref{swim,classic,toinfinityandbeyond}.

The remaining results in the article
rely on the observation that the category of 
graded $\cNB_t$-modules
has some useful additional
structure: it is an {\em affine lowest weight category} in a suitably generalized sense. In particular, there are certain graded $\cNB_t$-modules
$\Delta(n)$, $\nabla(n)$, $\bar\Delta(n)$ and $\bar\nabla(n)$, the
{\em standard}, {\em costandard}, {\em proper standard modules} and
{\em proper costandard modules}, all of which
are equipped with explicit bases.
The proper standard module $\bar\Delta(n)$ has a unique irreducible
quotient denoted $L(n)$,
which is also the unique irreducible submodule of $\bar\nabla(n)$.
The modules $L(n)\:(n \geq 0)$ give a complete set of graded
irreducible $\cNB_t$-modules up to isomorphism and grading
shift. There is a graded analog of the usual BGG reciprocity identifying
certain standard filtration multiplicities $(P(n):\Delta(m))_q$ with
the graded decomposition multiplicities $[\bar\Delta(m):L(n)]_q$; see \cref{wagner}. These assertions follow from an application of the general machinery of {\em graded triangular bases} developed in \cite{GTB}---the nil-Brauer category is a perfect example for this theory.

The minimal standard modules $\Delta(0)$
and $\Delta(1)$ are projective and therefore coincide with $P(0)$ and $P(1)$, respectively,
but after that the two families of modules diverge. At the decategorified level, the standard modules $\Delta(n)$ and costandard modules $\nabla(n)$ correspond to the {\em standard basis} $\delta_n\:(n \geq 0)$ and
the {\em costandard basis} $\atled_n\:(n \geq 0)$ for $\Uinone$, both of which are introduced in \cref{sect2}.
These two bases are interchanged by the bar involution, and the costandard basis
is an orthogonal basis with respect to the form $(\cdot,\cdot)^\imath$.
The standard basis elements satisfy the following recurrence relation:
\begin{align*}
\delta_0 &= 1,&
b \delta_n &= [n+1]\delta_{n+1} +
\frac{q^{1-n}}{1-q^{2}}
\delta_{n-1},
\end{align*}
interpreting $\delta_{-1}$ as $0$.
This should be compared with the following theorem describing the
effect of the endofunctor $B$ on a
standard module:

 \begin{alphatheorem} 
\label{introthm:ses}
For $n \geq 0$, there is a short exact sequence
of graded $\cNB_t$-modules
$$
0 \longrightarrow 
\bigoplus_{i \geq 0} 
q^{2i+1-n} \Delta(n-1)
\longrightarrow B \Delta(n)
\longrightarrow 
\Delta(n+1)^{\oplus [n+1]}
\longrightarrow 0.
$$
In the first term, $q$ denotes the upward grading shift functor, and this term should be interpreted as 0 in case $n=0$.
\end{alphatheorem}

An interesting feature of \cref{introthm:ses} is the presence of the infinite direct sum in the first term of the short exact sequence---the finitely generated $\cNB_t$-modules
$B \Delta(n)\:(n > 0)$ are {\em not} Noetherian. This corresponds to
the fact that the PBW basis $\delta_n\:(n \geq 0)$ is a basis for
$\Uinone$ over $\Q(q)$, but not for $\UAi$ over $\Z[q,q^{-1}]$. 
\cref{introthm:ses} is proved in
\cref{burger} in the main body of the text.
There is also a parallel result for proper standard modules; see \cref{burger2}.

For closed formulae for the transition matrices between the bases $b^{(m)}\:(m \geq 0)$
and $\delta_n\:(n\geq 0)$, see \cref{cornucopia}.
Translating to representation theory and using BGG reciprocity, we obtain the following explicit formula for graded decomposition numbers:

\begin{alphatheorem}\label{introthm:decomposition}
The irreducible subquotients of the
proper standard module $\bar\Delta(n)\:(n \geq 0)$
are isomorphic (up to grading shifts)
to $L(n+2m)$ for $m \geq 0$ with
\begin{align*}
    [\bar\Delta(n): L(n+2m)]_q  &= 
\begin{dcases}
q^{m(2m-1)}\:\big/\:(1-q^{4}) (1-q^{8})\cdots (1-q^{4m})
    &\text{if $n \equiv t\pmod{2}$}\\
q^{m(2m+1)}\:\big/\:(1-q^{4}) (1-q^{8})\cdots (1-q^{4m})
    &\text{if $n\not\equiv t\pmod{2}$.}
\end{dcases}
\end{align*}
\end{alphatheorem}

To formulate one more such combinatorial result, for a finitely generated 
graded $\cNB_t$-module $V$,
its {\em graded character} is the
series $$
\ch V = \sum_{n \geq 0} \dim_q(1_n V) \chi^n
\in \N\lround q\rround\llbracket \chi\rrbracket
$$
where $\chi$ is a formal variable and
$\dim_q(1_n V) \in \N\lround q \rround$
is the graded dimension of the
graded vector space obtained by evaluating the functor $V$ on the $n$th tensor power of the generating object $B$.

 \begin{alphatheorem}
\label{introthm:char}
For $n \geq 0$, we have that
$$
\ch L(n) = 
[n]!\, \chi^n\:\;\Bigg /\!
\prod_{\substack{1 \leq k \leq n+1 \\ k \equiv t \pmod{2}}}
(1-[k]^2 \chi^2) \:\:\in \N[q,q^{-1}]\llbracket \chi\rrbracket.
$$
\end{alphatheorem}

Finally,
we also prove {\em branching rules} which give complete information about 
the structure of the modules
$B L(n)\:(n \geq 0)$;
see \cref{lowerbound}.
Except in the case that
$n=t=0$ (when it is zero),
these branching rules show that $B L(n)$ is a self-dual uniserial
module with irreducible socle and cosocle isomorphic
(up to appropriate grading shifts) to $L(n-1)$ if $n \equiv t\pmod{2}$ or to $L(n+1)$ if $n \not \equiv t \pmod{2}$. 
Moreover,
$$
\End_{\cNB_t}(B L(n)) \cong \kk[x] / \big(x^{\beta(n)}\big)
$$
where $\beta(n) = n$ if $n \equiv t \pmod{2}$
or $n+1$ if $n \not\equiv t\pmod{2}$.
The combinatorics arising here is the same as the combinatorics of the underlying icrystal basis 
described in \cite[Ex. 4.1.4]{Wat23}.

\vspace{2mm}
\noindent{\em General conventions.}
Throughout the article, $t \in \{0,1\}$ will be a fixed parameter.
Given also $n \in \N$, we use the shorthand
$\delta_{n \equiv t}$
to denote 1 if $n \equiv t\pmod{2}$ or 0 otherwise. Similarly,
$\delta_{n \not\equiv t}$
denotes 1 if $n \not\equiv t\pmod{2}$ or 0 
otherwise.
We write $S_n$ for the symmetric group on $n$ letters. Let $s_i \in S_n$ be the simple transposition
$(i\:\:i\!+\!1)$, let 
$\ell:S_n \rightarrow \N$ be
the associated length function, and
let $w_n$ be the longest element of $S_n$.
 We denote the category of graded vector spaces over the field $\kk$ by $\gVec$, using $q$ for the {\em upward} grading shift functor.
So, for a graded vector space $V=\bigoplus_{d \in \Z} V_d$,
its grading shift $qV$ is the same underlying vector space
with new grading defined via
$(q V)_d := V_{d-1}$ for each $d \in \Z$.
For a graded vector space $V = \bigoplus_{d \in \Z} V_d$ with finite-dimensional graded pieces, we define its {\em graded dimension} to be
\begin{equation}
\dim_q V := \sum_{d \in \Z} (\dim V_d) q^{d}.
\end{equation}
For any formal series
$f(q) = \sum_{d \in \Z} a_d q^d$ with each $a_d \in \N$, we write $V^{\oplus f(q)}$
for $\bigoplus_{d \in \Z} q^d V^{\oplus a_d}$.
Also $\overline{f(q)}$ denotes $f(q^{-1})$.

%% file: s2-numerology.tex
\setcounter{section}{1}

\section{Bases of the split iquantum group of rank one}\label{sect2}

In this section, we recall some basic facts about the split iquantum group of rank 1 following 
\cite{BW18KL, BeW18}.  Then we introduce 
a new PBW-type basis, and
derive combinatorial formulae
for various 
transition matrices, including
between the PBW basis 
and the icanonical basis.
For all of this, we work over the field
$\Q(q)$ for an indeterminate $q$.
We write $[n]$ for the quantum
integer $\frac{q^n-q^{-n}}{q-q^{-1}}$,
$[n]!$ for the quantum factorial,
and $\sqbinom{n}{r}:=[n][n-1]\cdots[n-r+1] / [r]!$. The word {\em anti-linear} always means with respect to the bar involution $-:\Q(q)\rightarrow \Q(q)$ that is the field automorphism taking $q$ to $q^{-1}$.
We denote the limit of a convergent
 sequence $(f_\lambda)_{\lambda \geq 0}$
in $\Q\lround q^{-1} \rround$
by $\lim_{\lambda \rightarrow \infty} f_\lambda$.

\subsection{Quantum groups}
Our general conventions for quantum groups are the same as in \cite{Lubook}, except that
we write $q$ in place of Lusztig's $v$.
Let $\f$ be the polynomial algebra over $\Q(q)$ generated by one
element $\theta$. We write $\theta^{(n)}$ for the {\em divided power}
$\theta^n / [n]!$.
This is a fixed point for the {\em bar involution} $\psi:\f
\rightarrow \f$, which is the anti-linear involution defined from
$\psi(\theta) := \theta$.
Let $(\cdot,\cdot):\f
\times \f \rightarrow \Q(q)$ be the non-degenerate symmetric bilinear
form from
\cite[Sec.~1.2.5]{Lubook}. It satisfies
\begin{equation}\label{midnight}
\big( \theta^{(m)},\theta^{(n)}\big)
= \frac{\delta_{m,n}}{(1-q^{-2}) (1-q^{-4})
\cdots (1-q^{-2n})}
\end{equation}
for $m, n\geq 0$. 
Let $R:\f\rightarrow \f$
be the linear map defined by
\begin{align}\label{taxes}
R(1) &= 0,
&
R\big(\theta^{(n)}\big)
&= 
\frac{q^{n-1} \theta^{(n-1)}}{1-q^{-2}}
\end{align}
for $n \geq 1$.
This map arises naturally as 
the adjoint of left 
multiplication by $\theta$:
we have that 
\begin{equation}\label{adj1}
( \theta x, y )
=( x, R(y))
\end{equation}
for all $x,y \in \f$.
Equivalently, $R(x) = r(x) / (1-q^{-2})$
where $r$ is the map defined in either the first or the second paragraph of \cite[Sec.~1.2.13]{Lubook} (the two maps coincide
in rank one).

The quantum group 
$\U={\mathrm U}_q(\mathfrak{sl}_2)$ is the
$\Q(q)$-algebra
with generators $e, f, k, k^{-1}$
satisfying the relations
\begin{align*}
k e k^{-1} &= q^2 e,&
k f k^{-1} &= q^{-2} f,&
[e,f] &= \frac{k-k^{-1}}{q-q^{-1}}.
\end{align*}
Here, we have switched to using lower case for $e,f$ compared to
\cite{Lubook}
so that we can use the
upper case letters $E, F$ for corresponding functors in
categorification.  The subalgebras of $\U$ generated by
$f$ and by $e$ are denoted $\U^-$ and $\U^+$,
respectively. Both are isomorphic to $\f$ via the maps
$\f \rightarrow \U^+, x \mapsto x^+$ and $\f \rightarrow \U^-, x
\mapsto x^-$ defined so that $\theta^+ := e$ and $\theta^- := f$.
The {\em divided powers}
$e^{(n)} := e^n / [n]!$ and $f^{(n)} := f^n / [n]!$ are the images of
$\theta^{(n)}$ under these maps.
There are various useful symmetries:
\begin{itemize}
\item
Let $\psi:\U \rightarrow \U$
be the usual {\em bar involution} on $\U$, that is, the anti-linear algebra involution which
fixes $e$ and $f$ and takes $k$ to $k^{-1}$.
\item
Let $\rho:\U \rightarrow \U$ 
be the linear algebra anti-involution
such that
$\rho(k) = k$,
$\rho(e) = q kf$,
$\rho(f) = q^{-1} e k^{-1}$.
\end{itemize}
By \cite[Prop.~3.1.6(b)]{Lubook}
(or an easy induction exercise
using \cref{taxes}), we have that
 \begin{align}\label{Ey}
ex^- -x^- e = 
q^{-1} k R(x)^- - q^{-1} 
R(x)^- k^{-1}
\end{align}
for any $x \in \f$.

We denote the 
irreducible $\U$-module of highest weight
$\lambda \in \N$ by $V(\lambda)$.
This is generated by a vector 
$\eta_\lambda$ such that $e \eta_\lambda = 0$
and $k \eta_\lambda = q^\lambda \eta_\lambda$.
There is an anti-linear involution
$\psi_\lambda:V(\lambda)\rightarrow V(\lambda)$
such that $\psi_\lambda(\eta_\lambda) = \eta_\lambda$ and $\psi_\lambda(uv)
= \psi(u) \psi_\lambda(v)$
for $u \in \U, v \in V(\lambda)$.
Also let $(\cdot,\cdot)_\lambda:V(\lambda)\times V(\lambda)\rightarrow \Q(q)$ be 
the unique non-degenerate symmetric bilinear form on $V(\lambda)$
such that
\begin{align}\label{worse}
(\eta_\lambda,\eta_\lambda)_\lambda &= 1, &
(uv_1 ,v_2)_\lambda &= ( v_1, \rho(u) v_2)_\lambda
\end{align}
for $u \in \U, v_1,v_2 \in V(\lambda)$.
The form $(\cdot,\cdot)$
on $\f$ can be recovered from these forms on the modules
$V(\lambda)$ since we have that
\begin{equation}\label{amsterdam}
(x,y)= 
\lim_{\lambda \rightarrow \infty}
\big(x^- \eta_\lambda, y^- \eta_\lambda \big)_\lambda
\end{equation}
for all $x,y \in \f$
by a  
special case of \cite[Prop.~19.3.7]{Lubook}. 
The vectors
$f^{(n)} \eta_\lambda\:(0 \leq n \leq \lambda)$ give the {\em canonical basis} for $V(\lambda)$.
In fact, they give a basis for an
integral form
${_\Z}V(\lambda)$ over $\A$.
The 
anti-involution $\psi_\lambda$ restricts
to an anti-linear involution
of ${_\Z}V(\lambda)$,
and
the values of the form $(\cdot,\cdot)_\lambda$
on elements of 
 ${_\Z}V(\lambda)$ lie in $\Z[q,q^{-1}]$.

For the purposes of categorification, 
one usually replaces $\U$ by its
modified form $\dot\U$,
which is a locally unital algebra $\dot\U = \bigoplus_{\lambda,\mu \in \Z} 1_\mu \dot\U 1_\lambda$ with a distinguished system 
$1_\lambda (\lambda \in \Z)$ of mutually orthogonal idempotents 
replacing the diagonal generators $k, k^{-1}$.
The relationship between $\U$ and $\dot\U$
can be expressed either by saying that $\dot\U$ is a $(\U,\U)$-bimodule, or that
$\U$ embeds into the completion
 of $\dot\U$ consisting of matrices
$(a_{\mu,\lambda})_{\lambda,\mu \in \Z}
\in \prod_{\lambda,\mu\in\Z} 1_\mu \dot \U 1_\lambda$
such that there are only finitely many non-zero entries in each row and column.
The element $k \in \U$ corresponds to the
diagonal matrix with $q^\lambda 1_\lambda$
as its $\lambda$th diagonal entry,
while $e,f \in \U$ are identified with the matrices
whose only non-zero entries are $1_{\lambda+2} e 1_\lambda\:(\lambda \in \Z)$ 
and $1_{\lambda} f 1_{\lambda+2}\:(\lambda \in \Z)$, respectively.

\subsection{The iquantum group and its standard/costandard bases}
The {\em iquantum group}
${\mathrm U}^\imath(\sl_2)$ is the subalgebra $\Uinone$ of $\U$ 
 generated by 
\begin{equation}\label{gradingmagic}
b := f +\rho(f) = f+q^{-1} e k^{-1}.
\end{equation}
As an algebra, $\Uinone$ is uninteresting since
it is the free $\Q(q)$-algebra on $b$. However it is an interesting coideal subalgebra of $\U$ for an appropriate choice of comultiplication.

The symmetry $\rho$ of $\U$ restricts to 
a linear anti-involution 
$\rho:\Uinone\rightarrow \Uinone$
with $\rho(b) = b$.
Also, the {\em bar involution}
$\psi^\imath:\Uinone\rightarrow \Uinone$
is the unique anti-linear involution
such that $\psi^\imath(b) = b$.
We stress a key point: $\psi^\imath$ 
is {\em not} the restriction of the bar involution $\psi$ on $\U$, indeed, the latter does not leave  $\Uinone$ invariant.
For $\lambda\in\N$,
there is a unique anti-linear involution
$\psi_\lambda^\imath:V(\lambda)\rightarrow V(\lambda)$
such that 
\begin{align}\label{better}
\psi_\lambda^\imath(\eta_\lambda) &= \eta_\lambda,&
\psi_\lambda^\imath(uv) &=
\psi^\imath(u) \psi_\lambda^\imath(v)
\end{align}
for all 
$u \in \Uinone, v \in V(\lambda)$; see \cite[Cor. 3.11]{BW18KL} and \cite[Prop.~5.1]{BW18QSP}. 
Also, by \cite[Lem.~6.25]{BW18QSP},
there is a symmetric 
bilinear form $(\cdot,\cdot)^\imath:
\Uinone\times\Uinone \rightarrow \Q(q)$
such that
\begin{equation}\label{otherform}
( u_1, u_2)^\imath
= \lim_{\lambda\rightarrow \infty}
\big(u_1 \eta_\lambda, 
u_2 \eta_\lambda \big)_\lambda
\end{equation}
for all $u_1,u_2 \in \Uinone$.
From \cref{worse}, we get that
\begin{align}\label{eeee}
(b u_1, u_2)^\imath
&= (u_1, b u_2)^\imath
\end{align}
for any $u_1,u_2 \in \Uinone$.
In \cite[Th.~6.27]{BW18QSP},
it is shown that
$(\cdot,\cdot)^\imath$ is non-degenerate.
This also follows from the following
theorem together with the non-degeneracy
of the form $(\cdot,\cdot)^-$
on $\U^-$.

\begin{theorem}\label{jy2}
There is a unique 
isomorphism of $\Q(q)$-vector spaces
$\jmath:\Uinone \stackrel{\sim}{\rightarrow} \f$
such that
\begin{equation}\label{late}
\lim_{\lambda \rightarrow \infty}
\big(u \eta_\lambda, x^-\eta_\lambda\big)_\lambda
=(\jmath(u),x)
\end{equation}
for all $u \in \Uinone$ and $x \in \f$.
Moreover, the following hold for $u,u_1,u_2 \in \Uinone$:
\begin{enumerate}
\item $\jmath(bu) = \theta \jmath(u) + R(\jmath(u))$.
\item
$(u_1,u_2)^\imath
= \big(\jmath(u_1), \jmath(u_2)\big)$.
\end{enumerate}
\end{theorem}

\begin{proof}
Uniqueness of a linear map
$\jmath$ satisfying \cref{late}
follows easily from 
the non-degeneracy of the form $(\cdot,\cdot)$.
To prove existence, we can assume that $u$ is a power of $b$ and proceed by induction on degree.
Let $\jmath(1) := 1$, which clearly satisfies
\cref{late} for all $x \in \f$.
Now assume for some $u \in \Uinone$
that $\jmath(u)$ satisfying
\cref{late} for all $x$ has been constructed
inductively, and consider $\jmath(bu)$.
Using \cref{worse} and the identity \cref{Ey} multiplied on the left by $q k^{-1}$, we have that
\begin{align*}
\lim_{\lambda \rightarrow \infty} \big(
b u \eta_\lambda, x^- \eta_\lambda
\big)_\lambda
&\!\!\stackrel{\cref{worse}}{=}\!
\lim_{\lambda \rightarrow \infty} \big(
u \eta_\lambda, 
b x^- \eta_\lambda
\big)_\lambda
=
\lim_{\lambda \rightarrow \infty} \big(
u \eta_\lambda, f x^- \eta_\lambda
+ q k^{-1} e x^- \eta_\lambda
\big)_\lambda\\
&\!\!\!\stackrel{\cref{Ey}}{=}\!
\lim_{\lambda \rightarrow \infty} \big(
u \eta_\lambda, f x^- \eta_\lambda
+ R(x)^- \eta_\lambda
-
k^{-1} R(x)^- k^{-1} \eta_\lambda
\big)_\lambda\\
&=
\lim_{\lambda \rightarrow \infty} \big(
u \eta_\lambda, (\theta x)^- \eta_\lambda
+ R(x)^- \eta_\lambda
\big)_\lambda\\
&=
(\jmath(u), \theta x + R(x))
\stackrel{\cref{adj1}}{=} 
(\theta \jmath(u)+R(\jmath(u)),x).
\end{align*}
So $\jmath(bu) := \theta \jmath(u) + R(\jmath(u))$
satisfies \cref{late}. This proves the existence of a linear map $\jmath$ satisfying
\cref{late}, and at the same time we have established (1).
To see that $\jmath$ is a linear isomorphism,
it follows easily from (1) that $\jmath(b^n)$
is a monic polynomial of degree $n$ in $\theta$.
Since $\Uinone$ and $\f$
are free on $b$ and on $\theta$, respectively,
it is now clear that $\jmath$ is an isomorphism.

It remains to prove (2).
By the definition \cref{otherform}
and \cref{late}, we need to show that
$$
\lim_{\lambda \rightarrow \infty}
\big(u_1\eta_\lambda, \jmath(u_2)^-\eta_\lambda\big)_\lambda=
\lim_{\lambda \rightarrow \infty}
\big(u_1 \eta_\lambda, u_2 \eta_\lambda\big)_\lambda
$$
for all $u_1,u_2 \in \Uinone$.
Note
that the limit on the left hand side exists by what we have proved so far.
We assume that $u_2$ is a power of $b$ and proceed by induction on its degree.
The base case $u_2 = 1$ is clear.
Now assume the result has been proved for all $u_1$ and some
$u_2$, and consider $b u_2$.
Using (1), we have that
\begin{align*}
\lim_{\lambda \rightarrow\infty}
\big(u_1\eta_\lambda, \jmath(b u_2)^- \eta_\lambda \big)_\lambda&=\;
\lim_{\lambda \rightarrow\infty}
\big(u_1\eta_\lambda, f \jmath(u_2)^- \eta_\lambda + R(\jmath(u_2))^- \eta_\lambda\big)_\lambda\\
&=\;
\lim_{\lambda \rightarrow\infty}
\big(u_1\eta_\lambda, f \jmath(u_2)^- \eta_\lambda + R(\jmath(u_2))^- \eta_\lambda
-k^{-1} R(\jmath(u_2))^- k^{-1}\eta_\lambda\big)_\lambda\\
&\!\!\stackrel{\cref{Ey}}{=}\!
\lim_{\lambda \rightarrow\infty}
\big(u_1\eta_\lambda, f \jmath(u_2)^- \eta_\lambda + q k^{-1} e \jmath(u_2)^- \eta_\lambda \big)_\lambda=
\lim_{\lambda \rightarrow\infty}
\big(u_1\eta_\lambda, b \jmath(u_2)^- \eta_\lambda \big)_\lambda\\
&\!\!\stackrel{\cref{worse}}{=}\!
\lim_{\lambda \rightarrow\infty}
\big(b u_1\eta_\lambda, \jmath(u_2)^- \eta_\lambda \big)_\lambda
=
\lim_{\lambda \rightarrow\infty}
\big(b u_1\eta_\lambda, u_2 \eta_\lambda \big)_\lambda
\stackrel{\cref{worse}}{=}
\lim_{\lambda \rightarrow\infty}
\big(u_1\eta_\lambda, b u_2 \eta_\lambda \big)_\lambda.
\end{align*}
\end{proof}

Applying \cref{jy2}, we let $\atled_n \in \Uinone$ be
the unique element such that
$\jmath(\atled_n) = \theta^{(n)}$.
The elements
$\atled_n\:(n \geq 0)$
give a basis for $\Uinone$, which we call
the {\em costandard basis}.
From \cref{jy2}(2) and \cref{midnight}, 
we get that
\begin{equation}\label{bilform}
\big(
\atled_m,\atled_n\big)^\imath 
= \frac{\delta_{m,n}}{(1-q^{-2})(1-q^{-4})\cdots(1-q^{-2n})}
\end{equation}
for $m,n\geq 0$.
Thus, the costandard basis is an orthogonal basis. 
The following recurrence relation 
 is
easily deduced using \cref{jy2}(1) and \cref{taxes}:
\begin{align}\label{costdrecurrence}
\atled_0 &= 1,&
b \atled_n &= [n+1]\atled_{n+1} +
\frac{q^{n-1}}{1-q^{-2}}
\atled_{n-1}
\end{align}
for $n \geq 0$, interpreting 
$\atled_{-1}$ as $0$.
Applying the bar involution $\psi^\imath$ to the costandard basis
gives another basis $\delta_n := \psi^\imath(\atled_n)\:(n \geq 0)$
for $\Uinone$
which we call the {\em standard basis}. It satisfies the recurrence
\begin{align}\label{stdrecurrence}
\delta_0 &= 1,&
b \delta_n &= [n+1]\delta_{n+1} +
\frac{q^{1-n}}{1-q^{2}}
\delta_{n-1}.
\end{align}
Despite only differing by an application of the bar involution,
we generally prefer to work with $\delta_n$, although $\atled_n$ has
the advantage of being an orthogonal basis.

\begin{remark}
The linear isomorphism in \cref{jy2} is analogous to the
isomorphism $\U^+\otimes \U^- \cong \dot\U 1_\zeta$ in
\cite[Theorem 2.8]{PBW}. The costandard basis for $\Uinone$ with the orthogonality property
\cref{bilform} is analogous to the PBW bases for
modified quantum groups of finite type constructed in \cite{PBW}.
\end{remark}

\subsection{Combinatorics of chord diagrams}\label{chorddiagrams}
Next, we investigate the rational functions
$\K_{m,n}(q) \in \Q(q)$ defined
from the expansion
\begin{equation}\label{more}
b^m = \sum_{n=0}^m \K_{m,n}(q) \atled_n=
\sum_{n=0}^m \K_{m,n}(q^{-1}) \delta_n.
\end{equation}
One reason to be interested in these 
is that
\begin{equation}\label{interest}
(b^m,b^n)^\imath
\stackrel{\cref{eeee}}{=} 
(1, b^{m+n})^\imath
=(\atled_0, b^{m+n})^\imath
\stackrel{\cref{bilform}}{=} \K_{m+n,0}(q)
\end{equation}
for any $m,n \geq 0$.

\begin{lemma}\label{curiosity}
For $0 \leq n \leq m$, we have that
\begin{align*}
\K_{0,0}(q) &= 1,
&
\K_{m,n}(q) &= 
[n] \K_{m-1,n-1}(q) + \frac{q^n \K_{m-1,n+1}(q)}{1-q^{-2}},
\end{align*}
interpreting
$\K_{m,n}(q)$ as 0 if $n < 0$ or $n > m$.
\end{lemma}

\begin{proof}
Applying $\jmath$
to $b^m = \sum_{n=0}^m \K_{m,n}(q) \atled_n$
gives that 
$\jmath(b^m) = \sum_{n=0}^m \K_{m,n}(q) \theta^{(n)}$.
Thus, $\K_{m,n}(q)$ is the $\theta^{(n)}$-coefficient
of $\jmath(b^m)$.
Suppose that $m \geq 1$.
By \cref{jy2}(1), we have that
$\jmath(b^{m}) = \theta \jmath(b^{m-1}) + R(\jmath(b^{m-1}))$.
Then we observe using \cref{taxes} that the right hand side equals
$$
\sum_{n=1}^{m} [n] \K_{m-1,n-1}(q) \theta^{(n)}
+ \sum_{n=0}^{m-2}  
\frac{q^n \K_{m-1,n+1}(q)}{1 - q^{-2}}
\theta^{(n)}.
$$
From this, we see that the coefficient $\K_{m,n}(q)$ of $\theta^{(n)}$
in $\jmath(b^m)$ satisfies the 
recurrence relation in the statement of the lemma.
\end{proof}

We are going to give
an elementary combinatorial interpretation 
of $\K_{m,n}(q)$ 
in terms of certain chord diagrams with $n$ chords tethered to a fixed basepoint and $f=(m-n)/2$ free chords.
The notion of a chord diagram is quite standard and we will not give a formal definition, but note
in our setup that a 
pair of chords cannot intersect twice,
and a chord cannot have both ends tethered.
The following is an example of a chord diagram with $n=3$ tethered chords, $f=4$ free chords, and $c=11$ crossings:
\begin{equation}\label{orderfood}
\begin{tikzpicture}[scale=1.4,anchorbase]
\draw[ultra thick] (0,0) circle (1);
\draw (-.7,0.7) to (1,0);
\draw (-.88,-.5) to (.9,-.4);
\draw (-.45,-.9) to (.35,.95);
\draw (-.94,-.3) to (-.3,.97);
\draw (0,-1) to (.8,.6);
\draw (0,-1) to (-1,0.05);
\draw (0,-1) to (0,1);
\node at (.9,.67) {$\scriptstyle 3$};
\node at (-1.1,.09) {$\scriptstyle 1$};
\node at (0,1.12) {$\scriptstyle 2$};
\closeddot{0,-1};
\end{tikzpicture}
\end{equation}
The three tethered chords
are the ones attached to the
basepoint. 
We have also numbered the 
free endpoints of the tethered chords 
in order going clockwise around the circle.
Here is one more example 
with $n=4,f=3$ and $c=5$:
\begin{equation}\label{orderfood2}
\begin{tikzpicture}[scale=1.4,anchorbase]
\draw[ultra thick] (0,0) circle (1);
\draw (-.7,0.7) to (1,0);
\draw (.28,-.95) to (.95,.29);
\draw (0,-1) to (-.94,-0.3);
\draw (-1,-.05) to (-.3,.97);
\draw (0,-1) to (.35,.95);
\draw (0,-1) to (.8,.6);
\draw (0,-1) to (0,1);
\node at (.35,1.06) {$\scriptstyle 3$};
\node at (.9,.67) {$\scriptstyle 4$};
\node at (-1.05,-.32) {$\scriptstyle 1$};
\node at (0,1.12) {$\scriptstyle 2$};
\closeddot{0,-1};
\end{tikzpicture}
\end{equation}
In a chord diagram with $f$ free and $n$ tethered chords, the maximum possible number of crossings 
is $nf+\frac{1}{2} f(f-1)$.
Counting chord diagrams up to planar isotopy fixing the basepoint, 
let $N(f,n,c)$ be the number of chord diagrams with $f$ free chords,
$n$ tethered chords,
and $c$ crossings, and
\begin{equation}\label{high}
T_{f,n}(q) := \sum_{c=0}^{nf+\frac{1}{2}f(f-1)}
N(f,n,c) q^c \in \N[q]
\end{equation}
be the resulting generating function.
We obviously have that $T_{0,n}(q) = 1$,
and $T_{1,n}(q) = (n+1)+nq + (n-1)q^2+\cdots + q^n$. Other examples: $T_{2,0}(q) = 2+q$ and $T_{3,0}(q) = 5+6q+3q^2+q^3$.
Note also that $T_{f,n}(1) = \binom{2f+n}{n} (2f-1)!!$ (here, $n!!$ denotes the double factorial defined recursively by $n!!=n\cdot (n-2)!!$ and $0!!=
(-1)!!=1$).
Let $\{n\}$ be the classical $q$-integer
$1+q+q^2+\cdots+q^{n-1}$.

\begin{lemma}\label{swim}
The generating function $T_{f,n}(q)$ satisfies the recurrence relation
\begin{align}\label{breakfast}
T_{0,0}(q) &= 1,&
T_{f,n}(q) &= T_{f,n-1}(q) +\{n+1\} T_{f-1,n+1}(q),
\end{align}
interpreting $T_{f,n}(q)$ as 0 if $n$ or $f$ is negative.
\end{lemma}

\begin{proof}
It is clear that $T_{0,0}(q) = 1$.
Now suppose that $f+n > 0$.
Let $\mathrm{C}(f,n)$ be the set of chord diagrams
with $f$ free and $n$ tethered chords.
We are going to construct a set partition
$$
\mathrm{C}(f,n) = \overline{\mathrm{C}}(f,n)\sqcup
\coprod_{i=0}^n \mathrm{C}_i(f,n).
$$
Take a chord diagram $D \in \mathrm{C}(f,n)$.
Consider the chord $x$ in $D$
which has the nearest free endpoint to the basepoint measured in a clockwise
direction around the circumference of the circle.
There are two cases:
\begin{itemize}
\item
If $x$ is a tethered chord then
we put $D$ into the set $\overline{\mathrm{C}}(f,n)$
and let $\theta(D) \in \mathrm{C}(f,n-1)$
be the chord diagram obtained from $D$
by removing $x$.
Note that $\theta(D)$ has the same number of crossings as $D$.
An example of this situation is given by \cref{orderfood2}; for this $\theta(D)$ is
$$
\begin{tikzpicture}[scale=1.4,anchorbase]
\draw[ultra thick] (0,0) circle (1);
\draw (-.7,0.7) to (1,0);
\draw (.28,-.95) to (.95,.29);
\draw (-1,-.05) to (-.3,.97);
\draw (0,-1) to (.35,.95);
\draw (0,-1) to (.8,.6);
\draw (0,-1) to (0,1);
\node at (.35,1.06) {$\scriptstyle 2$};
\node at (.9,.67) {$\scriptstyle 3$};
\node at (0,1.12) {$\scriptstyle 1$};
\closeddot{0,-1};
\end{tikzpicture}
$$
\item
Otherwise, $x$ is a free chord.
Its furthest endpoint from the basepoint
lies between the free endpoints of the $i$th and $(i+1)$th 
tethered chords for some $0 \leq i \leq n$.
We put $D$ into the set $\mathrm{C}_i(f,n)$
and let $\theta_i(D) \in \mathrm{C}(f-1,n+1)$
be the chord diagram obtained from $D$
by replacing $x$ by a tethered chord $y$ with the same furthest endpoint as $x$.
Note that $\theta_i(D)$ has $i$ fewer crossings than $D$ since $y$ crosses $i$ fewer tethered chords compared to $x$.
An example is given by 
\cref{orderfood}; for this, 
we have that $i=2$ and
$\theta_2(D)$ is 
$$
\begin{tikzpicture}[scale=1.4,anchorbase]
\draw[ultra thick] (0,0) circle (1);
\draw (-.7,0.7) to (1,0);
\draw (-.88,-.5) to (.9,-.4);
\draw (0,-1) to (.35,.95);
\draw (-.94,-.3) to (-.3,.97);
\draw (0,-1) to (.8,.6);
\draw (0,-1) to (-1,0.05);
\draw (0,-1) to (0,1);
\node at (.35,1.06) {$\scriptstyle 3$};
\node at (.9,.67) {$\scriptstyle 4$};
\node at (-1.1,.09) {$\scriptstyle 1$};
\node at (0,1.12) {$\scriptstyle 2$};
\closeddot{0,-1};
\end{tikzpicture}
$$
\end{itemize}
We have now defined the partition of
$\mathrm{C}(f,n)$. It is also clear that 
$\theta:\overline{\mathrm{C}}(f,n)\stackrel{\sim}{\rightarrow} \mathrm{C}(f,n-1)$
and all $\theta_i:\mathrm{C}_i(f,n)
\stackrel{\sim}{\rightarrow}
\mathrm{C}(f-1,n+1)$ are bijections.
The lemma follows by
computing the generating function
$T_{f,n}(q)$ using this partition
to see that
$T_{f,n}(q) = T_{f,n-1}(q)
+ \sum_{i=0}^n q^i T_{f-1,n+1}(q)$.
\end{proof}

\begin{theorem}\label{unfair}
For $0 \leq n \leq m$ with $n \equiv m\pmod{2}$, we have that
$$
\K_{m,n}(q) =
\begin{dcases}
[n]!\frac{T_{f,n}(q^2)}{(1-q^{-2})^f}
&\text{if $m = n+2f$ for some $f \in \N$}\\
0&\text{otherwise.}
\end{dcases}
$$
\end{theorem}

\begin{proof}
It is clear from \cref{curiosity}
that $\K_{m,n}(q) = 0$ if $n \not \equiv m\pmod{2}$. 
Also using \cref{curiosity}
it follows that the rational function
$\widetilde{T}_{f,n}(q)$ defined from
$$
\widetilde{T}_{f,n}(q^2) := (1-q^{-2})^f \K_{n+2f,n}(q) / [n]!
$$
satisfies the recurrence relation
in \cref{swim}.
Hence, $\widetilde{T}_{f,n}(q^2) = 
T_{f,n}(q^2)$ and the result follows.
\end{proof}

\begin{corollary}\label{classic}
The bilinear form $(\cdot,\cdot)^\imath$
on $\Uinone$ satisfies
$$
(b^m, b^n)^\imath =
\begin{dcases}
\textstyle\frac{T_{f,0}(q^2)}{(1-q^{-2})^f}
&\text{if $m+n=2f$ for some $f\in \N$}\\
0&\text{otherwise.}
\end{dcases}
$$
\end{corollary}

\begin{proof}
This follows from the theorem using also \cref{interest}.
\end{proof}

For example, \cref{classic} implies the following:
\begin{align}  \label{BF012}
 (b,b)^\imath = 
 (1, b^2)^\imath &=\frac{1}{1-q^{-2}},
& (b^2, b^2)^\imath  =
 (b, b^3)^\imath  =
 (1, b^4)^\imath  &=
 \frac{2+q^2}{(1-q^{-2})^2}.
\end{align}
The generating function $T_{f,0}(q)$
for ordinary chord diagrams 
has been studied classically;
e.g., see \cite{chords}. 
Our more general tethered chord diagrams 
will show up again
in a slightly different guise later in the article; see \cref{reader}.

\subsection{The icanonical basis}\label{seccanbases}
So far we have not used the parameter
$t \in \{0,1\}$, but all subsequent results depend on it.
To avoid notational confusion, it is helpful to 
appeal to the construction from \cite[Chap.~4]{BW18KL} and \cite[Sec.~3.7]{BW18QSP}, which shows that $\Uinone$ 
has a modified form
$\dot\Uinone = \dot\Uinone 1_{\bar 0} \oplus \dot \Uinone 1_{\bar 1}$.
We will denote the summands here simply by
$\Uizero$ and $\Uione$ since they are actually unital algebras. In fact, 
the map $\Uinone \rightarrow \dot\Ui,
u \mapsto u 1_t$ is an algebra isomorphism.
We use this to transport
all of the results about $\Uinone$ established so far to $\Ui$, and work only with the latter from now on.
In particular, $\Ui$ is freely generated
by $b = b 1_t$, it
has the symmetries
$\rho$ and $\psi^\imath$ 
fixing $b$ as before,
it possesses a bilinear form $(\cdot,\cdot)^\imath$
as in \cref{otherform}, 
there is a linear isomorphism
$\jmath:\Ui\stackrel{\sim}{\rightarrow}\f$
as in \cref{jy2}, 
and we have the standard
basis $\delta_n\:(n \geq 0)$
for $\Ui$ satisfying \cref{stdrecurrence}.
However, one should 
have in mind that $\Ui$
is a subalgebra not of the original quantum group $\U$ but rather of
the summand of the completion of $\dot\U$
consisting of matrices
$(a_{\mu,\lambda})_{\mu,\lambda \in \Z}
\in \prod_{\lambda,\mu\in \Z} 1_\mu \dot \U 1_\lambda$
such that $a_{\mu,\lambda} = 0$ if 
$\lambda,\mu \not\equiv t \pmod{2}$.
This means that
$\Ui$ should only be allowed to act on $\U$-modules whose weights satisfy $\lambda \equiv t \pmod{2}$.
For example, the definition \cref{otherform}
of the form $(\cdot,\cdot)^\imath$
on $\Ui$ should really be written now as
\begin{equation}\label{otherformt}
( u_1, u_2)^\imath
= \lim_{\substack{\lambda\rightarrow \infty\\
\lambda\equiv t\pmod{2}}}
\big(u_1 \eta_\lambda, 
u_2 \eta_\lambda \big)_\lambda
\end{equation}
for all $u_1,u_2 \in \Ui$.

By the integrality properties from \cite[Th.~4.18]{BW18KL} and
\cite[Th.~5.3]{BW18QSP}, the symmetry $\psi_\lambda^\imath$
restricts to an anti-linear involution
on ${_\Z}V(\lambda)$.
Applying \cite[Th.~4.20]{BW18KL} and \cite[Th.~5.7]{BW18QSP},
we define the {\em icanonical basis} for $V(\lambda)$
to be the unique $\A$-basis 
$b^{(n)} \eta_\lambda\:(0 \leq n \leq \lambda)$ for ${_\Z}V(\lambda)$
such that 
each $b^{(n)} \eta_\lambda$ is $\psi_\lambda^\imath$-invariant
and 
$$
b^{(n)} \eta_\lambda -
f^{(n)}\eta_\lambda
\in\sum_{m=0}^\lambda q^{-1}\Z[q^{-1}] f^{(m)}
\eta_\lambda.
$$
As the notation suggests,
for $\lambda\equiv t\pmod{2}$,
the vector $b^{(n)} \eta_\lambda$ is obtained
by applying an element $b^{(n)} \in \Ui$
to $\eta_\lambda$. In fact, there is a
{\em unique} element $b^{(n)} \in \Ui\:(n \geq 0)$ such that $b^{(n)} \eta_\lambda$ is the icanonical basis element of $L(\lambda)$ for all $0\leq n \leq \lambda$ with $\lambda \equiv t \pmod{2}$; see \cite[Chap.~4]{BW18KL} and \cite[Th.~2.10, Th.~ 3.6]{BeW18}.
The elements $b^{(n)}\:(n \geq 0)$
thus defined give a remarkable basis for $\Ui$ again called the {\em icanonical basis}. 

Closed formulae for the icanonical basis elements were 
worked out 
in \cite{BeW18} (see also \cite{BW18KL}): 
for $n\geq 0$, we have that
\begin{equation}\label{idp}
b^{(n)} :=
\begin{dcases}
\frac{1}{[n]!}
\prod_{\substack{k=0\\k\equiv t\pmod{2}}}^{n-1}
\!\!\!\left(b^2 - [k]^2\right)
&\text{if $n$ is even}\\
\frac{b}{[n]!}
\prod_{\substack{k=1\\k\equiv t\pmod{2}}}^{n-1}
\!\!\!\left(b^2 - [k]^2\right)
&\text{if $n$ is odd.}
\end{dcases}
\end{equation}
This is also known as the {\em idivided power}.
It is straightforward to check 
from \cref{idp} that the icanonical basis satisfies the recurrence relation
\begin{align}\label{idprecurrence}
b^{(0)} &= 1,&
b b^{(n)} &= [n+1]b^{(n+1)} + \delta_{n\equiv t}
[n]b^{(n-1)},
\end{align}
for any $n \geq 0$.
\begin{theorem}\label{cornucopia}
For $n \geq 0$, we have that
\begin{align}\label{fairly}
b^{(n)} &= \sum_{m=0}^{\lfloor \frac{n}{2}\rfloor}
\frac{q^{m(2m+1-2\delta_{n\equiv t})}}{(1-q^{4}) (1-q^{8})\cdots (1-q^{4m})} \delta_{n-2m},\\\label{squarely}
\delta_n &= \sum_{m=0}^{\lfloor \frac{n}{2}\rfloor}(-1)^m
\frac{q^{m(2\delta_{n\not\equiv t}+1)}}{(1-q^{4}) (1-q^{8})\cdots (1-q^{4m})} b^{(n-2m)}.
\end{align}
\end{theorem}

\begin{proof}
To prove the first formula, use \cref{stdrecurrence} to verify that the expression on the right hand side satisfies the recurrence relation \cref{idprecurrence}.
Similarly, \cref{squarely} follows by using \cref{idprecurrence} to verify that the expression on the right hand side satisfy
the recurrence relation \cref{stdrecurrence}.
\end{proof}

\begin{corollary}\label{dontreallyneedthedetailedformula}
The icanonical basis of $\Ui$ is almost orthonormal in the sense that
$$
(b^{(m)}, b^{(n)})^\imath \in \delta_{m,n}+ \big(q^{-1} \Z \llbracket q^{-1}\rrbracket\cap \Q(q)\big)
$$
for $m,n \geq 0$.
\end{corollary}

\begin{proof}
This is clear from \cref{fairly} and \cref{bilform}.
\end{proof}

\begin{remark}
Using \cref{fairly}, one can derive the following explicit formula for the pairings between icanonical basis elements:
\begin{align*}
(b^{(n)},b^{(m)})^\imath &=
\sum_{\substack{0 \leq i \leq \min(m,n)\\
i\equiv n \equiv m\pmod{2}}}
\frac{q^{-\frac{1}{2}(n-i)
(n-i+1-2\delta_{n\equiv t})-\frac{1}{2}(m-i)
(m-i+1-2\delta_{m\equiv t})}}{
\prod_{j=1}^i
(1-q^{-2j})
\prod_{k=1}^{\frac{n-i}{2}}(1-q^{-4k})
\prod_{l=1}^{\frac{m-i}{2}}(1-q^{-4l})}
\end{align*}
for any $m,n \geq 0$.
This is 0 if $m \not\equiv n\pmod{2}$.
\end{remark}

The icanonical basis in fact 
gives a basis
for an {integral form} $\UAi$ of $\Ui$ over $\A$. Equivalently,
we have that
$$
\UAi = \left\{u \in \Ui\:\big|\:u \left({_\Z}V(\lambda)\right)
\subseteq {_\Z}V(\lambda)\text{ for all }\lambda \in \N\text{ with }\lambda\equiv t\pmod{2}\right\},
$$
from which one sees that $\UAi$ is a $\A$-subalgebra of $\Ui$.
Since both $\rho$ and $\psi^\imath$ 
fix each of the icanonical basis elements
$b^{(n)}$, they restrict to symmetries
on $\UAi$.
Also, the form on $\Ui$ restricts
to $(\cdot,\cdot)^\imath:\UAi\times\UAi\rightarrow \Z[q,q^{-1}]$.
From \cref{stdrecurrence}, it is apparent that $\delta_n \in \Ui$ does not lie in the
integral form. 

\subsection{The character ring}\label{charring}
The {\em character ring} 
is the ring 
$\Q(q) \llbracket \chi \rrbracket$
for a formal variable $\chi$.
It is natural to consider from a representation-theoretic perspective (see \cref{grch}).
We view $\Q(q)\llbracket \chi \rrbracket$ as a left $\Ui$-module so that
\begin{align}\label{dualaction}
b \chi^n &= \begin{cases}
\chi^{n-1}&\text{if $n > 0$}\\
0&\text{if $n=0$,}
\end{cases}
\end{align}
and extending in the natural way to power series.
We identify the character ring 
with the full linear dual $\left(\Ui\right)^*$ so that
$\sum_{n \geq 0} f_n(q) \chi^n$ is the function mapping
$b^n$ to $f_n(q)$.
Thus, the topological basis for $\Q(q)\llbracket \chi \rrbracket$ given by the monomials $\chi^n\:(n \geq 0)$ is dual to the basis $b^n\:(n \geq 0)$ of $\Ui$.
Then we let $\el_n, \bar\delta_n$ and $\bar\atled_n$
be the unique elements of the character ring
such that
\begin{align}\label{chmap}
\el_n\big(b^{(m)}\big) &= 
\bar\delta_n\big(\delta_m\big) = \bar\atled_n\big(\atled_m\big)
=\delta_{m,n}
\end{align}
for $m,n \geq 0$.
The topological bases $\el_n\:(n \geq 0)$,
$\bar\delta_n\:(n \geq 0)$ and $\bar\atled_n\:(n \geq 0)$
for $\Q(q) \llbracket \chi\rrbracket$ give the
{\em dual canonical basis}, {\em proper standard basis} and
{\em proper costandard basis}, respectively.

There is a {\em bar involution}
on the character ring, which is the
anti-linear map
\begin{align}\label{barinv}
\bar\psi^\imath : \Q(q)\llbracket\chi\rrbracket
&\rightarrow
\Q(q)\llbracket\chi\rrbracket,
&
\sum_{n \geq 0} f_n(q) \chi^n &\mapsto \sum_{n \geq 0} 
f_n(q^{-1}) \chi^n.
\end{align}
This is compatible with the bar involution on $\Ui$ in the sense
that $\bar\psi^\imath(ug) = \psi^\imath(u) \bar\psi^\imath(g)$
for $u \in \Ui$ and $g \in \Q(q)\llbracket \chi\rrbracket$.
Also, the bar involution on $\Q(q)\llbracket \chi\rrbracket$ is related to the bar involution on $\Ui$ by the formula
\begin{equation}\label{dualbar}
\big(\bar\psi^\imath(g)\big)(u) = \overline{g(\psi^\imath(u))}.
\end{equation}
This follows easily from \cref{barinv} as $\psi^\imath(b^n) = b^n$ for each $n \geq 0$.
Using \cref{dualbar} and the definitions, it 
follows that the dual canonical basis elements $\el_n$ are 
fixed by $\bar\psi^\imath$, and
$\bar\psi^\imath(\bar\delta_n) = \bar\atled_n$.

From \cref{cornucopia}, we get that
\begin{align}\label{frombefore}
\bar\delta_n &= \sum_{m=0}^\infty
\frac{q^{m(2m+1-2\delta_{n\equiv t})}}{(1-q^{4}) (1-q^{8})\cdots (1-q^{4m})} \el_{n+2m},\\
\el_n &= \sum_{m=0}^\infty(-1)^m
\frac{q^{m(2\delta_{n\not\equiv t}+1)}}{(1-q^{4}) (1-q^{8})\cdots (1-q^{4m})} \bar\delta_{n+2m}.
\end{align}
for $n \geq 0$.
Also the following recurrence relations
may be deduced from \cref{stdrecurrence,idprecurrence}:
\begin{align}\label{witcher2}
b \bar \delta_{n} &= [n]\bar\delta_{n-1} + \frac{q^{-n}}{1-q^{2}}\bar\delta_{n+1},\\
b \el_n &= [n] \el_{n-1} +  \delta_{n\not\equiv t} [n+1] \el_{n+1}
\label{witcher1}
\end{align}
for any $n \geq 0$.

We proceed to derive explicit formulae for
$\bar\delta_n$ and $\el_n$ as formal series in $\chi$.

\begin{lemma}\label{reaching}
For $n \geq 0$, 
we have that
$$
\bar\delta_n = 
[n]!  \sum_{f \geq 0} 
\frac{T_{f,n}(q^{-2})}{(1-q^{2})^f} \chi^{n+2f}.
$$
\end{lemma}

\begin{proof}
By \cref{more}, we have that
$b^m = \sum_{n=0}^m \K_{m,n}(q^{-1}) \delta_n$.
Applying the function $\bar\delta_n$ to this equation,
we deduce that the coefficient of $\chi^m$ in the expansion of $\bar\delta_n$
is equal to $\K_{m,n}(q^{-1})$.
It remains to apply \cref{unfair}.
\end{proof}

\begin{lemma}\label{americano}
$\el_0 = 
\begin{dcases}
1&\text{if $t=0$}\\
1+\chi^2+\chi^4+\chi^6+\cdots&\text{if $t=1$.}
\end{dcases}$
\end{lemma}

\begin{proof}
Suppose first that $t=0$.
We need to show that
$\el_0(b^n) = \delta_{n,0}$ for any 
$n \geq 0$.
This is clear for $n=0$ since $b^0 = b^{(0)}$
and $\el_0(b^{(0)}) = 1$. 
Also \cref{idp} shows that all $b^{(n)}\:(n > 0)$
are divisible by $b$, so we can
use \cref{idp} to express $b^n\:(n > 0)$ as a linear combination of $b^{(1)},\dots,b^{(n)}$. This implies
that $\el_0(b^n) = 0$ for $n > 0$ as required.

Now suppose that $t=1$.
We need to show that
$\el_0(b^{2n+1}) = 0$
and $\el_0(b^{2n}) = 1$
for $n \geq 0$.
By \cref{idp},
$b^{(2n+1)}$ is a linear combination of
$b^{2m+1}$ for $0 \leq m \leq n$,
and inverting obviously gives that
$b^{2n+1}$ is a linear combination of $b^{(2m+1)}$ for $0 \leq m \leq n$.
This implies that $\el_0(b^{2n+1}) = 0$.
Also, by \cref{idp} again, $b^{(0)} = 1$ and
$[2n][2n-1] b^{(2n)} = (b^2 - [2n-1]^2) b^{(2n-2)}$
for $n \geq 1$.
Using this, one shows by induction on $n\geq 0$ that
$b^{2n} = a_n b^{(2n)} + \cdots + a_1 b^{(2)}
+ b^{(0)}$ for some $a_1,\dots,a_n \in \Q(q)$.
It follows that $\el_0(b^{2n}) = 1$.
\end{proof}

\begin{theorem}\label{oops}
We have that 
\begin{equation}\label{oopsformula}
\el_n = [n]! \chi^n \prod_{\substack{1 \leq k\leq n+1 \\k\equiv t \pmod{2}}}
\frac{1}{1-[k]^2 \chi^2}
=
[n]! \sum_{m \geq 0}
\left(\sum_{\alpha \in \Par_t(m\times n)}
[\alpha_1+1]^2\cdots [\alpha_m+1]^2\right) \chi^{n+2m}
\end{equation}
where $\Par_t(m\times n)$
is the set of 
$\alpha \in \N^m$ with
$0 \leq \alpha_1 \leq \cdots \leq \alpha_m \leq n$ and $\alpha_i\not\equiv t \pmod{2}$ for each $i$.
\end{theorem}

\begin{proof}
The second equality follows by expanding the product. To prove the first equality, we proceed by induction on $n$.
The induction base follows from \cref{americano}.
For the induction step, take $n > 0$.
The constant term of $\el_n$ is 0
since $\el_n(1) = 
\el_n(b^{(0)}) = 0$
so we have that $b \el_n = \el_n / \chi$ by \cref{dualaction}.
Suppose first that $n\equiv t\pmod{2}$.
Then \cref{witcher1} shows that 
\begin{equation}\label{oops1}
\el_n = [n]\ \chi \el_{n-1}
\end{equation}
and we easily get done by induction in this case.
When $n\not\equiv t\pmod{2}$, \cref{witcher1} gives that 
$$
\el_n = [n]\ \chi \el_{n-1} + [n+1]\ \chi \el_{n+1}
= [n]\ \chi \el_{n-1} + [n+1]^2 \chi^2 \el_n.
$$
Hence, \begin{equation}\label{oops2}
\el_n = \frac{[n]\ \chi}{1-[n+1]^2 \chi^2} \el_{n-1},
\end{equation}
and again the result follows by induction.
\end{proof}

\begin{corollary}\label{texas}
For $n \geq 0$, we have that 
$$
b^n= \sum_{m=0}^{\lfloor \frac{n}{2}\rfloor}[n-2m]!
\left(\sum_{\alpha \in \Par_t(m\times (n-2m))}
[\alpha_1+1]^2\cdots [\alpha_m+1]^2\right) b^{(n-2m)}.
$$
\end{corollary}

\begin{proof}
The coefficient of $b^{(n-2m)}$
in the expansion of $b^n$ is
$\el_{n-2m}(b^n)$, i.e., it is the $\chi^n$-coefficient of 
$\el_{n-2m}$.
Now use \cref{oops}.
\end{proof}

%% file: s3-nilbrauer.tex
\setcounter{section}{2}

\section{The nil-Brauer category}\label{sect3}

{\em For the remainder of the article}, we will work over a field $\kk$ of characteristic different from 2.
All algebras, categories, functors, etc. will be assumed to be $\kk$-linear without further mention, and we reserve the symbol $\otimes$
for tensor products of vector spaces or algebras over $\kk$.
By a {\em graded} category,
{\em graded} monoidal category,
{\em graded} functor, etc.
we mean one that is enriched in the closed symmetric monoidal category $\gVec$ of graded vector spaces.

In this section, we first recall the definition of the nil-Brauer category $\cNB_t$ and the crucial basis theorem for its morphism spaces
from \cite{BWWbasis}. Then we relate the graded dimensions of these spaces to the bilinear form $(\cdot,\cdot)^\imath$ 
on the iquantum group $\Ui$. Finally, we discuss the center of $\cNB_t$, and prove a useful result about minimal polynomials.

\subsection{Definition and basic properties}\label{dbp}
We use the usual string calculus for morphisms in strict monoidal categories; our general convention is that $f \circ g$ denotes composition 
of $f$ drawn on top of $g$ 
(``vertical composition'')
and $f \star g$ denotes the 
tensor product of $f$ drawn to the left of $g$ (``horizontal composition''). We always draw string diagrams so that the underlying strings are smooth curves. Recall the following definition from \cite[Def.~2.1]{BWWbasis}. 

\begin{definition}\label{NBdef}
The {\em nil-Brauer category} $\cNB_t$ 
is the strict graded monoidal category
with one generating object $B$ (whose identity endomorphism will be
represented diagrammatically by the unlabeled string
$\;\begin{tikzpicture}[anchorbase]\draw[-] (0,-.2) to (0,.2);
\end{tikzpicture}\;$) 
and four generating morphisms
\begin{align}\label{gens}
\begin{tikzpicture}[anchorbase]
	\draw[-] (0.08,-.3) to (0.08,.4);
    \closeddot{0.08,.05};
\end{tikzpicture}
&:B\rightarrow B,&
\begin{tikzpicture}[anchorbase]
\draw[-] (0.28,-.3) to (-0.28,.4);
	\draw[-] (-0.28,-.3) to (0.28,.4);
\end{tikzpicture}
&:B\star B \rightarrow B \star B,
&
\begin{tikzpicture}[anchorbase]
	\draw[-] (0.4,0) to[out=90, in=0] (0.1,0.4);
	\draw[-] (0.1,0.4) to[out = 180, in = 90] (-0.2,0);
\end{tikzpicture}
&:B \star B\rightarrow\one,&
\begin{tikzpicture}[anchorbase]
	\draw[-] (0.4,0.4) to[out=-90, in=0] (0.1,0);
	\draw[-] (0.1,0) to[out = 180, in = -90] (-0.2,0.4);
\end{tikzpicture}
&:\one\rightarrow B\star B,
\\\notag
&\text{(degree $2$)}&
&\text{(degree $-2$)}&
&\text{(degree $0$)}&
&\text{(degree $0$)}
\end{align}
subject to the following relations:
\begin{align}\label{rels1}
\begin{tikzpicture}[anchorbase]
	\draw[-] (0.28,0) to[out=90,in=-90] (-0.28,.6);
	\draw[-] (-0.28,0) to[out=90,in=-90] (0.28,.6);
	\draw[-] (0.28,-.6) to[out=90,in=-90] (-0.28,0);
	\draw[-] (-0.28,-.6) to[out=90,in=-90] (0.28,0);
\end{tikzpicture}
&=0,
&\begin{tikzpicture}[anchorbase]
	\draw[-] (0.45,.6) to (-0.45,-.6);
	\draw[-] (0.45,-.6) to (-0.45,.6);
        \draw[-] (0,-.6) to[out=90,in=-90] (-.45,0);
        \draw[-] (-0.45,0) to[out=90,in=-90] (0,0.6);
\end{tikzpicture}
&=
\begin{tikzpicture}[anchorbase]
	\draw[-] (0.45,.6) to (-0.45,-.6);
	\draw[-] (0.45,-.6) to (-0.45,.6);
        \draw[-] (0,-.6) to[out=90,in=-90] (.45,0);
        \draw[-] (0.45,0) to[out=90,in=-90] (0,0.6);
\end{tikzpicture}\:,\\
\label{rels2}
\begin{tikzpicture}[baseline=-2.5mm]
\draw (0,-.15) circle (.3);
\end{tikzpicture}
&= t 1_\one,
&\begin{tikzpicture}[anchorbase]
  \draw[-] (0.3,0) to (0.3,.4);
	\draw[-] (0.3,0) to[out=-90, in=0] (0.1,-0.4);
	\draw[-] (0.1,-0.4) to[out = 180, in = -90] (-0.1,0);
	\draw[-] (-0.1,0) to[out=90, in=0] (-0.3,0.4);
	\draw[-] (-0.3,0.4) to[out = 180, in =90] (-0.5,0);
  \draw[-] (-0.5,0) to (-0.5,-.4);
\end{tikzpicture}
&=
\begin{tikzpicture}[anchorbase]
  \draw[-] (0,-0.4) to (0,.4);
\end{tikzpicture}
=\begin{tikzpicture}[anchorbase]
  \draw[-] (0.3,0) to (0.3,-.4);
	\draw[-] (0.3,0) to[out=90, in=0] (0.1,0.4);
	\draw[-] (0.1,0.4) to[out = 180, in = 90] (-0.1,0);
	\draw[-] (-0.1,0) to[out=-90, in=0] (-0.3,-0.4);
	\draw[-] (-0.3,-0.4) to[out = 180, in =-90] (-0.5,0);
  \draw[-] (-0.5,0) to (-0.5,.4);
\end{tikzpicture}
\:,\\\label{rels3}
\begin{tikzpicture}[anchorbase,scale=1.1]
	\draw[-] (0.35,.3)  to [out=90,in=-90] (-.1,.9) to [out=90,in=180] (.1,1.1);
 \draw[-] (-.15,.3)  to [out=90,in=-90] (.3,.9) to [out=90,in=0] (.1,1.1);
\end{tikzpicture}
&=0
\:,&
\begin{tikzpicture}[anchorbase,scale=1.1]
	\draw[-] (0.5,0) to[out=90, in=0] (0.1,0.5);
	\draw[-] (0.1,0.5) to[out = 180, in = 90] (-0.3,0);
 \draw[-] (0.1,0) to[out=90,in=-90] (-.3,.7);
\end{tikzpicture}
&=
\begin{tikzpicture}[anchorbase,scale=1.1]
	\draw[-] (0.5,0) to[out=90, in=0] (0.1,0.5);
	\draw[-] (0.1,0.5) to[out = 180, in = 90] (-0.3,0);
 \draw[-] (0.1,0) to[out=90,in=-90] (.5,.7);
\end{tikzpicture}\:,
\\
\label{rels4}
\begin{tikzpicture}[anchorbase,scale=1.4]
	\draw[-] (0.25,.3) to (-0.25,-.3);
	\draw[-] (0.25,-.3) to (-0.25,.3);
 \closeddot{-0.12,0.145};
\end{tikzpicture}
-
\begin{tikzpicture}[anchorbase,scale=1.4]
	\draw[-] (0.25,.3) to (-0.25,-.3);
	\draw[-] (0.25,-.3) to (-0.25,.3);
     \closeddot{.12,-0.145};
\end{tikzpicture}
&=
\begin{tikzpicture}[anchorbase,scale=1.4]
 	\draw[-] (0,-.3) to (0,.3);
	\draw[-] (-.3,-.3) to (-0.3,.3);
\end{tikzpicture}
-
\begin{tikzpicture}[anchorbase,scale=1.4]
 	\draw[-] (-0.15,-.3) to[out=90,in=180] (0,-.1) to[out=0,in=90] (0.15,-.3);
 	\draw[-] (-0.15,.3) to[out=-90,in=180] (0,.1) to[out=0,in=-90] (0.15,.3);
\end{tikzpicture}
\:,&
\begin{tikzpicture}[anchorbase,scale=1.1]
	\draw[-] (0.4,0) to[out=90, in=0] (0.1,0.5);
	\draw[-] (0.1,0.5) to[out = 180, in = 90] (-0.2,0);
 \closeddot{.38,.2};
\end{tikzpicture}
&=
-\begin{tikzpicture}[anchorbase,scale=1.1]
	\draw[-] (0.4,0) to[out=90, in=0] (0.1,0.5);
	\draw[-] (0.1,0.5) to[out = 180, in = 90] (-0.2,0);
 \closeddot{-.18,.2};
\end{tikzpicture}
\:.
\end{align}
\end{definition}

\begin{remark}
One source of motivation for 
\cref{NBdef} is the expected compatibility
of $\cNB_t$
with the bilinear form $(\cdot,\cdot)^\imath$ on $\Ui$, something which will be proved in general in \cref{toinfinityandbeyond}. 
From this perspective, the
formulae \cref{BF012}
suggest the existence of generators of the degrees specified in \cref{gens} and some of the basic relations. This is similar to Lauda's approach to categorification of
${\mathrm U}_q(\mathfrak{sl}_2)$ in 
\cite{Lauda}. 
\end{remark}

The following relations are easily derived from the defining relations in \cite[(2.6)--(2.8)]{BWWbasis}:
\begin{align}\label{rels5}
\begin{tikzpicture}[anchorbase,scale=1.1]
	\draw[-] (0.5,0) to[out=-90, in=0] (0.1,-0.5);
	\draw[-] (0.1,-0.5) to[out = 180, in = -90] (-0.3,0);
 \draw[-] (0.1,0) to[out=-90,in=90] (-.3,-.7);
\end{tikzpicture}
&=
\begin{tikzpicture}[anchorbase,scale=1.1]
	\draw[-] (0.5,0) to[out=-90, in=0] (0.1,-0.5);
	\draw[-] (0.1,-0.5) to[out = 180, in = -90] (-0.3,0);
 \draw[-] (0.1,0) to[out=-90,in=90] (.5,-.7);
\end{tikzpicture}\:,&
\begin{tikzpicture}[anchorbase]
	\draw[-] (0,-.3)  to (0,-.1) to [out=90,in=180] (.3,.3) to [out=0,in=90] (.5,.1) to[out=-90,in=0] (.3,-.1) to [out=180,in=-90] (0,.3) to (0,.5);
\end{tikzpicture}&=
0=\begin{tikzpicture}[anchorbase]
	\draw[-] (0,-.3)  to (0,-.1) to [out=90,in=0] (-.3,.3) to [out=180,in=90] (-.5,.1) to[out=-90,in=180] (-.3,-.1) to [out=0,in=-90] (0,.3) to (0,.5);
\end{tikzpicture}\:,\\
\label{rels6}
\begin{tikzpicture}[anchorbase,scale=1.1]
	\draw[-] (0.35,-.3)  to [out=-90,in=90] (-.1,-.9) to [out=-90,in=180] (.1,-1.1);
 \draw[-] (-.15,-.3)  to [out=-90,in=90] (.3,-.9) to [out=-90,in=0] (.1,-1.1);
\end{tikzpicture}
&=0,&
\begin{tikzpicture}[anchorbase,scale=1.4]
 	\draw[-] (-0.25,-.3) to[out=90,in=180] (0,.1) to[out=0,in=90] (0.25,-.3);
 	\draw[-] (-0.25,.3) to[out=-90,in=180] (0,-.1) to[out=0,in=-90] (0.25,.3);
\end{tikzpicture}&=0,
\\
\label{rels7}
\begin{tikzpicture}[anchorbase,scale=1.4]
	\draw[-] (0.25,.3) to (-0.25,-.3);
	\draw[-] (0.25,-.3) to (-0.25,.3);
 \closeddot{-0.12,-0.145};
\end{tikzpicture}
-
\begin{tikzpicture}[anchorbase,scale=1.4]
	\draw[-] (0.25,.3) to (-0.25,-.3);
	\draw[-] (0.25,-.3) to (-0.25,.3);
     \closeddot{.12,0.145};
\end{tikzpicture}
&=
\begin{tikzpicture}[anchorbase,scale=1.4]
 	\draw[-] (0,-.3) to (0,.3);
	\draw[-] (-.3,-.3) to (-0.3,.3);
\end{tikzpicture}
-
\begin{tikzpicture}[anchorbase,scale=1.4]
 	\draw[-] (-0.15,-.3) to[out=90,in=180] (0,-.1) to[out=0,in=90] (0.15,-.3);
 	\draw[-] (-0.15,.3) to[out=-90,in=180] (0,.1) to[out=0,in=-90] (0.15,.3);
\end{tikzpicture}
\:,&
\begin{tikzpicture}[anchorbase,scale=1.1]
	\draw[-] (0.4,0) to[out=-90, in=0] (0.1,-0.5);
	\draw[-] (0.1,-0.5) to[out = 180, in = -90] (-0.2,0);
 \closeddot{.38,-.2};
\end{tikzpicture}
&=
-\begin{tikzpicture}[anchorbase,scale=1.1]
	\draw[-] (0.4,0) to[out=-90, in=0] (0.1,-0.5);
	\draw[-] (0.1,-0.5) to[out = 180, in = -90] (-0.2,0);
 \closeddot{-.18,-.2};
\end{tikzpicture}
\:.
\end{align}
In view of the last relation from \cref{rels3}
and the first relation from \cref{rels5},
we can unambiguously 
denote the morphisms in these two equations by the ``pitchforks''
$\begin{tikzpicture}[anchorbase,scale=.6]
	\draw[-] (0.5,0) to[out=90, in=0] (0.1,0.5);
	\draw[-] (0.1,0.5) to[out = 180, in = 90] (-0.3,0);
 \draw[-] (0.1,0) to[out=90,in=-90] (.1,.8);
\end{tikzpicture}
$
and 
$\begin{tikzpicture}[anchorbase,scale=.6]
	\draw[-] (0.5,0) to[out=-90, in=0] (0.1,-0.5);
	\draw[-] (0.1,-0.5) to[out = 180, in = -90] (-0.3,0);
 \draw[-] (0.1,0) to[out=-90,in=90] (.1,-.8);
\end{tikzpicture}
$, respectively.
Together with the last relation of \cref{rels2}, it then follows that a string diagram with no dots can be deformed under planar isotopy without changing the morphism that it represents. This is not true in the presence of dots due to the sign in the last relations of \cref{rels4,rels7}---there is a sign change whenever a dot slides across the critical point of a cup or cap.

The relations discussed so far imply that there are strict graded monoidal functors
\begin{align}
\tR:\cNB_t &\rightarrow \cNB_t^{\rev}, &B &\mapsto B, & s &\mapsto (-1)^{\ndots(s)} s^{\leftrightarrow},
\label{R}\\\label{T}
\tT:\cNB_t &\rightarrow \cNB_t^{\op}, &B &\mapsto B, & s &\mapsto s^{\updownarrow}.
\end{align}
Here, for a string diagram $s$ 
we use $s^{\updownarrow}$ 
and $s^{\leftrightarrow}$ to denote
its reflection in a horizontal or vertical axis,
and $\ndots(s)$ denotes 
the total number of dots 
in the diagram.
The category $\cNB_t$ is pivotal with duality functor $\tD := \tR \circ \tT = \tT \circ \tR$, which rotates a string diagram $s$ through $180^\circ$ then scales by $(-1)^{\ndots(s)}$. 

\subsection{Generating functions for dots and bubbles}
Next we recall the generating function formalism from \cite[Sec.~2]{BWWbasis}. 
We denote the $r$th power of 
$\txtdot$ under vertical composition
simply by labeling the dot with $r$. More generally,
given a polynomial $f(x) = \sum_{r\geq 0} c_r x^r \in \kk[x]$
and a dot in some string diagram $s$,
we denote
$$
\sum_{r \geq 0} c_r \times 
(\text{the morphism obtained from $s$ by labeling the dot by }r)
$$
by attaching what we call a {\em pin} to the dot, labeling 
the node at the head of the pin by $f(x)$:
\begin{equation}\label{labelledpin}
\begin{tikzpicture}[anchorbase]
	\draw[-] (0,-.4) to (0,.4);
    \pinX{(0,0)}{}{(.8,0)}{f(x)};
\end{tikzpicture}\;
:=
\sum_{r \geq 0} c_r 
\:\begin{tikzpicture}[anchorbase]
	\draw[-] (0,-.4) to (0,.4);
    \closeddot{0,0};
    \node at (.2,0) {$\scriptstyle{r}$};
\end{tikzpicture}
\in \End_{\cNB_t}(B).
\end{equation}
In the drawing of a pin, the arm and the head of the pin can be moved freely around larger diagrams so long as the dot at the pointy end stays put---these are not part of the string calculus. More generally, $f(x)$ here could be a polynomial with coefficients in the algebra $\kk\lround u^{-1} \rround$ of formal Laurent series in an indeterminate $u^{-1}$; then the string $s$ decorated with pin labeled $f(x)$ defines a generating function of morphisms. 

We will use the following shorthands for the
generating functions of \cite[(2.14)--(2.15)]{BWWbasis}:
\begin{align}\label{ping}
\begin{tikzpicture}[anchorbase]
	\draw[-] (0,-.4) to (0,.4);
\circled{0,0}{u};
\end{tikzpicture}\;
&:=
\begin{tikzpicture}[anchorbase]
	\draw[-] (0,-.4) to (0,.4);
    \closeddot{0,0};
    \pinX{(0,0)}{}{(1.3,0)}{(u-x)^{-1}};
\end{tikzpicture}\;
= 
u^{-1}
\: \begin{tikzpicture}[anchorbase]
	\draw[-] (0.08,-.4) to (0.08,.4);
\end{tikzpicture}
+
u^{-2}\: \begin{tikzpicture}[anchorbase]
	\draw[-] (0.08,-.4) to (0.08,.4);
    \closeddot{0.08,0};
\end{tikzpicture}
+
u^{-3}
\: \begin{tikzpicture}[anchorbase]
	\draw[-] (0.08,-.4) to (0.08,.4);
    \closeddot{0.08,.1};
    \closeddot{0.08,-.1};
\end{tikzpicture}
+
u^{-4}
\: \begin{tikzpicture}[anchorbase]
	\draw[-] (0.08,-.4) to (0.08,.4);
    \closeddot{0.08,.2};
    \closeddot{0.08,0};
    \closeddot{0.08,-.2};
\end{tikzpicture}
+\cdots \in \End_{\cNB_t}(B)\llbracket u^{-1}\rrbracket\:,\\
\label{pong}
\begin{tikzpicture}[anchorbase]
	\draw[-] (0,-.4) to (0,.4);
\circledbar{0,0}{u};
\end{tikzpicture}\;
&:=
\begin{tikzpicture}[anchorbase]
	\draw (0,-.4) to (0,.4);
    \closeddot{0,0};
    \pinX{(0,0)}{}{(1.3,0)}{(u+x)^{-1}};
\end{tikzpicture}\;
= 
u^{-1}\: \begin{tikzpicture}[anchorbase]
	\draw[-] (0.08,-.4) to (0.08,.4);
\end{tikzpicture}
-
u^{-2}\: \begin{tikzpicture}[anchorbase]
	\draw[-] (0.08,-.4) to (0.08,.4);
    \closeddot{0.08,0};
\end{tikzpicture}
+
u^{-3}
\: \begin{tikzpicture}[anchorbase]
	\draw[-] (0.08,-.4) to (0.08,.4);
    \closeddot{0.08,.1};
    \closeddot{0.08,-.1};
\end{tikzpicture}
-
u^{-4}\: \begin{tikzpicture}[anchorbase]
	\draw[-] (0.08,-.4) to (0.08,.4);
    \closeddot{0.08,.2};
    \closeddot{0.08,0};
    \closeddot{0.08,-.2};
\end{tikzpicture}
+\cdots \in \End_{\cNB_t}(B)\llbracket u^{-1}\rrbracket\:.
\end{align}
The notation here is motivated by the following standard trick:
for any $f(x) \in \kk[x]$, we have that
\begin{align}\label{trick}
\left[f(u)\:\begin{tikzpicture}[anchorbase]
	\draw[-] (0.08,-.3) to (0.08,.4);
	\circled{0.08,.05}{u};
\end{tikzpicture}
 \right]_{u^{-1}}
&=
\begin{tikzpicture}[anchorbase]
	\draw[-] (0.08,-.3) to (0.08,.4);
    \closeddot{0.08,.05};
	\pinX{(0.08,.05)}{}{(1,.05)}{f(x)};
\end{tikzpicture}\:,
&
\left[f(u)\:\begin{tikzpicture}[anchorbase]
	\draw[-] (0.08,-.3) to (0.08,.4);
	\circledbar{0.08,.05}{u};
\end{tikzpicture}
 \right]_{u^{-1}}
&=
\begin{tikzpicture}[anchorbase]
	\draw[-] (0.08,-.3) to (0.08,.4);
    \closeddot{0.08,.05};
	\pinX{(0.08,.05)}{}{(1,.05)}{f(-x)};
\end{tikzpicture}\:,
\end{align}
where
$[-]_{u^r}$ denotes the $u^r$-coefficient of the formal Laurent series inside the brackets.
These identities follow by using linearity to reduce to the case that $f(x) = x^n$ for $n \geq 0$, then explicitly computing coefficients on both sides.
As we do with ordinary dots, we denote the $n$th power of one of these 
``dot generating functions'' 
by labeling them also by $n$.
This makes sense for any $n \in \Z$
since we have by the definitions that
\begin{align*}
\begin{tikzpicture}[anchorbase]
	\draw[-] (0,-.4) to (0,.4);
    \circled{0,0}{u};
    \node at (.35,0) {$\scriptstyle -1$};
\end{tikzpicture}
&:=
\left(\;\begin{tikzpicture}[anchorbase]
	\draw[-] (0,-.4) to (0,.4);
    \circled{0,0}{u};
\end{tikzpicture}\;\right)^{-1}
=
\begin{tikzpicture}[anchorbase]
	\draw[-] (0,-.4) to (0,.4);
    \closeddot{0,0};
    \pinX{(0,0)}{}{(1,0)}{u-x};
\end{tikzpicture}\;
= 
u
\: \begin{tikzpicture}[anchorbase]
	\draw[-] (0.08,-.4) to (0.08,.4);
\end{tikzpicture}\;
-
\begin{tikzpicture}[anchorbase]
	\draw[-] (0.08,-.4) to (0.08,.4);
    \closeddot{0.08,0};
\end{tikzpicture}\;,&
\begin{tikzpicture}[anchorbase]
	\draw[-] (0,-.4) to (0,.4);
    \circledbar{0,0}{u};
    \node at (.35,0) {$\scriptstyle -1$};
\end{tikzpicture}&:=
\left(\;\begin{tikzpicture}[anchorbase]
	\draw[-] (0,-.4) to (0,.4);
    \circledbar{0,0}{u};
\end{tikzpicture}\;\right)^{-1}
=
\begin{tikzpicture}[anchorbase]
	\draw (0,-.4) to (0,.4);
    \closeddot{0,0};
    \pinX{(0,0)}{}{(1,0)}{u+x};
\end{tikzpicture}\;
= 
u\: \begin{tikzpicture}[anchorbase]
	\draw[-] (0.08,-.4) to (0.08,.4);
\end{tikzpicture}\;
+
\begin{tikzpicture}[anchorbase]
	\draw[-] (0.08,-.4) to (0.08,.4);
    \closeddot{0.08,0};
\end{tikzpicture}\;.
\end{align*}
The endomorphisms 
\cref{ping,pong} obviously commute with each 
other and all other pins.
Note also that $\tT$ and $\tR$ satisfy
\begin{align}\label{RTdots}
\tR\left(
\begin{tikzpicture}[anchorbase]
	\draw[-] (0.08,-.3) to (0.08,.4);
\circledbar{0.08,0.05}{u};
\end{tikzpicture}\right) &=
\begin{tikzpicture}[anchorbase]
	\draw[-] (0.08,-.3) to (0.08,.4);
\circled{0.08,0.05}{u};
\end{tikzpicture}\:,&
\tR\left(
\begin{tikzpicture}[anchorbase]
	\draw[-] (0.08,-.3) to (0.08,.4);
\circled{0.08,0.05}{u};
\end{tikzpicture}\right) &=
\begin{tikzpicture}[anchorbase]
	\draw[-] (0.08,-.3) to (0.08,.4);
\circledbar{0.08,0.05}{u};
\end{tikzpicture}\:,&
\tT\left(
\begin{tikzpicture}[anchorbase]
	\draw[-] (0.08,-.3) to (0.08,.4);
\circledbar{0.08,0.05}{u};
\end{tikzpicture}\right) &=
\begin{tikzpicture}[anchorbase]
	\draw[-] (0.08,-.3) to (0.08,.4);
\circledbar{0.08,0.05}{u};
\end{tikzpicture}\:,&
\tT\left(
\begin{tikzpicture}[anchorbase]
	\draw[-] (0.08,-.3) to (0.08,.4);
\circled{0.08,0.05}{u};
\end{tikzpicture}\right) &=
\begin{tikzpicture}[anchorbase]
	\draw[-] (0.08,-.3) to (0.08,.4);
\circled{0.08,0.05}{u};
\end{tikzpicture}\:.
\end{align}
 Another useful trick is to apply the substitution
 $u \mapsto -u$;
 this interchanges $\begin{tikzpicture}[anchorbase,scale=.8]
\draw[-] (-.2,.3) to (-.2,-.3);
\circled{-.2,0}{u};
\end{tikzpicture}$ and 
$-\begin{tikzpicture}[anchorbase,scale=.8]
\draw[-] (-.2,.3) to (-.2,-.3);
\circledbar{-.2,0}{u};
\end{tikzpicture}$.

It is clear from the last relation in \cref{rels3} that
$\begin{tikzpicture}[anchorbase,scale=.9]
	\draw[-] (0.4,0) to[out=90, in=0] (0.1,0.5);
	\draw[-] (0.1,0.5) to[out = 180, in = 90] (-0.2,0);
\pinX{(.4,.2)}{}{(1.2,.2)}{f(x)};
\end{tikzpicture}
=
\begin{tikzpicture}[anchorbase,scale=.9]
	\draw[-] (0.4,0) to[out=90, in=0] (0.1,0.5);
	\draw[-] (0.1,0.5) to[out = 180, in = 90] (-0.2,0);
\pinX{(-.2,.2)}{}{(-1.2,.2)}{f(-x)};
\end{tikzpicture}\:$ and similarly for cups,
hence, we have that
\begin{align}
\label{rels9}
\begin{tikzpicture}[anchorbase,scale=1]
	\draw[-] (0.4,0) to[out=90, in=0] (0.1,0.5);
	\draw[-] (0.1,0.5) to[out = 180, in = 90] (-0.2,0);
\circledbar{0.4,.2}{u};
\end{tikzpicture}
&=
\begin{tikzpicture}[anchorbase,scale=1]
	\draw[-] (0.4,0) to[out=90, in=0] (0.1,0.5);
	\draw[-] (0.1,0.5) to[out = 180, in = 90] (-0.2,0);
\circled{-0.2,.2}{u};
\end{tikzpicture}\:,&
\begin{tikzpicture}[anchorbase,scale=1]
	\draw[-] (0.4,0) to[out=90, in=0] (0.1,0.5);
	\draw[-] (0.1,0.5) to[out = 180, in = 90] (-0.2,0);
\circled{0.4,.2}{u};
\end{tikzpicture}
&=
\begin{tikzpicture}[anchorbase,scale=1]
	\draw[-] (0.4,0) to[out=90, in=0] (0.1,0.5);
	\draw[-] (0.1,0.5) to[out = 180, in = 90] (-0.2,0);
\circledbar{-0.2,.2}{u};
\end{tikzpicture}\:,&
\begin{tikzpicture}[anchorbase,scale=1]
	\draw[-] (0.4,0) to[out=-90, in=0] (0.1,-0.5);
	\draw[-] (0.1,-0.5) to[out = 180, in = -90] (-0.2,0);
\circledbar{0.4,-.2}{u};
\end{tikzpicture}
&=
\begin{tikzpicture}[anchorbase,scale=1]
	\draw[-] (0.4,0) to[out=-90, in=0] (0.1,-0.5);
	\draw[-] (0.1,-0.5) to[out = 180, in = -90] (-0.2,0);
\circled{-0.2,-.2}{u};
\end{tikzpicture}\:,&
\begin{tikzpicture}[anchorbase,scale=1]
	\draw[-] (0.4,0) to[out=-90, in=0] (0.1,-0.5);
	\draw[-] (0.1,-0.5) to[out = 180, in = -90] (-0.2,0);
\circled{0.4,-.2}{u};
\end{tikzpicture}
&=
\begin{tikzpicture}[anchorbase,scale=1]
	\draw[-] (0.4,0) to[out=-90, in=0] (0.1,-0.5);
	\draw[-] (0.1,-0.5) to[out = 180, in = -90] (-0.2,0);
\circledbar{-0.2,-.2}{u};
\end{tikzpicture}\:.
\end{align}
Further
useful relations involving these generating functions are
\begin{align}
\begin{tikzpicture}[anchorbase,scale=1.4]
	\draw[-] (0.25,.3) to (-0.25,-.3);
	\draw[-] (0.25,-.3) to (-0.25,.3);
	\circled{-.12,.145}{u};
\end{tikzpicture}
-\begin{tikzpicture}[anchorbase,scale=1.4]
	\draw[-] (0.25,.3) to (-0.25,-.3);
	\draw[-] (0.25,-.3) to (-0.25,.3);
	\circled{.12,-.145}{u};
\end{tikzpicture}
&=
\:\begin{tikzpicture}[anchorbase,scale=1.4]
 	\draw[-] (0,-.3) to (0,.3);
	\draw[-] (-.3,-.3) to (-0.3,.3);
	\circled{0,-.15}{u};
	\circled{-.3,.15}{u};
\end{tikzpicture}
\:-\:\begin{tikzpicture}[anchorbase,scale=1.4]
 	\draw[-] (-0.15,-.3) to[out=90,in=180] (0,-.05) to[out=0,in=90] (0.15,-.3);
 	\draw[-] (-0.15,.3) to[out=-90,in=180] (0,.05) to[out=0,in=-90] (0.15,.3);
	\circled{-.12,.15}{u};
	\circled{.12,-.15}{u};
	\end{tikzpicture}\;,
\label{rels10}&
\begin{tikzpicture}[anchorbase,scale=1.4]
	\draw[-] (0.25,.3) to (-0.25,-.3);
	\draw[-] (0.25,-.3) to (-0.25,.3);
	\circled{-.12,-.145}{u};
\end{tikzpicture}
-\begin{tikzpicture}[anchorbase,scale=1.4]
	\draw[-] (0.25,.3) to (-0.25,-.3);
	\draw[-] (0.25,-.3) to (-0.25,.3);
	\circled{.12,.145}{u};
\end{tikzpicture}
&=
\:\begin{tikzpicture}[anchorbase,scale=1.4]
 	\draw[-] (0,-.3) to (0,.3);
	\draw[-] (-.3,-.3) to (-0.3,.3);
	\circled{0,.15}{u};
	\circled{-.3,-.15}{u};
\end{tikzpicture}\:-\:\begin{tikzpicture}[anchorbase,scale=1.4]
 	\draw[-] (-0.15,-.3) to[out=90,in=180] (0,-.05) to[out=0,in=90] (0.15,-.3);
 	\draw[-] (-0.15,.3) to[out=-90,in=180] (0,.05) to[out=0,in=-90] (0.15,.3);
	\circled{-.12,-.15}{u};
	\circled{.12,.15}{u};
	\end{tikzpicture}\:.
\end{align}
These are also noted in \cite[(2.19)--(2.20)]{BWWbasis}.
Equating the coefficients of $u^{-n-1}$, we obtain
\begin{align}
\label{dotslide1}
\begin{tikzpicture}[anchorbase,scale=1.4]
	\draw[-] (0.25,.3) to (-0.25,-.3);
	\draw[-] (0.25,-.3) to (-0.25,.3);
 \closeddot{-.11,-.135};
 \node at (-.28,-.135) {$\scriptstyle n$};
\end{tikzpicture}
-\begin{tikzpicture}[anchorbase,scale=1.4]
	\draw[-] (0.25,.3) to (-0.25,-.3);
	\draw[-] (0.25,-.3) to (-0.25,.3);
 \closeddot{.12,.135};
 \node at (.28,.135) {$\scriptstyle n$};
\end{tikzpicture}
&=
\sum_{\substack{i,j \geq 0 \\ i+j=n-1}}\left(\begin{tikzpicture}[anchorbase,scale=1.4]
 	\draw[-] (0,-.3) to (0,.3);
	\draw[-] (-.3,-.3) to (-0.3,.3);
	\closeddot{-0.3,-.15};\node at (-.45,-.15) {$\scriptstyle i$};	
 \closeddot{0,.15};\node at (.15,.15) {$\scriptstyle j$};
\end{tikzpicture}\:-\:\begin{tikzpicture}[anchorbase,scale=1.4]
 	\draw[-] (-0.15,-.3) to[out=90,in=180] (0,-.05) to[out=0,in=90] (0.15,-.3);
 	\draw[-] (-0.15,.3) to[out=-90,in=180] (0,.05) to[out=0,in=-90] (0.15,.3);
	\closeddot{-.12,-.15};\node at (-.27,-.15) {$\scriptstyle i$};	
 \closeddot{0.12,.15};\node at (.27,.15) {$\scriptstyle j$};	\end{tikzpicture}\right),\\\label{dotslide2}
 \begin{tikzpicture}[anchorbase,scale=1.4]
	\draw[-] (0.25,.3) to (-0.25,-.3);
	\draw[-] (0.25,-.3) to (-0.25,.3);
 \closeddot{-.12,.135};
 \node at (-.28,.135) {$\scriptstyle n$};
\end{tikzpicture}
-\begin{tikzpicture}[anchorbase,scale=1.4]
	\draw[-] (0.25,.3) to (-0.25,-.3);
	\draw[-] (0.25,-.3) to (-0.25,.3);
 \closeddot{.12,-.135};
 \node at (.28,-.135) {$\scriptstyle n$};
\end{tikzpicture}
&=\sum_{\substack{i,j \geq 0\\i+j=n-1}}
\left(\begin{tikzpicture}[anchorbase,scale=1.4]
 	\draw[-] (0,-.3) to (0,.3);
	\draw[-] (-.3,-.3) to (-0.3,.3);
	\closeddot{-0.3,.15};\node at (-.45,.15) {$\scriptstyle i$};	
 \closeddot{0,-.15};\node at (.15,-.15) {$\scriptstyle j$};
\end{tikzpicture}\:-\:\begin{tikzpicture}[anchorbase,scale=1.4]
 	\draw[-] (-0.15,-.3) to[out=90,in=180] (0,-.05) to[out=0,in=90] (0.15,-.3);
 	\draw[-] (-0.15,.3) to[out=-90,in=180] (0,.05) to[out=0,in=-90] (0.15,.3);
	\closeddot{-.12,.15};\node at (-.27,.15) {$\scriptstyle i$};	
 \closeddot{0.12,-.15};\node at (.27,-.15) {$\scriptstyle j$};
	\end{tikzpicture}\right).
\end{align}

Now consider the ``dotted bubble generating function''
\begin{equation}
\begin{tikzpicture}[anchorbase]
\draw (-.2,0) circle (.2);
\circled{0,0}{u};
\end{tikzpicture}
= 
\sum_{r \geq 0} 
u^{-r-1}\:
\begin{tikzpicture}[anchorbase]
\draw (0,0) circle (.2);
\closeddot{0.2,0};
\node at (0.4,0) {$\scriptstyle r$};
\end{tikzpicture}
\in t u^{-1}1_\one + u^{-2} \End_{\cNB_t}(\one)\llbracket u^{-1}\rrbracket.
\end{equation}
This is often useful, but even more important will be
the renormalization
\begin{equation}\label{deltadef}
\O(u)
 = \sum_{r \geq 0}  u^{-r}\O_r := (-1)^{t}\left(1_\one - 2u \;
\begin{tikzpicture}[anchorbase]
\draw (0,0) circle (.2);
\circled{.2,0}{u};
\end{tikzpicture}
\right)
\in 1_\one + u^{-1}\End_{\cNB_t}(\one)\llbracket u^{-1}\rrbracket.
\end{equation}
Its $u^{-r}$-coefficients $\O_r$ are given explicitly by
\begin{align}\label{Or}
\O_0 &= 1_\one,&
\O_r &= -2(-1)^{t} \:\begin{tikzpicture}[anchorbase]
\draw (0,0) circle (.2);
\closeddot{.2,0};
\node at (.5,0) {$\scriptstyle r$};
\end{tikzpicture}
\end{align}
for $r \geq 1$.
Note also by \cref{RTdots,rels9} that $\O(u)$
is invariant under $\tR$ and $\tT$.

\begin{theorem}[{\cite[Th.~2.5]{BWWbasis}}]\label{nbrelations}
The following relations hold in $\cNB_t$:
\begin{align}
\label{rels11}
2u \:\begin{tikzpicture}[anchorbase]
	\draw[-] (0,-.3)  to (0,-.1) to [out=90,in=180] (.3,.3) to [out=0,in=90] (.5,.1) to[out=-90,in=0] (.3,-.1) to [out=180,in=-90] (0,.3) to (0,.5);
	\circled{.5,.1}{u};
\end{tikzpicture}&=
2u\:
\begin{tikzpicture}[anchorbase]
	\draw[-] (0,-.3) to (0,.5);
 \draw (.5,0.1) circle (.2);
\circled{.7,.1}{u};
\circledbar{0,.1}{u};
\end{tikzpicture}
\,
-\,\begin{tikzpicture}[anchorbase]
	\draw[-] (0,-.3) to (0,.5);
	\circled{0,.1}{u};
\end{tikzpicture}
\:-\:
\begin{tikzpicture}[anchorbase]
	\draw[-] (0,-.3) to (0,.5);
	\circledbar{0,.1}{u};
\end{tikzpicture}
\,
\:,\\\label{rels11a}
\begin{tikzpicture}[anchorbase]
 \draw (-.2,0.1) circle (.2);
	\circled{0,.1}{u};
\end{tikzpicture}
\:+\:
\begin{tikzpicture}[anchorbase]
 \draw (-.2,0.1) circle (.2);
	\circledbar{0,.1}{u};
\end{tikzpicture}
&=
2u\: \begin{tikzpicture}[anchorbase]
 \draw (-.2,0.1) circle (.2);
	\circledbar{0,.1}{u};
\end{tikzpicture}
\:\:\begin{tikzpicture}[anchorbase]
 \draw (-.2,0.1) circle (.2);
	\circled{0,.1}{u};
\end{tikzpicture}
\:,\\
\label{rels12}
\O(u)\O(-u)&=1_\one,\\\label{rels13}
\O(u)\;
\begin{tikzpicture}[anchorbase]
	\draw[-] (0,-.3) to (0,.5);
\end{tikzpicture}&=
\begin{tikzpicture}[anchorbase]
	\draw[-] (0,-.3) to (0,.5);
	\pinX{(0,.1)}{}{(-1,.1)}{\left(\frac{u-x}{u+x}\right)^2};
\end{tikzpicture}\;
\O(u)\:.
\end{align}
\end{theorem}

\begin{corollary}\label{allcurlrels}
The following relations hold in $\cNB_t$:
\begin{align}\label{allcurl1}
2u \:\begin{tikzpicture}[anchorbase]
	\draw[-] (0,-.3)  to (0,-.1) to [out=90,in=180] (.3,.3) to [out=0,in=90] (.5,.1) to[out=-90,in=0] (.3,-.1) to [out=180,in=-90] (0,.3) to (0,.5);
	\circled{.5,.1}{u};
\end{tikzpicture}&=
- (-1)^t \,\begin{tikzpicture}[anchorbase]
	\draw[-] (0,-.3) to (0,.5);
\circledbar{0,.13}{u};
\end{tikzpicture}\:\,
\O(u)-\,\begin{tikzpicture}[anchorbase]
	\draw[-] (0,-.3) to (0,.5);
	\circled{0,.1}{u};
\end{tikzpicture}\:,&
2u \:\begin{tikzpicture}[anchorbase]
	\draw[-] (0,-.3)  to (0,-.1) to [out=90,in=180] (.3,.3) to [out=0,in=90] (.5,.1) to[out=-90,in=0] (.3,-.1) to [out=180,in=-90] (0,.3) to (0,.5);
	\circledbar{.5,.1}{u};
\end{tikzpicture}&=
(-1)^t \,\begin{tikzpicture}[anchorbase]
	\draw[-] (0,-.3) to (0,.5);
\circled{0,.13}{u};
\end{tikzpicture}\:\,
\O(-u)\,+\begin{tikzpicture}[anchorbase]
	\draw[-] (0,-.3) to (0,.5);
	\circledbar{0,.1}{u};
\end{tikzpicture}\:,\\
2u \:\begin{tikzpicture}[anchorbase]
	\draw[-] (0,-.3)  to (0,-.1) to [out=90,in=0] (-.3,.3) to [out=180,in=90] (-.5,.1) to[out=-90,in=180] (-.3,-.1) to [out=0,in=-90] (0,.3) to (0,.5);
	\circledbar{-.5,.1}{u};
\end{tikzpicture}&=
- (-1)^t
\O(u)\,\begin{tikzpicture}[anchorbase]
	\draw[-] (0,-.3) to (0,.5);
\circled{0,.13}{u};
\end{tikzpicture}\;-\,\begin{tikzpicture}[anchorbase]
	\draw[-] (0,-.3) to (0,.5);
	\circledbar{0,.1}{u};
\end{tikzpicture}\:,&
2u \:\begin{tikzpicture}[anchorbase]
	\draw[-] (0,-.3)  to (0,-.1) to [out=90,in=0] (-.3,.3) to [out=180,in=90] (-.5,.1) to[out=-90,in=180] (-.3,-.1) to [out=0,in=-90] (0,.3) to (0,.5);
	\circled{-.5,.1}{u};
\end{tikzpicture}&=
(-1)^t
\O(-u)\,\begin{tikzpicture}[anchorbase]
	\draw[-] (0,-.3) to (0,.5);
\circledbar{0,.13}{u};
\end{tikzpicture}\,+\,\begin{tikzpicture}[anchorbase]
	\draw[-] (0,-.3) to (0,.5);
	\circled{0,.1}{u};
\end{tikzpicture}\:.\label{allcurl2}
\end{align}
\end{corollary}

\begin{proof}
The first equality follows from \cref{rels11}
and the definition \cref{deltadef}.
The others follow by applying $\tR$
or using the substitution $u \mapsto -u$.
\end{proof}

\begin{corollary}\label{curlstuff}
For $n \geq 0$, we have that 
\begin{align*}
\begin{tikzpicture}[anchorbase]
	\draw[-] (0,-.3)  to (0,-.1) to [out=90,in=180] (.3,.3) to [out=0,in=90] (.5,.1) to[out=-90,in=0] (.3,-.1) to [out=180,in=-90] (0,.3) to (0,.5);
  \closeddot{.5,.1};\node at (.9,0.1) {$\scriptstyle n+1$};
\end{tikzpicture}&=
\sum_{r=0}^{n-1}(-1)^{r}
\begin{tikzpicture}[anchorbase]
	\draw[-] (-.5,-.3) to (-.5,.3);
\draw (0,0) circle (.2);
\closeddot{0.2,0};
\closeddot{-0.5,0};
\node at (-0.7,0) {$\scriptstyle r$};
\node at (0.6,0) {$\scriptstyle n-r$};
\end{tikzpicture}
-\delta_{n\equiv t}\begin{tikzpicture}[anchorbase]
	\draw[-] (0,-.3) to (0,.3);
 \closeddot{0,0};\node at (.2,0) {$\scriptstyle n$};
\end{tikzpicture}\:.
\end{align*}
\end{corollary}

\begin{proof}
This follows by equating the coefficients of $u^{-n-1}$ in the first equation from \cref{allcurl1}.
\end{proof}

\subsection{The basis theorem}\label{sqfunc}
Let $\Lambda$ be the graded algebra of 
symmetric functions over $\kk$. Adopting standard notation, this is freely generated
either by the elementary symmetric functions $e_r\:(r > 0)$
or by the complete symmetric functions $h_r\:(r > 0)$; our convention for the grading puts
these in degree $2r$.
The two families of generators are related by the identity
\begin{equation}\label{grassmannian}
e(-u) h(u) = 1
\end{equation}
where 
\begin{align}\label{genfuncs}
e(u) &= \sum_{r \geq 0} u^{-r} e_r,&
h(u) &= \sum_{r \geq 0} u^{-r} h_r
\end{align}
are the corresponding generating functions, and $e_0 = h_0 = 1$ by convention.
It is also convenient to interpret $e_r$ and $h_r$ as $0$ when $r < 0$.

Following \cite[Ch.~III, Sec.~8]{Mac}, we define a power series
$q(u) \in \Lambda \llbracket u^{-1}\rrbracket$ 
and elements $q_r\;(r \geq 0)$ of $\Lambda$ so that
\begin{equation}\label{qdef}
q(u) = \sum_{r\geq 0} u^{-r} q_r := e(u)h(u).
\end{equation}
By \cref{grassmannian},
we have that 
\begin{equation}\label{qgrassmannian}
q(u) q(-u) = 1.
\end{equation}
Equivalently, $q_0 = 1$ and
\begin{align}\label{muscles}
q_{2r} &= (-1)^{r-1}{\textstyle\frac{1}{2}}q_r^2 +
\sum_{s=1}^{r-1}(-1)^{s-1}q_s q_{2r-s}
\end{align}
for $r \geq 1$; cf. \cite[(III.8.2$'$)]{Mac}.
As with $e_r$ and $h_r$, we adopt the convention that $q_r = 0$ for $r < 0$.

The graded subalgebra of $\Lambda$ generated by all $q_r\:(r \geq 0)$ is denoted
$\Gamma$. As explained in \cite{Mac}, 
$\Gamma$ is freely generated by $q_1,q_3,q_5,\dots$
(and it has a distinguished basis given by the {\em Schur $Q$-functions} $Q_\lambda$ indexed by all strict partitions).
It follows that $\Gamma$ is generated by 
the elements $q_r\:(r \geq 0)$ subject only to the relations
\cref{qgrassmannian}.
Hence, \cref{rels12} is all that is needed to establish the existence of a
graded algebra homomorphism
\begin{align}\label{gammacor}
\gamma_t:\Gamma &\rightarrow \End_{\cNB_t}(\one),
&q_r
&\mapsto \O_r.
\end{align}
By \cite[Cor.~5.4]{BWWbasis}, this is actually an
{\em isomorphism}.

Now we recall the basis theorem for morphism spaces in $\cNB_t$, which is the main result of \cite{BWWbasis}. For $m, n \geq 0$,
any morphism $f:B^{\star n} \rightarrow B^{\star m}$
is represented by a linear combination of {\em $m \times n$ string diagrams}, i.e., string diagrams with $m$ boundary points at the top and $n$ boundary points at the bottom
that are obtained by composing the generating
morphisms from \cref{gens}. 
It follows that $\Hom_{\cNB_t}(B^{\star n}, B^{\star m})$ is 0 unless $m \equiv n \pmod{2}$.
The individual strings in 
an $m \times n$ string diagram $s$
are of four basic types:
cups (with two boundary points on the top edge),
caps (with two boundary points on the bottom edge),
propagating strings (with one boundary point at the top and one at the bottom),
and internal bubbles (no boundary points).
We define an equivalence relation $\sim$
on the set of $m \times n$ string diagrams
by declaring that $s \sim s'$ if
their strings define the same matching on the set of $m+n$ boundary points. 
We say that $s$ is
{\em reduced} if the following properties hold:
\begin{itemize}
\item 
There are no internal bubbles.
\item Propagating strings have no critical points ($=$points of slope 0).
\item Cups and caps each have exactly one critical point.
\item There are no {\em double crossings} ($=$ two different strings which cross each other at least twice).
\end{itemize}
These assumptions imply in particular that there are no {\em self-intersections} ($=$ crossings of a string with itself).
Fix a set $\overline{\RSD}(m\times n)$ of representatives for the
$\sim$-equivalence classes
of {\em undotted} reduced $m \times n$ string diagrams.
The total number of such diagrams is $(m+n-1)!!$ if $m \equiv n \pmod{2}$,
and there are none otherwise.
For each of these $\sim$-equivalence
class representatives, we also choose distinguished points in the interior of each of its strings that are away from points of intersection. Then
let $\RSD(m\times n)$ be the
set of all morphisms
$f:B^{\star n} \rightarrow B^{\star m}$ which can be obtained
by taking an element of $\overline{\RSD}(m\times n)$ then adding dots
labeled by non-negative multiplicities
at each of the distinguished points on the strings.

\begin{theorem}[{\cite[Th.~5.1]{BWWbasis}}]
\label{basisthm} 
Viewed as a graded $\Gamma$-module so that
$p \in \Gamma$ acts on $f:B^{\star n}\rightarrow B^{\star m}$
by $f \cdot p := f \star \gamma_t(p)$,
the space $\Hom_{\cNB_t}(B^{\star n}, B^{\star m})$
is free with basis $\RSD(m\times n)$.
\end{theorem}

Now we can make the first significant connection
between $\cNB_t$ and the iquantum group.
Recall the bilinear form 
$(\cdot,\cdot)^\imath:\Ui \times \Ui \rightarrow \Q(q)$ from \cref{otherformt}.

\begin{theorem}\label{toinfinityandbeyond}
For $m,n \in \N$,
we have that
$\Hom_{\cNB_t}(B^{\star n}, B^{\star m}) \cong \Gamma^{\oplus \overline{(b^m, b^n)^\imath}}$
as a graded $\Gamma$-module.
\end{theorem}

\begin{proof}
We compare the explicit combinatorial formula for
$(b^m, b^n)^\imath$ from \cref{classic}
with the graded rank of 
$\Hom_{\cNB_t}(B^{\star n}, B^{\star m})$ as a free graded $\Gamma$-module computed via \cref{basisthm}.
If $m \not \equiv n \pmod{2}$ then
$(b^m,b^n)^\imath = 0$ and $\RSD(m\times n)$
is empty, and the result is clear.
Now assume that $m \equiv n \pmod{2}$
and let $f := (m+n)/2$.
There is an obvious 
bijection between equivalence classes of $m \times n$ string diagrams
and chord diagrams with $f$ free chords and no tethered chords. This just arises
by identifying the $(m+n)$ boundary points of
strings in an $m \times n$ string diagram
with the $(m+n)$ endpoints of chords in a chord diagram in some fixed way that preserves the clockwise ordering, then replacing strings by chords so that the underlying matching of these points is preserved.
In a string diagram, each crossing is of degree $-2$, so it
contributes $q^{-2}$ to the graded rank. The dots placed at the $f$ distinguished points produce
the factor $1 / (1-q^{2})^f$, this being 
$\dim_{q} \kk[x_1,\dots,x_f]$ with 
$x_i$ in degree 2.
Recalling the definition of
the generating function
$T_{f,0}(q)$ from \cref{high},
we deduce that 
$$
\rank_q \Hom_{\cNB_t}(B^{\star n}, B^{\star m})
 = 
 \sum_{s \in \RSD(m\times n)}
 q^{\deg(s)}
 = T_{f,0}(q^{-2}) / (1-q^{2})^f,
 $$
which is $\overline{(b^m,b^n)^\imath}$
according to \cref{classic}.
\end{proof}

\subsection{Central elements}\label{sscentral}
Recall that the {\em center} $Z(\cA)$
of a category $\cA$
means the algebra of endomorphisms of its identity 
endofunctor.
Thus, elements of $Z(\cNB_t)$ consist of
tuples $(z_n)_{n \geq 0}$ 
for elements $z_n \in \End_{\cNB_t}(B^{\star n})$
such that $z_m \circ f = f \circ z_n$
for all $m,n \geq 0$ and $f \in \Hom_{\cNB_t}(B^{\star n},
B^{\star m})$. In this subsection, we are going to use the dotted bubbles to construct many---conjecturally, all---elements of $Z(\cNB_t)$.

Since $\O(\pm u) \in 1_\one + u^{-1} \End_{\cNB_t}(\one)\llbracket u^{-1} \rrbracket$ and $2$ is invertible in $\kk$, it makes sense to take the square roots
$\sqrt{\O(\pm u)}$; we choose the ones that are positive in the sense that they again lie in
$1_\one + u^{-1} \End_{\cNB_t}(\one)\llbracket u^{-1} \rrbracket$.
We have that $\sqrt{\O(-u)} = \left(\sqrt{\O(u)}\right)^{-1}$ by \cref{rels12}.
Taking the square roots of
both sides of \cref{rels13}, both of which lie in $1_B + u^{-1} \End_{\cNB_t}(B) \llbracket u^{-1} \rrbracket$, we obtain
\begin{align}
\label{relssqrt1}
{\scriptstyle\sqrt{\O(u)}}\:\:
\begin{tikzpicture}[anchorbase]
	\draw[-] (0,-.3) to (0,.5);
	\circled{0,.1}{u};
\end{tikzpicture}&=
\begin{tikzpicture}[anchorbase]
	\draw[-] (0,-.3) to (0,.5);
	\circledbar{0,.1}{u};
\end{tikzpicture}\;
{\scriptstyle\sqrt{\O(u)}}\:,&
{\scriptstyle\sqrt{\O(-u)}}\:\:
\begin{tikzpicture}[anchorbase]
	\draw[-] (0,-.3) to (0,.5);
	\circledbar{0,.1}{u};
\end{tikzpicture}&=
\begin{tikzpicture}[anchorbase]
	\draw[-] (0,-.3) to (0,.5);
	\circled{0,.1}{u};
\end{tikzpicture}\;
{\scriptstyle\sqrt{\O(-u)}}\:.
\end{align}

Let $e_{r,n}, h_{r,n}, q_{r,n} \in \kk[x_1,\dots,x_n]^{S_n}$
be the symmetric polynomials in $n$ variables obtained by specializing the symmetric functions
$e_r, h_r, q_r$ from \cref{genfuncs,qdef}.
We have 
that
\begin{equation}
q_{r,n} = \sum_{s=0}^r e_{s,n} h_{r-s,n}.
\end{equation}
Moreover,
\begin{equation}\label{question}
\sum_{r \geq 0} u^{-r}q_{r,n} =
\prod_{i=1}^n \frac{u+x_i}{u-x_i} 
\in 1+u^{-1} \kk[x_1,\dots,x_n]\llbracket u^{-1} \rrbracket.
\end{equation}
In the statement of the next theorem,
for a polynomial $f \in \kk[x_1,\dots,x_n]$,
we use the notation $f 1_n = 1_n f$ to denote the endomorphism of
$B^{\star n}$ defined by
interpreting $x_i$ as 
$|^{\star (i-1)}\! \star\!\!\smalltxtdot
\star |^{\star (n-i)}$, i.e., the dot 
on the $i$th string.

\begin{theorem}\label{bankrupt}
For any $r \geq 0$, we have that $\left(q_{r,n} 1_n\right)_{n \geq 0}\in Z(\cNB_t)$.
\end{theorem}

\begin{proof}
We need to show that
$q_{r,m} 
1_m \circ f = f \circ q_{r,n} 1_n$
for any $f \in \Hom_{\cNB_t}\left(B^{\star n},
B^{\star m}\right)$.
By \cref{relssqrt1}, we have that
\begin{equation}\label{trp}
\sum_{r \geq 0} u^{-r} q_{r,n} 1_n
=\prod_{i=1}^n \frac{u+x_i}{u-x_i} 1_n=
\begin{tikzpicture}[anchorbase]
	\draw[-] (0,-.5) to (0,.5);
 \node at (-0.33,-.25) {$\scriptstyle -1$};
	\draw[-] (.4,-.5) to (.4,.5);
 \node at (0.73,-.25) {$\scriptstyle -1$};
	\draw[-] (1.25,-.5) to (1.25,.5);
\node at (.9,0) {$\cdots$};
\node at (1.58,-.25) {$\scriptstyle -1$};
	\circledbar{0,-.25}{u};
	\circledbar{.4,-.25}{u};
	\circledbar{1.25,-.25}{u};
	\circled{.4,.25}{u};
	\circled{1.25,.25}{u};
	\circled{0,.25}{u};
\end{tikzpicture}=
{\scriptstyle\sqrt{\O(-u)}}\:
 \star  |^{\star n} \star {\scriptstyle\sqrt{\O(u)}}.
\end{equation}
The result follows from this since the expression on the right hand side clearly has the desired property by the interchange law.
\end{proof}

\begin{corollary}\label{centralstuff}
Let
$p_{r,n} := \sum_{i=1}^n x_i^r \in \kk[x_1,\dots,x_n]^{S_n}$
be the $r$th power sum.
For any {\em odd} $r \geq 1$,
we have that
$(p_{r,n} 1_n)_{n \geq 0}
\in Z(\cNB_t)$.
\end{corollary}

\begin{proof}
It suffices to note that any odd power sum can be written as a polynomial in the symmetric polynomials $q_{r,n}$. This can be proved by taking the
logarithmic derivative of \cref{question}.
\end{proof}

\subsection{Cyclotomic quotients}\label{minpols}
In this subsection, the grading plays no role---we 
view $\cNB_t$ simply as a strict $\kk$-linear monoidal category. 
Let $\cV$ be a strict $\cNB_t$-module category. 
This means that $\cV$ is a $\kk$-linear category, and
we are given a strict $\kk$-linear
monoidal functor
$\mu$ from $\cNB_t$ to the strict $\kk$-linear 
monoidal category $\END_\kk(\cV)$ 
whose objects are $\kk$-linear endofunctors of $\cV$
and whose morphisms are natural transformations.
We denote the endofunctor $\mu(B):\cV\rightarrow \cV$ simply by $B$. The string calculus for monoidal categories extends to the module category $\cV$ 
in the usual way. We will represent the identity
endomorphism of an object $V$ of $\cV$ by the labeled string 
$\;\begin{tikzpicture}[anchorbase]
\draw[-,objcolor,thick] (0.68,-.2)\botlabel{V} to (0.68,.2);
\end{tikzpicture}\ $ on the right hand boundary of the diagram.
For a string diagram $s$ representing
a morphism in $\Hom_{\cNB_t}(B^{\star n},B^{\star m})$,
we represent the morphism
$\mu(s)_V:B^n V \rightarrow B^m V$ diagrammatically
by $s\;\begin{tikzpicture}[baseline=-2pt]
\draw[-,objcolor,thick] (0.68,-.2)\botlabel{V} to (0.68,.2);
\end{tikzpicture}$.

Let $V \in \cV$ be any object. We write $Z_V$ for the center of its endomorphism algebra.
Any dotted bubble evaluates to an element of $Z_V$. In particular,
the generating function 
$\O(u)$ for dotted bubbles from \cref{deltadef} 
evaluates to a power series $\O_V(u) \in Z_V\llbracket u^{-1}\rrbracket$. Let
\begin{equation}\label{washburne}
\theta_V: Z_V[u] \rightarrow \End_{\cV}(BV)
\end{equation}
be the $\kk$-algebra homomorphism
mapping $z \in Z_V$ to $B z : BV \rightarrow BV$ and
$u$ to the endomorphism of $BV$ defined by the dot.
We represent $\theta_V(f(u))$ diagrammatically by
$\ \begin{tikzpicture}[anchorbase]
\draw (1.5,.3) to (1.5,-.3);
\pinX{(1.5,0)}{}{(.86,0)}{f(x)};
\draw[-,objcolor,thick] (1.85,.3) to (1.85,-.3)\botlabel{V};
\end{tikzpicture}\ $.
The coefficients of $f(x)$ are elements of $Z_V$ which could be represented by coupons on the right-hand string, but we will not attempt to do that here.
Let $I_V := \ker \theta_V \unlhd Z_V[u]$. 

\begin{theorem}\label{flash}
Let $V$ be an object of 
an $\cNB_t$-module category $\cV$ as above,
and $f(u) \in I_V$.
The formal Laurent series
$f^*(u) := (-1)^t f(-u) \O_V(-u) \in Z_V\lround u^{-1}\rround$
is a polynomial in $I_V$ with the same constant term as $f(u)$, and $f^{**}(u) = f(u)$. 
Moreover,
$\tfrac{f(u)+f^*(u)}{2}$
and $\tfrac{f(u)-f^*(u)}{2u}$ are $*$-invariant elements of $I_V$.
Assuming in addition that $f(u) \in I_V$ is monic,
$$
g(u) := \begin{dcases}
\tfrac{f(u)+f^*(u)}{2}&\text{if $\deg f(u) \equiv t \pmod{2}$}\\
\tfrac{f(u)-f^*(u)}{2u}&\text{if $\deg f(u) \not\equiv t \pmod{2}$}
\end{dcases}
$$
is a monic polynomial in $I_V$
with $\deg g(u) = \deg f(u)$ or $\deg g(u) = \deg f(u)-1$ in the two cases, such that $\O_V(u) = (-1)^t g(-u) / g(u)$.
\end{theorem}

\begin{proof}
For $r \geq 0$, we have that
\begin{align*}
\left[f(-u) - f^*(u)\right]_{u^{-r}} &= 
\left[f(-u) - (-1)^t f(-u) \O_V(-u)\right]_{u^{-r}}
=
\left[u^r f(-u) \tfrac{\id_V - (-1)^t \O_V(-u)}{u}\right]_{u^{-1}}\\
&\!\!\stackrel{\cref{deltadef}}{=}
\left[u^r f(-u)\ \ 
\begin{tikzpicture}[anchorbase]
\draw[-] (-0.25,0) arc(180:-180:0.25);
\circledbar{.25,0}{u};
\draw[-,objcolor,thick] (.7,.4) to (.7,-.4)\botlabel{V};
\end{tikzpicture}\ \right]_{u^{-1}}
\stackrel{\cref{trick}}{=}
\begin{tikzpicture}[anchorbase]
\draw[-] (-0.25,0) arc(180:-180:0.25);
\pinX{(.25,0)}{}{(1.3,0)}{(-x)^r f(x)};
\draw[-,objcolor,thick] (2.2,.4) to (2.2,-.4)\botlabel{V};
\end{tikzpicture}
=0.
\end{align*}
Thus, for $r \geq 0$, we have that
$\left[f^*(u)\right]_{u^{-r}}
=\left[f(-u)\right]_{u^{-r}} = \delta_{r,0} f(0)$.
This shows
that $f^*(u)$ is a polynomial in $Z_V[u]$ with the same constant term as $f(u)$. 
Also $f^{**}(u) = f(u) \O_V(u)\O_V(-u)$, which equals $f(u)$ by \cref{rels12}.

Since $f^*(u)$ and $f(u)$ are polynomials with the same constant terms,
$\frac{(-u)^r f^*(u) - u^r f(u)}{2u}$
is a polynomial in $u$ for any $r \geq 0$.
So it is an element of $Z_V[u]$.
Its image under $\theta_V$ is
\begin{align*}
\begin{tikzpicture}[anchorbase]
\draw[-] (0.2,-.4) to (0.2,.4);
\pinX{(.2,0)}{}{(-1.6,0)}{\tfrac{(-x)^r f^*(x)-x^r f(x)}{2x}};
\draw[-,objcolor,thick] (0.6,.4) to (0.6,-.4)\botlabel{V};
\end{tikzpicture}
\ \ &\!\!\!\stackrel{\cref{trick}}{=}
\left[
\tfrac{(-u)^r f^*(u)-u^r f(u)}{2u}\ \ 
\begin{tikzpicture}[anchorbase]
\draw[-] (0.2,-.4) to (0.2,.4);
\circled{.2,0}{u};
\draw[-,objcolor,thick] (0.6,.4) to (0.6,-.4)\botlabel{V};
\end{tikzpicture}\ 
\right]_{u:-1}=
\left[
\tfrac{(-u)^r f^*(u)}{2u}\ \ 
\begin{tikzpicture}[anchorbase]
\draw[-] (0.2,-.4) to (0.2,.4);
\circled{.2,0}{u};
\draw[-,objcolor,thick] (0.6,.4) to (0.6,-.4)\botlabel{V};
\end{tikzpicture}
-\tfrac{u^r f(u)}{2u}\ \ 
\begin{tikzpicture}[anchorbase]
\draw[-] (0.2,-.4) to (0.2,.4);
\circled{.2,0}{u};
\draw[-,objcolor,thick] (0.6,.4) to (0.6,-.4)\botlabel{V};
\end{tikzpicture}\ 
\right]_{u^{-1}}\\
&=\left[
\tfrac{(-u)^r (-1)^t f(-u) \O_V(-u)}{2u}\ \ 
\begin{tikzpicture}[anchorbase]
\draw[-] (0.2,-.4) to (0.2,.4);
\circled{.2,0}{u};
\draw[-,objcolor,thick] (0.6,.4) to (0.6,-.4)\botlabel{V};
\end{tikzpicture}
+\tfrac{(-u)^r f(-u)}{2u}\ \ 
\begin{tikzpicture}[anchorbase]
\draw[-] (0.2,-.4) to (0.2,.4);
\circledbar{.2,0}{u};
\draw[-,objcolor,thick] (0.6,.4) to (0.6,-.4)\botlabel{V};
\end{tikzpicture}\ 
\right]_{u^{-1}}\\
&=\left[
(-u)^r f(-u)\left(\tfrac{(-1)^t}{2u} \ \ 
\begin{tikzpicture}[anchorbase]
\draw[-] (-0.6,-.4) to (-0.6,.4);
\circled{-.6,0}{u};
\node at (0.35,0) {$\O_V(-u)$};
\draw[-,objcolor,thick] (1.15,.4) to (1.15,-.4)\botlabel{V};
\end{tikzpicture}
+\tfrac{1}{2u}\ \ 
\begin{tikzpicture}[anchorbase]
\draw[-] (0.2,-.4) to (0.2,.4);
\circledbar{.2,0}{u};
\draw[-,objcolor,thick] (0.6,.4) to (0.6,-.4)\botlabel{V};
\end{tikzpicture}\right)
\right]_{u^{-1}}\\
&\!\!\!\stackrel{\cref{allcurl1}}{=}\left[
(-u)^r f(-u)\ \ 
\begin{tikzpicture}[anchorbase]
\draw[-] (0,-0.5) to[out=up,in=180] (0.3,0.2) to[out=0,in=up] (0.45,0) to[out=down,in=0] (0.3,-0.2) to[out=180,in=down] (0,0.5);
\circledbar{.45,0}{u};
\draw[-,objcolor,thick] (.8,.5) to (.8,-.5) \botlabel{V};
\end{tikzpicture}\ \right]_{u^{-1}}
\stackrel{\cref{trick}}{=}\ 
\begin{tikzpicture}[anchorbase]
\draw[-] (0,-0.5) to[out=up,in=180] (0.3,0.2) to[out=0,in=up] (0.45,0) to[out=down,in=0] (0.3,-0.2) to[out=180,in=down] (0,0.5);
\pinX{(.45,0)}{}{(1.3,0)}{x^r f(x)};
\draw[-,objcolor,thick] (2,.4) to (2,-.4) \botlabel{V};
\end{tikzpicture} = 0.
\end{align*}
Thus, we have proved that 
$\frac{(-u)^r f^*(u) - u^r f(u)}{2u} \in I_V$ for 
any $r \geq 0$. The $r=0$ and $r=1$ cases show that
$\tfrac{f(u)+f^*(u)}{2}$
and $\tfrac{f(u)-f^*(u)}{2u}$ both belong to $I_V$, which implies that $f^*(u) \in I_V$. Moreover, it is easy to check using $f^{**}(u) = f(u)$
that both of these elements are $*$-invariant.

Finally, we assume that $f(u)$ is monic
and consider the two cases separately:
\begin{itemize}
\item If $\deg f(u) \equiv t \pmod{2}$, 
then both $f(u)$ and $f^*(u)$ are monic, so
$g(u) := \frac{f(u) + f^*(u)}{2} \in I_V$ is monic of the same degree as $f(u)$.
\item
If $\deg f(u) \not\equiv t \pmod{2}$
then both $f(u)$ and $-f^*(u)$ are monic, so
 $g(u) := \frac{f(u)-f^*(u)}{2u} \in I_V$
is monic 
of one less degree than $f(u)$.
\end{itemize}
The equality $g^*(u) = g(u)$ means that
$(-1)^t g(-u) \O_V(-u) = g(u)$, so 
$g(u) \O_V(u) = (-1)^t g(-u)$.
Since $g(u)$ is monic, it is invertible in $Z_V\lround u^{-1}\rround$, and 
the final assertion of the theorem follows.
\end{proof}

Some of the ideas in the proof of \cref{flash} were suggested by Anthropic's Claude AI. 
In the first version of this paper, 
we were only able to prove the slightly weaker result in the following corollary, which 
is all is needed later on in the article.

\begin{corollary}\label{lynch}
Suppose that $V$ is an object of $\cV$
such that $\End_\cV(V) = \kk$
and $\dim \End_{\cV}(B V) < \infty$.
The ideal $I_V$ of the PID $\kk[u]$ is non-zero,
so it is generated by the unique
monic polynomial $m_V(u)$ of smallest degree belonging to $I_V$.
Moreover, 
$\deg m_V(u) \equiv t \pmod{2}$ and
$\O_V(u) = (-1)^t m_V(-u) / m_V(u)$.
\end{corollary}

\begin{proof}
The kernel is non-zero since
$\theta_V$ is a homomorphism from the infinite-dimensional algebra
$\kk[u]$ to $\End_{\cV}(BV)$, which is finite dimensional.
Now we apply \cref{flash} taking $f(u) := m_V(u)$.
If $\deg m_V(u) \not\equiv t \pmod{2}$, the theorem shows that
$\frac{1}{2u}(m_V(u) - m_V^*(u))$ is a monic polynomial in $I_V$
of smaller degree than $m_V(u)$, which is a contradiction.
Thus, $\deg m_V(u) \equiv t \pmod{2}$. Now the theorem shows that
$g(u) := \frac{1}{2}(m_V(u)+m_V^*(u))$ is a monic polynomial of the same degree as $m_V(u)$. Since $m_V(u)$ is the unique monic polynomial in $I_V$ of this degree, it follows that $g(u) = m_V(u)$. Consequently, the last part of \cref{flash} shows that
$\O_V(u) = (-1)^t m_V(-u) / m_V(u)$.
\end{proof}

\begin{remark}
In this remark---not needed subsequently in the article---we briefly explain why \cref{flash} is significant. Suppose that
$l \geq 0$ with $l \equiv t \pmod{2}$.
Let $\Lambda_l := \kk[x_1,\dots,x_l]^{S_l}$ be the algebra of symmetric polynomials and $e_l(u) := (u-x_1)\cdots (u-x_l)
\in \Lambda_l[u]$.
Let $\CNB_l$ be the cyclotomic nil-Brauer category
from \cite[Def.~2.11]{BBE}. 
It is a strict $\kk$-linear monoidal category in its own right, but all that is important here is that $\CNB_l$ is 
a strict $\cNB_t$-module category.
An equivalent way to formulate the 
definition of $\CNB_l$ is that it is 
the quotient of the $\Lambda_l$-linear 
$\cNB_t$-module category
$\cNB_t \otimes_\kk \Lambda_l$
by the $\Lambda_l$-linear submodule category
generated by the relations
$\O(u) = (-1)^t e_l(-u) / e_l(u)$
and 
$\begin{tikzpicture}[anchorbase]
\draw (1.5,.3) to (1.5,-.3);
\pinX{(1.5,0)}{}{(.8,0)}{e_l(x)};
\end{tikzpicture}=0$.
Moreover, we have that
$\End_{\CNB_l}(\one)=\Lambda_l$.
Indeed, it is easy to see from the definition
that the canonical map $\iota:\Lambda_l \rightarrow
\End_{\CNB_l}(\one)$ is surjective, and $\iota$ is injective because
its composition with 
the algebra homomorphism $\End_{\CNB_l}(\one) \rightarrow \Lambda_l$ defined by the functor $\Theta$ from
\cite[Th.~5.1]{BBE} is the identity.

Now let $V \in \cV$ be as above,
and assume that $\dim \End_{\cV}(BV) < \infty$. This assumption ensures that there exist monic polynomials $f(u) \in I_V$. Hence, by \cref{flash},
there exists a 
monic polynomial $g(u) \in I_V$
such that $l := \deg g(u) \equiv t \pmod{2}$ and 
$\O_V(u) = (-1)^t g(-u) / g(u)$.
In view of the defining relations for $\CNB_l$ described in the previous paragraph, it follows that there is
a unique morphism of strict $\cNB_t$-module categories
$\operatorname{ev}_V:
\CNB_l \rightarrow \cV$ mapping the generating object $\one$
of $\CNB_l$ to $V$ and the polynomial 
$e_l(u) \in \End_{\CNB_l}(\one)[u]$ 
to $g(u) \in \End_\cV(V)[u]$. This makes precise the idea that $\CNB_l$ is the {\em universal cyclotomic quotient} of $\cNB_t$  of level $l$. 

For a description of the Grothendieck group $K_0(\CNB_l)$ as a left $K_0(\cNB_t)$-module, see \cite{BBE}.
\end{remark}

%% file: s4-idempotents.tex
\setcounter{section}{3}

\section{Primitive idempotents}\label{sect4}

In this section, we work out the structure of the primitive homogeneous idempotents in $\cNB_t$
and prove \cref{introthm:PIM,introthm:iso}.
We continue to work over the field $\kk$
of characteristic different from 2.

\subsection{Extended graphical calculus}
We begin by introducing some further diagrammatical shorthands
in the spirit of the
``thick calculus'' of \cite{KLMS}.
We denote
the tensor product $|^{\star a}$
of $a$ strings by a single thick string labeled by $a$.
A thick cup or cap labeled by $a$ denotes that number of nested ordinary cups or caps (no crossings).
Sometimes it is notationally convenient to be able to split thick strings into thinner ones or to merge thinner strings to obtain thicker ones:
the diagrams
\begin{align*}

\end{align*}
simply represent the identity morphisms $B^{\star n} \rightarrow B^{\star a} \star B^{\star b}$
and $B^{\star a} \star B^{\star b} 
\rightarrow B^{\star n}$
for $a+b=n$.
We will often omit a thickness label on a thick
string when it can be inferred from others in the diagram.

For $a+b=n$, the thick crossing
$$
%
$$
denotes the morphism
$B^{\star a} \star B^{\star b}
\rightarrow B^{\star b} \star B^{\star a}$
obtained by composing ordinary 
crossings according
to a reduced expression for the longest
of the minimal length $S_n / (S_a \times S_b)$-coset representatives.
We use a thick string decorated with a cross
to denote the composition of thin crossings
corresponding to a reduced expression for the longest element $w_n$.
For example:
\begin{align*}
%
\:\;.
\end{align*}
When working with these morphisms, 
we will often make implicit use of 
various obvious consequences of the braid
relations, 
such as
\begin{align*}
%
\:.
\end{align*}
In view of the pitchfork relations, one can also draw this cross at the critical point of a thick cup or cap without there being any ambiguity as to the meaning:
\begin{align*}
%
\:.
\end{align*}

We use a dot on a string of thickness $n$ labeled by 
$\alpha \in \N^n$
to denote
the tensor product of dots on ordinary strings
labeled by the parts of $\alpha$:
$$
%
$$
The $n$-tuples
$\rho_n := (n-1,n-2,\dots,1,0) \in \N^n$
and $\varpi_{r,n} := (1,\dots,1,0,\dots,0) \in \N^n$ with $r$ entries equal to $1$ followed by
$(n-r)$ entries equal to $0$
will appear often.
To simplify notation, we allow the subscript $n$ to be omitted in these when 
used to label a node on a string of thickness $n$:
\begin{align*}
%
\:.
\end{align*}
Generalizing the notation \cref{labelledpin},
given a polynomial $f=
\sum_{\alpha \in \N^n} c_\alpha x_1^{\alpha_1} \cdots x_n^{\alpha_n}
\in \kk[x_1,\dots,x_n]$, the pin
$$
%
$$
denotes the endomorphism $f 1_n = 1_n f$
of $B^{\star n}$.
Often for this $f$ will be the elementary symmetric polynomial 
$e_{r,n} := \sum_{1 \leq i_1 < \cdots < i_r \leq n} x_{i_1} \cdots x_{i_r}$.
Again, if this is pinned to a string of thickness $n$, we allow the subscript $n$ to be dropped,
writing simply
\begin{align*}
%
\end{align*}
since the number $n$ of variables in the elementary symmetric polynomial can be inferred from the thickness of the string.

\begin{lemma}\label{birdie}
For $0 \leq r \leq n$, we have that
\begin{align}\label{birdier}
%
\:.
\end{align}
\end{lemma}

\begin{proof}
The first equality is the well-known
identity $e_{r,n} = \sum_{s=0}^r (-1)^{r-s} 
x_{n+1}^{r-s}e_{s,n+1}$.
Then the second equality follows on applying $\tR$.
\end{proof}

\begin{lemma}\label{toobusy}
For $0 \leq i \leq n$, we have that
\begin{align}
%
\:.
\end{align}
\end{lemma}

\begin{proof}
We just prove the first identity; the second then follows on applying $\tR$. 
By \cref{basisthm}, the lowest non-zero degree of $\End_{\cNB_t}(B^{\star (n+1)})$ is 
$-n(n+1)$, and the diagram on the left hand side of the identity is of degree 
$-n(n-1)-4n+2i$. If $i < n$ 
then $-n(n-1)-4n+2i < -n(n+1)$ so the expression is 0.

To prove the result in the remaining case that $i=n$, we proceed by induction on $n$. Assume the result is true for $n$ and consider the next case
$$
%
\right)\:.
$$
In this expression, the term before the summation is 0 by the degree argument given already, the first term in the summation is 0 unless $a=0$ by the defining relations \cref{rels1}, and similarly the second term in the summation is 0 unless $b=0$. So
$$
%
\:,
$$
where we used the induction hypothesis for the second equality.
\end{proof}

\begin{corollary}\label{curly}
For $0 \leq i \leq n+1$, we have that
\begin{align}\label{curlylabel}
%
\:.
\end{align}
\end{corollary}

\begin{proof}
As usual, we just prove the first equality.
By the braid relation then \cref{curlstuff} and \cref{toobusy}, we get that
$$
%
\:.
$$
\end{proof}

\begin{corollary}\label{imanant}
For any $n \geq 1$, we have that $
%
$.
\end{corollary}

\begin{proof}
This follows by induction on $n$. For the induction step,
we have that
$$
%
\:,
$$
using \cref{toobusy} for the first equality and the induction hypothesis for the second one.
\end{proof}

\begin{corollary}\label{oslo}
For $0 \leq r \leq n$, we have that
\begin{align}\label{threeblocks}
%
\:.
\end{align}
\end{corollary}

\begin{proof}
We just prove the first equality.
If $r < n$ then the expression on the left hand side is 0 by degree considerations
like in the first paragraph of the proof of \cref{toobusy}. If $r = n$ then the left hand side is equal to
$%
\:$,
and the conclusion follows from \cref{imanant}.
\end{proof}

The remaining relations to be established in this subsection
are more complicated. The guiding principle here is that 
relations in the nil-Hecke algebra can be ported to the nil-Brauer category
providing there are enough additional strings to eliminate the cup/cap term in the dot sliding relation \cref{rels7}.

\begin{lemma}\label{tryingtoohard}
For $0 \leq i \leq n+1$, we have that
\begin{align}
%
.
\end{align}
\end{lemma}

\begin{proof}
We first 
slide both sets of $i$ dots downwards past the crossing
using \cref{dotslide1,dotslide2}
to see that
\begin{align*}
%
\right).
$$
Now the lemma follows 
using also the identities
\begin{align*}
%
\quad\text{for $j \leq n+1$}.
\end{align*}
These are consequences of \cref{toobusy,curly}.
\end{proof}

\begin{lemma}\label{camy}
For $i,j \geq 0$ with $i+j \leq 2n+3$, we have that
\begin{align}
%
.
\end{align}
\end{lemma}

\begin{proof}
We assume that $i \leq j$, and proceed by induction on $j-i$.
The base case $j-i=0$ follows by \cref{tryingtoohard}.
For the induction step, suppose that $i < j$ and $i+j \leq 2n+3$.
By induction, we have that
$$
%
.
$$
Then we vertically compose on top with $e_{1,2n+2}=\frac{1}{2}q_{1,2n+2}(x_1,\dots,x_{2n+2})$,
using the centrality from \cref{bankrupt} to 
commute this down to the middle; it becomes
$e_{1,2} = x_1+x_2$ in the middle on the left hand side and $e_{1,0} = 0$ in the middle on the right hand side.
We deduce that
$$
%
=0.
$$
If $j-i=1$, the last two terms are the same as the first two terms, and the result follows on dividing by 2.
If $j-i > 1$ we use the induction hypothesis
to simplify the last two terms to obtain
$$
%
=0.
$$
The result follows.
\end{proof}

\begin{corollary}\label{sam}
For $\alpha \in \N^{n+1}$ and $1 \leq i \leq n$ such that
$\alpha_i + \alpha_{i+1} \leq 2n+1$, we have that
\begin{align}\label{samlabel}
%
\:,
\end{align}
where $\widehat\alpha := (\alpha_1,\dots,\alpha_{i-1},\alpha_{i+2},\dots,\alpha_{n+1}) \in \N^{n-1}$
and $s_i \alpha$ is the tuple obtained from 
$\alpha$ by permuting the $i$th and $(i+1)$th entries.
\end{corollary}

\begin{proof}
Let $\beta := (\alpha_1,\dots,\alpha_{i-1})$
and $\gamma := (\alpha_{i+2},\dots,\alpha_{n+1})$.
By \cref{camy}, we have that
\begin{align*}
%
\:.
\end{align*}
\end{proof}

\begin{corollary}\label{sammer}
For $\alpha \in \N^{n+1}$ and $1 \leq i \leq n$ such that
$\alpha_i =\alpha_{i+1} \leq n$, we have that
\begin{align}\label{sammerlabel}
%
\:,
\end{align}
where $\widehat\alpha := (\alpha_1,\dots,\alpha_{i-1},\alpha_{i+2},\dots,\alpha_{n+1}) \in \N^{n-1}$.
\end{corollary}

\begin{proof}
This follows from \cref{sam}.
\end{proof}

\begin{lemma}\label{lemmabelow}
The following relation holds for
any $n \geq 1$ and $0 \leq r \leq n$:
\begin{align}\label{lbelow}
%
\:,
\end{align}
interpreting the final term as 0 in case $r \leq 1$.
\end{lemma}

\begin{proof}
We proceed by induction on $n$. The result is trivial when 
$n = 1$. It is also clear when $r=0$
thanks to \cref{imanant}. Now 
suppose that $n \geq 1$ and $0 \leq r \leq n$, and consider
$$
%
\right).
$$
Here, we commuted the single dot upward through the thick string. In the summation, the second term is 0 always,
and the first term is 0 unless $a=0$. So the expression simplifies to give
\begin{equation}\label{hay}
%
\:.
\end{equation}
If $r=0$, we simplify this using \cref{imanant}, then \cref{toobusy}, then induction to obtain
\begin{align*}
%
\:,
\end{align*}
as required for the induction step.
Now suppose that $r \geq 1$ and consider
\cref{hay} again.
Letting $\alpha := \varpi_{r+1,n+1} - \varpi_{1,n+1}+\rho_{n+1} = (n,n,\dots)$,
we use \cref{sam} and induction
to simplify the first term:
$$
%
\:,
$$
which is the second term we need to prove the induction step. Turning our attention to the second term on the right hand side of \cref{hay},
it remains to show that
$$
%
$$
assuming $r \geq 1$.
By the induction hypothesis plus the identities
$e_{r,n} = e_{r,n-1} + e_{r-1,n-1} x_n$
then
$e_{r+1,n} + e_{r,n} x_{n+1}=e_{r+1,n+1}$, we have that
\begin{align*}
%
\:.
\end{align*}
\end{proof}

\begin{corollary}\label{corollarybelow}
The following relation holds for
any $n \geq 1$ and $0 \leq r \leq n$:
\begin{align}\label{cbelow}
%
\:.
\end{align}
\end{corollary}

\begin{proof}
Add a cap at the bottom of the relation
from \cref{lemmabelow}. The second term then disappears.
\end{proof}

\subsection{Recurrence relation for idempotents}
\cref{imanant} obviously implies that 
\begin{equation}\label{endef}
\e_n :=
%
\end{equation}
is a homogeneous idempotent 
for each $n \geq 0$.
For example:
\begin{align*}
\e_0 &= 1_\one,&\e_1
&= 
%
\:\;.
\end{align*}
These are likely already familiar expressions, since the same diagrams are often used to 
represent distinguished primitive idempotents in the nil-Hecke algebra.

In the remainder of the section, 
we are going to show that
the idempotents $\e_n\:(n \geq 0)$
give a full set of primitive homogeneous idempotents in $\cNB_t$.
The first step, accomplished in this subsection, is to decompose
$B \star \e_n$ as a sum of mutually orthogonal
conjugates of $\e_{n+1}$ and $\e_{n-1}$.
We begin by introducing two more families of endomorphisms of $B^{\star (n+1)}$:
for $0 \leq r \leq n$ let
\begin{align}
\e_{r,n} &:= 
(-1)^r \:
\:\right)\:,
\label{fdef}
\end{align}
Recalling the convention that
the elementary symmetric function
$e_r = 0$ for $r < 0$ and, of course, $e_0=1$,
we have that
\begin{align}
\e_{0,n}
&=\e_{n+1},&
\f_{0,n} &= 0.
\end{align}
By \cref{lemmabelow,corollarybelow}, the definitions \cref{edef,fdef} can be written equivalently as
\begin{align}
\e_{r,n} &=(-1)^r \:%
,
\end{align}
where we interpret terms involving the undefined
symbols
$\varpi_{r-1}$ for $r = 0$
and $\varpi_{r-2}$ for $r=0$ or $1$ as 0.

\begin{example}
If $n=0$ then $\e_{0,0} = \;%
\:.
\end{align*}
\end{example}

\begin{lemma}\label{hug}
For $n \geq 0$, we have that
$B \star \e_n = \sum_{r=0}^n (\e_{r,n} + \f_{r,n})$.
\end{lemma}

\begin{proof}
For this calculation, 
it is convenient to drop the $\rho$
from the top of the diagrams, so we set
\begin{align*}
\mathring{\e}_n &:=
%
\:.
\end{align*}
We in fact show that $B \star \mathring{\e}_n = \sum_{r=0}^n (\mathring{\e}_{r,n}+\mathring{\f}_{r,n})$.
The first step is the same as in the proof of \cite[Lem.~2.13]{KLMS}:
\begin{align*}
B \star \mathring{\e}_n &= 
%
.
\end{align*}
The last two terms in this expression are equal to
$\mathring{\e}_{n,n}+\mathring{\f}_{n,n}$.
It remains to show that the first term is equal to
$\sum_{r=0}^{n-1} (\mathring{\e}_{r,n}+\mathring{\f}_{r,n})$:
\begin{align*}
(-1)^{n-1}\:
&%
\:.
\end{align*}
The first summation gives the remaining terms
$\sum_{r=0}^{n-1}(\mathring{\e}_{r,n}+\mathring{\f}_{r,n})$
that we want,
and the second summation is 0 thanks to \cref{corollarybelow,oslo}.
\end{proof}

Now we introduce several more families of morphisms in $\cNB_t$ for $0 \leq r\leq n$
and $1 \leq s \leq n$:
\begin{align}\label{uv}
\u_{r,n} &:=(-1)^r
%
,
\end{align}
again interpreting the undefined term involving $\varpi_{s-2}$
when $s=1$ as 0.
Note that $\u_{0,n} = \vv_{0,n} = \e_{n+1}$
thanks to \cref{imanant}, hence, $\w_{0,n} = 0$.
The same corollary also implies easily that 
$\e_{n+1}\circ \u_{r,n} = \u_{r,n}$,
$\vv_{r,n} \circ \e_{n+1} = \vv_{r,n}$,
$\e_{n-1} \circ \x_{s,n} = \x_{s,n}$ and
$\y_{s,n} \circ \e_{n-1} = \y_{s,n}$.

\begin{lemma}\label{id0}
For $0 \leq r \leq n$ and $1 \leq s \leq n$, we have that $\vv_{r,n} \circ \u_{r,n} = \e_{r,n}$
 and $\y_{s,n} \circ \x_{s,n} = \f_{s,n}$.
\end{lemma}

\begin{proof}
This follows from the definitions just given, 
using \cref{imanant} and
the alternative forms of the definitions of $\e_{r,n}$
and $\f_{s,n}$
from \cref{edefalt,fdefalt2}.
\end{proof}

\begin{lemma}\label{id1}
For $0 \leq r,s \leq n$, we have that
\begin{equation}
\u_{r,n} \circ \vv_{s,n} = 
\begin{dcases}
-\begin{tikzpicture}[anchorbase]
\draw (-.1,.4) to (-.1,-.4);
\draw[ultra thick] (.3,.4) to (.3,-.4);
\closeddot{-.1,0};\node at (-.28,0) {$\scriptstyle r$};
\node at (.31,-.52) {$\stringlabel{n}$};
\node at (.31,.5) {$\stringlabel{\phantom n}$};
\end{tikzpicture}\circ \f_{r,n}&\text{if $s=0 < r$ and $n\not\equiv t\pmod{2}$}\\
\delta_{r,s} \e_{n+1}&\text{otherwise.}
\end{dcases}
\end{equation}
\end{lemma}

\begin{proof}
This is clear for $n=0$ so assume $n \geq 1$.
By the definitions and \cref{imanant}, we have that
\begin{align*}
\u_{r,n} \circ \vv_{s,n} &=
(-1)^r
\begin{tikzpicture}[anchorbase,scale=1.2]
\draw[ultra thick] (-.2,-.7) to[out=45,in=-45,looseness=1.1] (0,.3) to (0,.73);
\draw (.2,-.7) to [out=135,in=-135,looseness=1.1] (-.01,.29);
\smallercross[8]{(.15,-.02)};
\closeddot{0.14,-.19};\node at (0.3,-.2) {$\scriptstyle\rho$};
\smallercross[-25]{(.21,-.29)};
\closeddot{0.1,.17};\node at (0.5,.2) {$\scriptstyle\varpi_r+\rho$};
\smallercross{(0,.4)};
\closeddot{0.01,.6};\node at (0.15,.6) {$\scriptstyle\rho$};
\closeddot{-0.15,-0.1};\node at (-0.5,-0.1) {$\scriptstyle n-s$};
\node at (0,.85) {$\stringlabel{n+1}$};
\node at (-.2,-.82) {$\stringlabel{n}$};
\end{tikzpicture}=
(-1)^r
\begin{tikzpicture}[anchorbase,scale=1.2]
\draw[ultra thick] (-.2,-.7) to[out=45,in=-45,looseness=1.1] (0,.3) to (0,.73);
\draw (.2,-.7) to [out=135,in=-135,looseness=1.1] (-.01,.29);
\smallercross[-20]{(.19,-.2)};
\closeddot{0.15,.05};\node at (0.52,.1) {$\scriptstyle\varpi_r+\rho$};
\smallercross{(0,.4)};
\closeddot{0.01,.6};\node at (0.15,.6) {$\scriptstyle\rho$};
\closeddot{-0.15,-0.1};\node at (-0.5,-0.1) {$\scriptstyle n-s$};
\node at (0,.85) {$\stringlabel{n+1}$};
\node at (-.2,-.82) {$\stringlabel{n}$};
\end{tikzpicture}
=
(-1)^r \begin{tikzpicture}[anchorbase,scale=1.2]
\draw[ultra thick] (0,.73) to (0,-.7);
\closeddot{0.005,-.1};
\closeddot{0.005,.32};
\smallercross{(0,.1)};
\smallercross{(0,-.3)};
\node at (0,-.82) {$\stringlabel{n+1}$};
\node at (0.22,-.1) {$\scriptstyle{\alpha}$};
\node at (0.22,.32) {$\scriptstyle{\rho}$};
\node at (0,.85){$\stringlabel{\phantom{n}}$};
\end{tikzpicture}
\end{align*}
where
$\alpha = (n-s,n,n-1,\dots,n-r+1,n-r-1,\dots,1,0) \in \N^{n+1}$.
If $s=r$ then $\alpha$ is a rearrangement of $\rho_{n+1}$,
so this is equal to $\e_{n+1}$ thanks to \cref{sam}.
If $0 < s \neq r$
then $\alpha$
has two entries equal to $n-s < n$, so this is 0
by \cref{sam,sammer}.
Finally if $0=s \neq r$
then $\alpha = (n,n,n-1,\dots,n-r+1,n-r-1,\dots,1,0)$
and \cref{sammer} gives the exceptional formula in this case,
referring to \cref{fdefalt} to see the appropriate form of
$\f_{r,n}$.
\end{proof}

\begin{corollary}\label{id2}
For $0 \leq r,s \leq n$, we have that
\begin{equation}
\e_{r,n} \circ \e_{s,n} = 
\begin{dcases}
-\f_{r,n}&\text{if $s=0 < r$ and $n\not\equiv t\pmod{2}$}\\
\delta_{r,s} \e_{r,n}&\text{otherwise.}
\end{dcases}
\end{equation}
\end{corollary}

\begin{proof}
By \cref{id0}, 
we have that $\e_{r,n} \circ \e_{s,n}
= \vv_{r,n} \circ \u_{r,n} \circ \vv_{s,n} \circ \u_{s,n}$.
Except in the case $s=0< r$ and $n\not\equiv t\pmod{2}$,
we have that $\u_{r,n} \circ \vv_{s,n} = \delta_{r,s} 
\e_{n+1}$ by \cref{id1}, and 
$\e_{n+1} \circ \u_{r,n} = \u_{r,n}$ by \cref{imanant}.
The conclusion then follows
using that $\vv_{r,n} \circ \u_{r,n} = \e_{r,n}$ once again.
Suppose from now on that $s=0 < r$ and $n\not\equiv t\pmod{2}$.
Then, using the form of $\f_{r,n}$
from \cref{fdefalt2}, 
\cref{id1} gives instead that
\begin{align*}
\e_{r,n} \circ \e_{0,n} 
&=\vv_{r,n} \circ (\u_{r,n} \circ \vv_{0,n}) \circ \u_{0,n}=
(-1)^{r}\!\!\!\!\!\!\!
\begin{tikzpicture}[anchorbase,scale=1.4]
\draw[ultra thick] (0,-1.2) to (0,1.2) to (.2,1.4);
\draw (-.2,1.4) to (-0.01,1.2);
\draw (0,-.95) to[out=135,in=-90,looseness=1] (-.3,-.5) to[out=90,in=135,looseness=1] (-.02,-.2);
\draw (0,.95) to[out=-135,in=90,looseness=1] (-.3,.5) to[out=-90,in=-135,looseness=1] (-.02,.2);
\smallestcross{(0,1.08)};
\smallestcross{(0,-1.08)};
\closeddot{-.29,-.6};\node at (-.45,-.6) {$\scriptstyle n$};
\closeddot{-.29,.6};\node at (-.45,.6) {$\scriptstyle n$};
\closeddot{0,.6};\node at (.15,.6) {$\scriptstyle\rho$};
\closeddot{0,-.6};\node at (.15,-.6) {$\scriptstyle\rho$};
\closeddot{-.1,1.3};
\closeddot{.1,1.3};
\closeddot{0.01,0};
\smallestcross{(0,0.4)};
\smallestcross{(0,-.4)};
\node at (.46,0) {$\scriptstyle \varpi_{r-1}+\rho$};
\node at (-.45,1.3) {$\scriptstyle n-r$};
\node at (.28,1.28) {$\scriptstyle \rho$};
\node at (0,-1.3) {$\stringlabel{n+1}$};
\node at (0.2,1.5) {$\stringlabel{n}$};
\end{tikzpicture}\!
=(-1)^{r}\!\!\!\!\!\!
\begin{tikzpicture}[anchorbase,scale=1.4]
\draw[ultra thick] (0,-1.2) to (0,1.1) to (.3,1.4);
\draw (-.3,1.4) to (-0.01,1.1);
\draw (.2,-1.2) to[out=90,in=-45] (0,-.95) to[out=135,in=-90,looseness=1] (-.3,-.5) to[out=90,in=135,looseness=1] (-.02,-.2);
\draw (.1,1.2) to[out=-45,in=90] (.2,1.1) to [out=-90,in=45] (0,.85) to[out=-135,in=90,looseness=1] (-.3,.4) to[out=-90,in=-135,looseness=1] (-.02,.1);
\smallestcross{(0,.65)};
\smallestcross{(0,-.75)};
\closeddot{-.3,-.55};\node at (-.45,-.55) {$\scriptstyle n$};
\closeddot{-.3,.45};\node at (-.45,.45) {$\scriptstyle n$};
\closeddot{0,.45};\node at (.15,.45) {$\scriptstyle\rho$};
\closeddot{0,-.55};\node at (.15,-.55) {$\scriptstyle\rho$};
\closeddot{-.2,1.3};
\closeddot{.2,1.3};
\closeddot{0.01,-.05};
\smallestcross{(0,0.25)};
\smallestcross{(0,-.35)};
\node at (.46,-.05) {$\scriptstyle \varpi_{r-1}+\rho$};
\node at (-.55,1.3) {$\scriptstyle n-r$};
\node at (.38,1.28) {$\scriptstyle \rho$};
\node at (0,-1.3) {$\stringlabel{n}$};
\node at (0.3,1.5) {$\stringlabel{n}$};
\end{tikzpicture}\!\!
=(-1)^{r}\!\!\!\!\!\!\!
\begin{tikzpicture}[anchorbase,scale=1.4]
\draw[ultra thick] (0,-1.2) to (0,1.1) to (.3,1.4);
\draw (-.3,1.4) to (-0.01,1.1);
\draw (.4,-1.2) to[out=90,in=-90,looseness=1] (-.3,-.5) to[out=90,in=180,looseness=1] (0,-.2) to[out=0,in=90,looseness=1](.2,-1.2);
\draw (.15,1.25) to[out=-45,in=90] (.3,1) to [out=-90,in=90] (-.3,.4) to[out=-90,in=180,looseness=1] (0,.1) to[out=0,in=-45,looseness=1] (0.05,1.15);
\smallestcross{(0,.45)};
\smallestcross{(0,-.55)};
\closeddot{-.3,-.55};\node at (-.45,-.55) {$\scriptstyle n$};
\closeddot{-.3,.45};\node at (-.45,.45) {$\scriptstyle n$};
\closeddot{-.2,1.3};
\closeddot{.2,1.3};
\closeddot{0.01,-.05};
\node at (.46,-.05) {$\scriptstyle \varpi_{r-1}+\rho$};
\node at (-.55,1.3) {$\scriptstyle n-r$};
\node at (.38,1.28) {$\scriptstyle \rho$};
\node at (0,-1.35) {$\stringlabel{n-1}$};
\node at (0.3,1.5) {$\stringlabel{n}$};
\end{tikzpicture}\!.
\end{align*}
It remains to apply \cref{curly} to see that this is equal to
$-\f_{r,n}$; for this \cref{fdefalt} is most convenient.
\end{proof}

\begin{lemma}\label{id2b}
Assume that $n\equiv t\pmod{2}$. For $1 \leq r,s \leq n$,
we have that
$\x_{r,n}\circ\y_{s,n} = \delta_{r,s} \e_{n-1}$.
\end{lemma}

\begin{proof}
When $n=t=1$ this follows immediately from the first relation from \cref{rels2}. Now suppose that $n \geq 2$.
Since $\x_{r,n}$ is a sum of two terms (the second being 0 in case $r=1$), so too is $\x_{r,n}\circ \y_{s,n}$. We compute the two terms separately. The first term is
\begin{align*}
(-1)^{r-1}\!\!\!\begin{tikzpicture}[anchorbase]
\draw[ultra thick] (0,-1.2) to (0,.4);
\draw (-.02,-.25) to [out=140,in=90,looseness=1.5] (-.5,-.6) to [out=-90,in=-140,looseness=1.5] (-.02,-.95);
\closeddot{0.01,.3};
\node at (.22,.3) {$\scriptstyle \rho$};
\closeddot{0.01,-.1};
\cross{(0,-.4)};
\cross{(0,.1)};
\node at (.62,-.1) {$\scriptstyle \varpi_{r-1}+\rho$};
\node at (0,.53) {$\stringlabel{n-1}$};
\closeddot{0.01,-.6};
\closeddot{-.47,-.6};
\cross{(0,-.8)};
\node at (.23,-.6) {$\scriptstyle \rho$};
\node at (-.87,-.6) {$\scriptstyle n-s$};
\node at (0,-1.33) {$\stringlabel{n-1}$};
\end{tikzpicture}
\!&=
(-1)^{r-1}\!\!\!\!\begin{tikzpicture}[anchorbase]
\draw[ultra thick] (0,-.9) to (0,.7);
\draw (-.02,-.15) to [out=140,in=90,looseness=1.5] (-.5,-.4) to [out=-90,in=-140,looseness=1.5] (-.02,-.65);
\closeddot{0.01,.55};
\node at (.22,.55) {$\scriptstyle \rho$};
\closeddot{0.01,.05};
\cross{(0,-.4)};
\cross{(0,.3)};
\node at (.62,.05) {$\scriptstyle \varpi_{r-1}+\rho$};
\node at (0,.83) {$\stringlabel{n-1}$};
\closeddot{-.47,-.4};
\node at (-.87,-.4) {$\scriptstyle n-s$};
\node at (0,-1.03) {$\stringlabel{n-1}$};
\end{tikzpicture}
=
(-1)^{r-1}\!\!\!\!\!\begin{tikzpicture}[anchorbase]
\draw[ultra thick] (0,-.9) to (0,.7);
\draw (.25,-.9) to [out=90,in=0] (-.1,-.3) to [out=180,in=90] (-.3,-.5) to [out=-90,in=180] (-.1,-.7) to [out=0,in=-90] (.2,-.4) to [out=90,in=-45] (0,-.1);
\closeddot{0.01,.55};
\node at (.22,.55) {$\scriptstyle \rho$};
\closeddot{0.01,.05};
\cross{(0,-.5)};
\cross{(0,.3)};
\node at (.62,.05) {$\scriptstyle \varpi_{r-1}+\rho$};
\node at (0,.83) {$\stringlabel{n-1}$};
\closeddot{-.3,-.5};
\node at (-.7,-.5) {$\scriptstyle n-s$};
\node at (-0.05,-1.03) {$\stringlabel{n-2}$};
\end{tikzpicture}
\!=(-1)^{r-1}\delta_{s,1}\!\begin{tikzpicture}[anchorbase]
\draw[ultra thick] (0,-.8) to (0,.8);
\closeddot{0.01,.4};
\node at (.22,.4) {$\scriptstyle \rho$};
\closeddot{0.01,-.15};
\cross{(0,.15)};
\cross{(0,-.4)};
\node at (.62,-.15) {$\scriptstyle \varpi_{r-1}+\rho$};
\node at (0,-.93) {$\stringlabel{n-1}$};
\end{tikzpicture},
\end{align*}
where we used \cref{imanant} for the first equality and \cref{curly} for the last one.
If $r=1$ (when we already know that the second term is 0)
this is $\delta_{s,1} \e_{n-1}$ by \cref{imanant}, and we are done.
Assuming from now on that $r \geq 2$, the second term is
\begin{align*}
(-1)^{r}\!\begin{tikzpicture}[anchorbase,scale=1.2]
\draw[ultra thick] (0,-1.3) to (0,.5);
\draw (-.02,-1.15) to [out=-155,in=-90,looseness=1.25] (-.9,-.65) to [out=90,in=-160,looseness=1.25] (-.02,0);
\draw (-.02,-.3) to[out=145,in=90,looseness=1.5] 
(-.65,-.65) to [out=-90,in=90] (.25,-1.3);
\closeddot{0.01,.35};
\node at (.2,.35) {$\scriptstyle \rho$};
\closeddot{0.01,-.15};
\smallercross{(0,-.45)};
\smallercross{(0,.15)};
\node at (.55,-.15) {$\scriptstyle \varpi_{r-2}+\rho$};
\node at (0,.61) {$\stringlabel{n-1}$};
\closeddot{0.01,-.65};
\closeddot{-.65,-.65};
\closeddot{-.9,-.65};
\smallercross{(0,-.85)};
\node at (.2,-.65) {$\scriptstyle \rho$};
\node at (-.35,-.65) {$\scriptstyle n-1$};
\node at (-1.2,-.65) {$\scriptstyle n-s$};
\node at (-0.05,-1.41) {$\stringlabel{n-2}$};
\end{tikzpicture}
&=
(-1)^{r}\!\begin{tikzpicture}[anchorbase,scale=1.2]
\draw[ultra thick] (0,-1.3) to (0,.5);
\draw (-.02,-1.15) to [out=-155,in=-90,looseness=1.25] (-.9,-.65) to [out=90,in=-160,looseness=1.25] (-.02,0);
\draw (-.02,-.4) to[out=145,in=90,looseness=1.5] 
(-.65,-.65) to [out=-90,in=90] (.25,-1.2) to (.25,-1.3);
\closeddot{0.01,.35};
\node at (.2,.35) {$\scriptstyle \rho$};
\closeddot{0.01,-.2};
\smallercross{(0,.15)};
\node at (.55,-.2) {$\scriptstyle \varpi_{r-2}+\rho$};
\node at (0,.61) {$\stringlabel{n-1}$};
\closeddot{-.65,-.65};
\closeddot{-.9,-.65};
\smallercross{(0,-.65)};
\node at (-.35,-.65) {$\scriptstyle n-1$};
\node at (-1.2,-.65) {$\scriptstyle n-s$};
\node at (-0.05,-1.41) {$\stringlabel{n-2}$};
\end{tikzpicture}
=(-1)^{r}\!\begin{tikzpicture}[anchorbase,scale=1.2]
\draw[ultra thick] (0,-1.3) to (0,.5);
\draw (-.02,-1.15) to [out=-155,in=-90,looseness=1.25] (-.9,-.65) to [out=90,in=-160,looseness=1.25] (-.02,0);
\draw (.5,-1.3) to[out=90,in=-90,looseness=1] 
(-.65,-.65) to [out=90,in=180,looseness=1] (-.2,-.3) to [out=0,in=90,looseness=1]  (.25,-1.3);
\closeddot{0.01,.35};
\node at (.2,.35) {$\scriptstyle \rho$};
\closeddot{0.01,-.2};
\smallercross{(0,.15)};
\node at (.55,-.2) {$\scriptstyle \varpi_{r-2}+\rho$};
\node at (0,.61) {$\stringlabel{n-1}$};
\closeddot{-.65,-.65};
\closeddot{-.9,-.65};
\smallercross{(0,-.65)};
\node at (-.35,-.65) {$\scriptstyle n-1$};
\node at (-1.2,-.65) {$\scriptstyle n-s$};
\node at (-0.05,-1.41) {$\stringlabel{n-3}$};
\end{tikzpicture}\\
&\!\!\stackrel{\cref{curlylabel}}{=}
(-1)^{r}\!\begin{tikzpicture}[anchorbase,scale=1.2]
\draw[ultra thick] (0,-1.3) to (0,.5);
\draw (-.02,-1) to [out=-140,in=-90,looseness=1.25] (-.5,-.65) to [out=90,in=-160,looseness=1.25] (-.02,0);
\draw (.3,-1.3) to[out=90,in=30,looseness=1] 
(0,-.4);
\closeddot{0.01,.35};
\node at (.2,.35) {$\scriptstyle \rho$};
\closeddot{0.01,-.2};
\smallercross{(0,.15)};
\node at (.55,-.2) {$\scriptstyle \varpi_{r-2}+\rho$};
\node at (0,.61) {$\stringlabel{n-1}$};
\closeddot{-.5,-.65};
\smallercross{(0,-.65)};
\node at (-.9,-.65) {$\scriptstyle n-s$};
\node at (-0.05,-1.41) {$\stringlabel{n-2}$};
\end{tikzpicture}
=
(-1)^{r}\!\begin{tikzpicture}[anchorbase,scale=1.2]
\draw[ultra thick] (0,-1.3) to (0,.5);
\draw (.3,-1.3) to [out=90,in=-90,looseness=1.25] (-.5,-.55) to [out=90,in=-160,looseness=1.25] (-.02,0);
\closeddot{0.01,.35};
\node at (.2,.35) {$\scriptstyle \rho$};
\closeddot{0.01,-.2};
\smallercross{(0,.15)};
\node at (.55,-.2) {$\scriptstyle \varpi_{r-2}+\rho$};
\node at (0,.61) {$\stringlabel{n-1}$};
\closeddot{-.5,-.55};
\smallercross{(0,-.55)};
\node at (-.9,-.55) {$\scriptstyle n-s$};
\node at (-0.05,-1.41) {$\stringlabel{n-2}$};
\end{tikzpicture}
=(-1)^{r}\!\begin{tikzpicture}[anchorbase,scale=1.2]
\draw[ultra thick] (0,-1.15) to (0,.65);
\closeddot{0.01,.4};
\node at (.2,.4) {$\scriptstyle \rho$};
\smallercross{(0,0)};
\node at (0,.76) {$\stringlabel{\phantom{n-1}}$};
\closeddot{0,-.4};
\smallercross{(0,-.8)};
\node at (.2,-.4) {$\scriptstyle \alpha$};
\node at (-0.05,-1.26) {$\stringlabel{n-1}$};
\end{tikzpicture}
\end{align*}
where $\alpha = (n-s,n-2,\dots,n-r+1,n-r-1,\dots,1,0) \in \N^{n-1}$.
If $s=1$ this cancels with the first term to give 0, and we are done.
Assuming from now on that $s \geq 2$,
the first term is 0, and it just remains to 
apply \cref{sam,sammer} to rewrite the second term, noting that $n\equiv t\pmod{2}$ so the first term on the right hand side of \cref{samlabel} is 0, as is the right hand side of \cref{sammerlabel}.
We get 0 if $r \neq s$ and, after one more application of \cref{imanant}, we get $\e_{n-1}$ if $r=s$, as claimed.
\end{proof}

\begin{corollary}\label{id2bcor}
Assume that $n\equiv t\pmod{2}$. 
For $1 \leq r,s \leq n$,
we have that
$\f_{r,n}\circ\f_{s,n} = \delta_{r,s} \f_{r,n}$.
\end{corollary}

\begin{proof}
This follows by \cref{id0,id2b}.
\end{proof}

\begin{lemma}\label{id2c}
Assume that $n\equiv t\pmod{2}$. For $0 \leq r\leq n$ and $1\leq s \leq n$, we have that
$\u_{r,n} \circ \y_{s,n} = \x_{s,n} \circ \vv_{r,n} = 0$.
\end{lemma}

\begin{proof}
We first consider $\x_{s,n} \circ \vv_{r,n}$.
Since $\x_{s,n}$ is a sum of two terms, so 
too is $\x_{s,n} \circ \vv_{r,n}$. We show that both of these terms are 0. The first term is
\begin{align*}
(-1)^{s-1}\!\!\begin{tikzpicture}[anchorbase]
\draw[ultra thick] (0,-1.2) to (0,.4);
\draw (-.02,-.25) to [out=140,in=90,looseness=1.5] (-.5,-.4) to [out=-90,in=155,looseness=1.5] (-.02,-.82);
\closeddot{0.01,.3};
\node at (.22,.3) {$\scriptstyle \rho$};
\closeddot{0.01,-.1};
\cross{(0,-.4)};
\cross{(0,.1)};
\node at (.62,-.1) {$\scriptstyle \varpi_{s-1}+\rho$};
\node at (0,.53) {$\stringlabel{n-1}$};
\closeddot{0.01,-.65};
\closeddot{-.5,-.4};
\cross{(0,-1)};
\node at (.23,-.65) {$\scriptstyle \rho$};
\node at (-.87,-.4) {$\scriptstyle n-r$};
\node at (0,-1.33) {$\stringlabel{n+1}$};
\end{tikzpicture}\!\!\!\!\!
&=
(-1)^{s-1}\!\!\!\!\begin{tikzpicture}[anchorbase]
\draw[ultra thick] (0,-1.2) to (0,.4);
\draw (-.02,-.25) to [out=140,in=90,looseness=1.25] (-.5,-.5) to [out=-90,in=120,looseness=1.25] (.25,-1.2);
\closeddot{0.01,.3};
\node at (.22,.3) {$\scriptstyle \rho$};
\closeddot{0.01,-.1};
\cross{(0,-.4)};
\cross{(0,.1)};
\node at (.62,-.1) {$\scriptstyle \varpi_{s-1}+\rho$};
\node at (0,.53) {$\stringlabel{n-1}$};
\closeddot{0.01,-.6};
\closeddot{-.5,-.5};
\cross{(0,-.8)};
\node at (.23,-.6) {$\scriptstyle \rho$};
\node at (-.87,-.5) {$\scriptstyle n-r$};
\node at (0,-1.33) {$\stringlabel{n}$};
\end{tikzpicture}\!\!\!
=(-1)^{s-1}\!\!\!\!\begin{tikzpicture}[anchorbase]
\draw[ultra thick] (0,-1.2) to (0,.4);
\draw (-.02,-.45) to [out=140,in=90,looseness=1.25] (-.5,-.55) to [out=-90,in=120,looseness=1.25] (.25,-1.2);
\closeddot{0.01,.25};
\node at (.22,.25) {$\scriptstyle \rho$};
\closeddot{0.01,-.25};
\cross{(0,0)};
\node at (.62,-.25) {$\scriptstyle \varpi_{s-1}+\rho$};
\node at (0,.53) {$\stringlabel{n-1}$};
\closeddot{-.5,-.55};
\cross{(0,-.7)};
\node at (-.87,-.55) {$\scriptstyle n-r$};
\node at (0,-1.33) {$\stringlabel{n}$};
\end{tikzpicture}\!
=(-1)^{s-1}\!\!\!\!\begin{tikzpicture}[anchorbase]
\draw[ultra thick] (0,-1.2) to (0,.4);
\draw (.25,-1.2) to [out=90,in=0,looseness=1] (-.2,-.35) to [out=180,in=90,looseness=1] (-.5,-.65) to [out=-90,in=90,looseness=1.5] (.5,-1.2);
\closeddot{0.01,.25};
\node at (.22,.25) {$\scriptstyle \rho$};
\closeddot{0.01,-.25};
\cross{(0,0)};
\node at (.62,-.25) {$\scriptstyle \varpi_{s-1}+\rho$};
\node at (0,.53) {$\stringlabel{n-1}$};
\closeddot{-.49,-.65};
\cross{(0,-.7)};
\node at (-.87,-.65) {$\scriptstyle n-r$};
\node at (-0.05,-1.33) {$\stringlabel{n-1}$};
\end{tikzpicture}\!\!.
\end{align*}
This is 0 by \cref{curly} since $n-1\not\equiv t\pmod{2}$.
The second term is 0 automatically if $s=1$, so we are done in this case.
When $s \geq 2$, the second term equals
\begin{align*}
(-1)^{s}\!\begin{tikzpicture}[anchorbase]
\draw[ultra thick] (0,-1.6) to (0,.7);
\draw (-.02,.1) to [out=-140,in=90,looseness=1] (-1.1,-.55) to [out=-90,in=140,looseness=1] (-.02,-1.2);
\draw (-.02,-.32) to [out=140,in=90,looseness=1.5] (-.4,-.6) to [out=-90,in=140,looseness=1.5] (-.02,-.97);
\closeddot{0.01,.5};
\node at (.22,.5) {$\scriptstyle \rho$};
\closeddot{0.01,-.15};
\cross{(0,-.5)};
\cross{(0,.3)};
\node at (.62,-.15) {$\scriptstyle \varpi_{s-1}+\rho$};
\node at (0,.83) {$\stringlabel{n-1}$};
\closeddot{0.01,-.75};
\closeddot{-.4,-.6};
\closeddot{-1.08,-.45};
\cross{(0,-1.4)};
\node at (.23,-.75) {$\scriptstyle \rho$};
\node at (-.77,-.61) {$\scriptstyle n-1$};
\node at (-1.45,-.44) {$\scriptstyle n-r$};
\node at (0,-1.73) {$\stringlabel{n+1}$};
\end{tikzpicture}\!
&=
(-1)^{s}\!\begin{tikzpicture}[anchorbase]
\draw[ultra thick] (0,-1.6) to (0,.7);
\draw (-.02,.1) to [out=-140,in=90,looseness=1] (-1.1,-.65) 
to [out=-90,in=90] (-1.1,-.75) to [out=-90,in=130,looseness=.75] (.5,-1.6);
\draw (-.02,-.32) to [out=140,in=90,looseness=1.5] (-.4,-.7) to [out=-90,in=90,looseness=1.75] (.25,-1.6);
\closeddot{0.01,.5};
\node at (.22,.5) {$\scriptstyle \rho$};
\closeddot{0.01,-.15};
\cross{(0,-.5)};
\cross{(0,.3)};
\node at (.62,-.15) {$\scriptstyle \varpi_{s-1}+\rho$};
\node at (0,.83) {$\stringlabel{n-1}$};
\closeddot{0.01,-.725};
\closeddot{-.4,-.725};
\closeddot{-1.05,-.45};
\cross{(0,-.95)};
\node at (.23,-.725) {$\scriptstyle \rho$};
\node at (-.77,-.725) {$\scriptstyle n-1$};
\node at (-1.43,-.44) {$\scriptstyle n-r$};
\node at (-0.05,-1.73) {$\stringlabel{n-1}$};
\end{tikzpicture}\!
=
(-1)^{s}\!\begin{tikzpicture}[anchorbase]
\draw[ultra thick] (0,-1.6) to (0,.7);
\draw (-.02,.1) to [out=-140,in=90,looseness=1] (-1.1,-.65) 
to [out=-90,in=90] (-1.1,-.85) to [out=-90,in=130,looseness=.75] (.75,-1.6);
\draw (.25,-1.6) to [out=90,in=0,looseness=1.5] (-0.1,-.4) to [out=180,in=90] (-.4,-.7) to [out=-90,in=90,looseness=1.75] (.5,-1.6);
\closeddot{0.01,.5};
\node at (.22,.5) {$\scriptstyle \rho$};
\closeddot{0.01,-.15};
\cross{(0,.3)};
\node at (.62,-.15) {$\scriptstyle \varpi_{s-1}+\rho$};
\node at (0,.83) {$\stringlabel{n-1}$};
\cross{(0,-.725)};
\closeddot{-.4,-.725};
\closeddot{-1.05,-.45};
\node at (-.77,-.725) {$\scriptstyle n-1$};
\node at (-1.43,-.44) {$\scriptstyle n-r$};
\node at (-0.05,-1.73) {$\stringlabel{n-2}$};
\end{tikzpicture}\\
&\!\!\stackrel{\cref{curlylabel}}{=}
(-1)^{s}\!\begin{tikzpicture}[anchorbase]
\draw[ultra thick] (0,-1.6) to (0,.7);
\draw (-.02,.1) to [out=-140,in=90,looseness=1] (-.8,-.65) 
to [out=-90,in=90] (-.8,-.75) to [out=-90,in=130,looseness=.75] (.6,-1.6);
\draw (.35,-1.6) to [out=90,in=40,looseness=1.5] (0,-.5);
\closeddot{0.01,.5};
\node at (.22,.5) {$\scriptstyle \rho$};
\closeddot{0.01,-.18};
\cross{(0,.3)};
\node at (.62,-.18) {$\scriptstyle \varpi_{s-1}+\rho$};
\node at (0,.83) {$\stringlabel{n-1}$};
\cross{(0,-.725)};
\closeddot{-.75,-.45};
\node at (-1.15,-.44) {$\scriptstyle n-r$};
\node at (-0.05,-1.73) {$\stringlabel{n-1}$};
\end{tikzpicture}
=
(-1)^{s}\!\begin{tikzpicture}[anchorbase]
\draw[ultra thick] (0,-1.6) to (0,.7);
\draw (-.02,.1) to [out=-140,in=90,looseness=1] (-.8,-.65) 
to [out=-90,in=90] (-.8,-.75) to [out=-90,in=130,looseness=.75] (.6,-1.6);
\draw (.35,-1.6) to [out=90,in=0,looseness=1] (0,-.8)
to [out=180,in=90,looseness=1] (-.35,-1.6);
\closeddot{0.01,.5};
\node at (.22,.5) {$\scriptstyle \rho$};
\closeddot{0.01,-.18};
\cross{(0,.3)};
\node at (.62,-.18) {$\scriptstyle \varpi_{s-1}+\rho$};
\node at (0,.83) {$\stringlabel{n-1}$};
\cross{(0,-.5)};
\closeddot{-.75,-.45};
\node at (-1.15,-.44) {$\scriptstyle n-r$};
\node at (0,-1.73) {$\stringlabel{n-2}$};
\end{tikzpicture},
\end{align*}
which is 0 by the second relation from \cref{rels6}.

Now consider $\u_{r,n} \circ \y_{s,n}$
for $0 \leq r \leq n$ and $1 \leq s \leq n$.
For notational convenience, we in fact show that 
$\mathring{\u}_{r,n} \circ \y_{s,n} = 0$, where
$\mathring{\u}_{r,n} := 
(-1)^r
\begin{tikzpicture}[baseline=-.75mm,scale=1]
\draw[ultra thick] (0,.3) to (0,0) to (.3,-.4);
\draw (-.01,0) to (-.3,-.4);
\cross{(0,.13)};
\closeddot{0.09,-.12};\node at (.55,-.05) {$\scriptstyle\varpi_r+\rho$};
\cross[40]{(-.02,-.33)};
\node at (0,.42) {$\stringlabel{n+1}$};
\node at (.4,-.52) {$\stringlabel{n}$};
\end{tikzpicture}$.
Applying \cref{imanant} as usual, we have that
\begin{align*}
\mathring{\u}_{r,n} \circ \y_{s,n} &=
(-1)^r
\begin{tikzpicture}[anchorbase,scale=1.2]
\draw[ultra thick] (0,.3) to (0,-.7);
\draw (-.02,0) to[out=-140,in=90,looseness=1] (-.4,-.45) to [out=-90,in=-120,looseness=1.5] (-.02,-.55);
\smallercross{(0,.11)};
\closeddot{0,-.2};\node at (.42,-.2)  {$\scriptstyle \varpi_r+\rho$};
\closeddot{-.4,-.45};\node at (-.75,-.45) {$\scriptstyle n-s$};
\smallercross{(0,-.4)};
\node at (0,.41) {$\stringlabel{n+1}$};
\node at (0,-.81) {$\stringlabel{n-1}$};
\end{tikzpicture}=
(-1)^r
\begin{tikzpicture}[anchorbase,scale=1.2]
\draw[ultra thick] (0,.3) to (0,-.7);
\draw (.2,.3) to[out=-120,in=90,looseness=1] (-.4,-.2) to [out=-90,in=-120,looseness=1.5] (-.02,-.55);
\smallercross{(0,0)};
\closeddot{0,-.2};\node at (.42,-.2)  {$\scriptstyle \varpi_r+\rho$};
\closeddot{-.4,-.2};\node at (-.75,-.2) {$\scriptstyle n-s$};
\smallercross{(0,-.4)};
\node at (0,.41) {$\stringlabel{n}$};
\node at (0,-.81) {$\stringlabel{n-1}$};
\end{tikzpicture}.
\end{align*}
This is of degree $2(r-s)-n(n-1)$ while by \cref{basisthm} the lowest non-zero degree of the graded vector space $\Hom_{\cNB_t}(B^{\star(n-1)}, B^{\star(n+1)})$ is $-n(n-1)$, so it is automatically 0 if
$r < s$. Assume henceforth that $r \geq s$.
When $n=t=1$, so $r=s=1$, it is easy to see that we get 0
using \cref{curlstuff}, so assume also that $n \geq 2$.

In this paragraph, we treat the case that $r > s$.
We have that $\varpi_{r,n}+\rho_n = (n,n-1,\dots,n-s,\dots,n-r+1,n-r-1,\dots,1,0) \in \N^n$.
Let $\alpha := (n-s,n,n-1,\dots,\widehat{n-s},\dots,n-r+1,n-r-1,\dots,1,0)\in \N^n$, i.e., we have moved the entry $n-s$
to the beginning.
Let $\beta := (n-s,\alpha_1,\dots,\alpha_{n-1})$.
We have that 
\begin{align*}
\mathring{\u}_{r,n} \circ \y_{s,n} &=
(-1)^r
\begin{tikzpicture}[anchorbase,scale=1.2]
\draw[ultra thick] (0,.3) to (0,-.7);
\draw (.2,.3) to[out=-120,in=90,looseness=1] (-.4,-.2) to [out=-90,in=-120,looseness=1.5] (-.02,-.55);
\smallercross{(0,0)};
\closeddot{0,-.2};\node at (.42,-.2)  {$\scriptstyle \varpi_r+\rho$};
\closeddot{-.4,-.2};\node at (-.75,-.2) {$\scriptstyle n-s$};
\smallercross{(0,-.4)};
\node at (0,.41) {$\stringlabel{n}$};
\node at (0,-.81) {$\stringlabel{n-1}$};
\end{tikzpicture}\stackrel{\cref{samlabel}}{=}
(-1)^{r+s}
\begin{tikzpicture}[anchorbase,scale=1.2]
\draw[ultra thick] (0,.3) to (0,-.7);
\draw (.2,.3) to[out=-120,in=90,looseness=1] (-.4,-.2) to [out=-90,in=-120,looseness=1.5] (-.02,-.55);
\smallercross{(0,0)};
\closeddot{0,-.2};\node at (.2,-.2)  {$\scriptstyle\alpha$};
\closeddot{-.4,-.2};\node at (-.75,-.2) {$\scriptstyle n-s$};
\smallercross{(0,-.4)};
\node at (0,.41) {$\stringlabel{n}$};
\node at (0,-.81) {$\stringlabel{n-1}$};
\end{tikzpicture}=(-1)^{r+s}
\begin{tikzpicture}[anchorbase,scale=1.2]
\draw[ultra thick] (0,.3) to (0,-.7);
\draw (-.2,.3) to[out=-60,in=90,looseness=1] (.4,-.2) to [out=-90,in=-60,looseness=1.5] (.02,-.55);
\smallercross{(0,0)};
\closeddot{0,-.2};\node at (-.2,-.2)  {$\scriptstyle\beta$};
\closeddot{.4,-.2};\node at (.64,-.2) {$\scriptstyle \alpha_n$};
\smallercross{(0,-.4)};
\node at (0,.41) {$\stringlabel{n}$};
\node at (0,-.81) {$\stringlabel{n-1}$};
\end{tikzpicture}.
\end{align*}
In checking the second equality here, one also needs to 
observe that the 
term arising from the first term on the right hand side of \cref{samlabel} 
(which can definitely appear as $n-1\not\equiv t\pmod{2}$) 
is 0 due to the second relation from \cref{rels6}.
Now we have that $\beta_1=\beta_2=n-s$, so this is 0
by \cref{sammer}; again, when $s=1$, the term arising from the right hand side of \cref{sammerlabel} vanishes due to \cref{rels6}.

Finally, we need to treat the case that $r=s$ (and $n \geq 2$ still).
We let $\alpha := \varpi_{r,n}+\rho_n = (n,n-1,\dots,n-r+1,n-r-1,\dots,1,0) \in \N^n$, 
$\beta := (n-s,\alpha_1,\dots,\alpha_{n-1})$,
and $\gamma := (n-s,\alpha_2,\dots,\alpha_n)$.
As $r=s\geq 1$, 
the tuple $\gamma$ is a permutation of $\rho_n$,
and $\alpha_1 = n$.
Using \cref{sam} several more times like in the previous paragraph, we get that
\begin{align*}
\mathring{\u}_{r,n} \circ \y_{s,n} &=
(-1)^{r}
\begin{tikzpicture}[anchorbase,scale=1.2]
\draw[ultra thick] (0,.3) to (0,-.7);
\draw (.2,.3) to[out=-120,in=90,looseness=1] (-.4,-.2) to [out=-90,in=-120,looseness=1.5] (-.02,-.55);
\smallercross{(0,0)};
\closeddot{0,-.2};\node at (.2,-.2)  {$\scriptstyle\alpha$};
\closeddot{-.4,-.2};\node at (-.75,-.2) {$\scriptstyle n-s$};
\smallercross{(0,-.4)};
\node at (0,.41) {$\stringlabel{n}$};
\node at (0,-.81) {$\stringlabel{n-1}$};
\end{tikzpicture}=(-1)^{r}
\begin{tikzpicture}[anchorbase,scale=1.2]
\draw[ultra thick] (0,.3) to (0,-.7);
\draw (-.2,.3) to[out=-60,in=90,looseness=1] (.4,-.2) to [out=-90,in=-60,looseness=1.5] (.02,-.55);
\smallercross{(0,0)};
\closeddot{0,-.2};\node at (-.2,-.2)  {$\scriptstyle\beta$};
\closeddot{.4,-.2};\node at (.62,-.2) {$\scriptstyle \alpha_n$};
\smallercross{(0,-.4)};
\node at (0,.41) {$\stringlabel{n}$};
\node at (0,-.81) {$\stringlabel{n-1}$};
\end{tikzpicture}
=(-1)^{r+1}
\begin{tikzpicture}[anchorbase,scale=1.2]
\draw[ultra thick] (0,.3) to (0,-.7);
\draw (-.2,.3) to[out=-60,in=90,looseness=1] (.4,-.2) to [out=-90,in=-60,looseness=1.5] (.02,-.55);
\smallercross{(0,0)};
\closeddot{0,-.2};\node at (-.29,-.2)  {$\scriptstyle s_1 \beta$};
\closeddot{.4,-.2};\node at (.62,-.2) {$\scriptstyle \alpha_n$};
\smallercross{(0,-.4)};
\node at (0,.41) {$\stringlabel{n}$};
\node at (0,-.81) {$\stringlabel{n-1}$};
\end{tikzpicture}\\
&=(-1)^{r+1}
\begin{tikzpicture}[anchorbase,scale=1.2]
\draw[ultra thick] (0,.3) to (0,-.7);
\draw (.2,.3) to[out=-120,in=90,looseness=1] (-.4,-.2) to [out=-90,in=-120,looseness=1.5] (-.02,-.55);
\smallercross{(0,0)};
\closeddot{0,-.2};\node at (.2,-.2)  {$\scriptstyle\gamma$};
\closeddot{-.4,-.2};\node at (-.57,-.2) {$\scriptstyle n$};
\smallercross{(0,-.4)};
\node at (0,.41) {$\stringlabel{n}$};
\node at (0,-.81) {$\stringlabel{n-1}$};
\end{tikzpicture}
=
\begin{tikzpicture}[anchorbase,scale=1.2]
\draw[ultra thick] (0,.3) to (0,-.7);
\draw (.2,.3) to[out=-120,in=90,looseness=1] (-.4,-.2) to [out=-90,in=-120,looseness=1.5] (-.02,-.55);
\smallercross{(0,0)};
\closeddot{0,-.2};\node at (.2,-.2)  {$\scriptstyle\rho$};
\closeddot{-.4,-.2};\node at (-.57,-.2) {$\scriptstyle n$};
\smallercross{(0,-.4)};
\node at (0,.41) {$\stringlabel{n}$};
\node at (0,-.81) {$\stringlabel{n-1}$};
\end{tikzpicture}=
\begin{tikzpicture}[anchorbase,scale=1.2]
\draw[ultra thick] (0,.3) to (0,-.7);
\draw (.2,.3) to[out=-120,in=90,looseness=1] (-.4,-.3) to [out=-90,in=-120,looseness=1.5] (-.02,-.5);
\smallercross{(0,-.25)};
\closeddot{-.4,-.3};\node at (-.57,-.3) {$\scriptstyle n$};
\node at (0,.41) {$\stringlabel{n}$};
\node at (0,-.81) {$\stringlabel{n-1}$};
\end{tikzpicture}=
\begin{tikzpicture}[anchorbase,scale=1.2]
\draw[ultra thick] (0,.3) to (0,-.7);
\draw (.5,.3) to[out=-90,in=90,looseness=1] (-.4,-.2) to [out=-90,in=180] (-.1,-.5) to [out=0,in=-90,looseness=1.5] (.25,.3);
\smallercross{(0,-.2)};
\closeddot{-.4,-.2};\node at (-.57,-.2) {$\scriptstyle n$};
\node at (-0.05,.41) {$\stringlabel{n-1}$};
\node at (0,-.81) {$\stringlabel{n-1}$};
\end{tikzpicture}.
\end{align*}
This is 0 by \cref{curly}, using that $n-1\not\equiv t\pmod{2}$.
\end{proof}

\begin{corollary}\label{id2ccor}
Assume that $n\equiv t\pmod{2}$. For $0 \leq r\leq n$ and $1\leq s \leq n$, we have that
$\e_{r,n} \circ \f_{s,n} = \f_{s,n} \circ \e_{r,n} = 0$.
\end{corollary}

\begin{proof}
This is clear from \cref{id0,id2c}.
\end{proof}

\begin{theorem}\label{id3}
The following hold for $n \geq 0$:
\begin{enumerate}
\item
If $n\equiv t\pmod{2}$ then
$\{\e_{r,n}, \f_{s,n}\:|\:0 \leq r\leq n, 1\leq s \leq n\}$
is a set of mutually orthogonal homogeneous idempotents
whose sum is $B \star \e_n$.
Each of the idempotents $\e_{r,n}\:(0 \leq r \leq n)$ is conjugate
to $\e_{n+1} = \e_{0,n}$
since 
$\e_{n+1} = \u_{r,n} \circ \vv_{r,n}$
and $\e_{r,n} = \vv_{r,n} \circ \u_{r,n}$
for $r=1,\dots,n$.
Each of the idempotents
$\f_{s,n}\:(1 \leq s \leq n)$ is conjugate to
$\e_{n-1}$ since $\e_{n-1} = \x_{s,n} \circ \y_{s,n}$
and $\f_{s,n} = \y_{s,n} \circ \x_{s,n}$
for $s=1,\dots,n$.
\item
If $n\not\equiv t\pmod{2}$ 
then $\{\e_{r,n}+\f_{r,n}\:|\:0 \leq r \leq n\}$
is a set of mutually orthogonal homogeneous idempotents
whose sum is $B \star \e_n$.
Each
of these idempotents is conjugate to $\e_{n+1}=\e_{0,n}$ since, recalling that $\w_{r,n} = \u_{r,n}-\u_{r,n} \circ \vv_{0,n}$,
we have that
$\e_{n+1} = \w_{r,n} \circ \vv_{r,n}$
and $\e_{r,n}+\f_{r,n} = \vv_{r,n} \circ \w_{r,n}$ for $r=1,\dots,n$.
\end{enumerate}
\end{theorem}

\begin{proof}
(1)
The fact that $\e_{r,n}\:(0 \leq r \leq n)$ are mutually orthogonal idempotents
follows from \cref{id2}. 
The fact that $\f_{s,n}\:(1 \leq s \leq n)$ are mutually orthogonal idempotents
follows from \cref{id2bcor}.
The orthogonality of each $\e_{r,n}\:(0 \leq r \leq n)$
with each $\f_{s,n}\:(1 \leq s \leq n)$ follows from
\cref{id2ccor}.
These idempotents sum to $B \star \e_n$ by \cref{hug}.
Also
$\u_{r,n} \circ \vv_{r,n} = \e_{n+1}$ by \cref{id1},
and $\vv_{r,n} \circ \u_{r,n} = \e_{r,n}$ by \cref{id0}.
Finally, $\x_{s,n} \circ \y_{s,n} = \e_{n-1}$
by \cref{id2b}, and
$\y_{s,n} \circ \x_{s,n} = \f_{s,n}$ by \cref{id0}.

\vspace{2mm}
\noindent
(2)
We first show that 
$\e_{r,n}+\f_{r,n}\:(0 \leq r \leq n)$
are mutually orthogonal idempotents by checking that
$$
(\e_{r,n}+\f_{r,n}) \circ
(\e_{s,n}+\f_{s,n})
 = \delta_{r,s} (\e_{r,n}+\f_{r,n})
 $$
 for $0 \leq r,s \leq n$.
If $r=0$ this follows because
$\f_{0,n}=0$, $\e_{0,n} \circ \e_{s,n}
= \delta_{0,s} \e_{0,n}$
and, assuming $s > 0$, we have that
$\e_{0,n} \circ \f_{s,n} 
= - \e_{0,n} \circ \e_{s,n} \circ \e_{0,n}
= 0$, all by \cref{id2}.
If $r > 0$ and $s=0$ it follows because
$\e_{r,n} \circ \e_{0,n} = -\f_{r,n}$
and $\f_{r,n} \circ \e_{0,n} = - \e_{r,n} \circ \e_{0,n} \circ \e_{0,n} = - \e_{r,n} \circ \e_{0,n} = \f_{r,n}$ by \cref{id2}.
Finally suppose that $1 \leq r,s \leq n$. Then
by \cref{id2} we have that
\begin{align*}
(\e_{r,n}+\f_{r,n})\circ (\e_{s,n}+\f_{s,n})&=
\e_{r,n} \circ \e_{s,n}
+ \e_{r,n} \circ \f_{s,n}
+ \f_{r,n} \circ \e_{s,n}
+ \f_{r,n} \circ \f_{s,n}\\
&=
\e_{r,n} \circ \e_{s,n}
- \e_{r,n} \circ \e_{s,n} \circ \e_{0,n}
- \e_{r,n} \circ \e_{0,n} \circ\e_{s,n}
+ \e_{r,n} \circ \e_{0,n} \circ \e_{s,n} 
\circ \e_{0,n}\\ &= \delta_{r,s}
\e_{r,n} - \delta_{r,s} \e_{r,n} \circ \e_{0,n}
= \delta_{r,s}(\e_{r,n} + \f_{r,n}).
\end{align*}
We have that $\sum_{r=0}^n (\e_{r,n}+\f_{r,n}) = B \star \e_n$ by \cref{hug}.
Finally, using \cref{id0,id1}, \cref{id2}
and $\u_{0,n} = \vv_{0,n} = \e_{0,n}$, 
we have that
\begin{align*}
\w_{r,n} \circ \vv_{r,n}
&=
\u_{r,n} \circ \vv_{r,n} - \u_{r,n} \circ \u_{0,n}
\circ \vv_{r,n} = \e_{n+1},\\
\vv_{r,n} \circ \w_{r,n}
&=
\vv_{r,n} \circ \u_{r,n} - \vv_{r,n} \circ \u_{r,n} \circ \e_{0,n}
=\e_{r,n} - \e_{r,n} \circ \e_{0,n}
= \e_{r,n}+\f_{r,n}
\end{align*}
for $1 \leq r \leq n$.
\end{proof}

\subsection{Locally unital graded algebras and modules}\label{blah}
Before explaining the full significance of 
\cref{id3}, we need to review some basic terminology.
Suppose that $\cA$ is any small 
graded category and let $\I$ be its object set.
The {\em path algebra} of $\cA$ is the
graded algebra
$$
A = \bigoplus_{i,j \in \I} 1_i\, A 1_j
\qquad\text{where}\qquad 1_i\, A 1_j := \Hom_\cA(j,i),
$$
with 
multiplication induced by composition in $\cA$. 
In general, this is {\em locally unital} rather than unital, equipped with the distinguished system $1_i\:(i \in \I)$ of mutually orthogonal idempotents 
arising from the identity endomorphisms of the objects of $\cA$. 

The {\em center} $Z(A)$ is the commutative 
subalgebra of the unital graded algebra
$\prod_{i \in \I} 1_i A 1_i$
consisting of tuples
$(z_i)_{i \in \I}$ such that $\theta\, z_j = z_i\, \theta$
for all $i,j \in \I$ and $\theta \in 1_i A 1_j$. 
This is a direct translation of the definition of the center of the category $\cA$.
Given a (unital) commutative graded algebra $R$, we say that $A$ is a {\em graded $R$-algebra}
if we are given a unital graded algebra 
homomorphism
 $\eta:R \rightarrow Z(A)$. Then each
 subspace
 $1_i A 1_j$ is naturally a graded $R$-module.

By a {\em graded left $A$-module}, we 
mean a module $V$ as usual which is itself locally unital in the sense that
$V = \bigoplus_{i \in \I} 1_i V$.
We sometimes refer to $1_i V$ as the
{\em $i$-weight space} of $V$.
There are also the obvious notions of 
graded right $A$-modules
and, given another locally unital graded algebra $B$, graded $(A,B)$-bimodules. 
For graded left $A$-modules $V$ and $W$
and $d \in \Z$,
we write $\Hom_{A}(V,W)_d$
for the vector space of all ordinary $A$-module
homomorphisms $f:V \rightarrow W$
such that $f(V_n) \subseteq W_{n+d}$ for each $n \in \Z$.
Then the graded vector space
$$
\Hom_{A}(V,W) := \bigoplus_{d \in \Z} \Hom_{A}(V,W)_d
$$
is a morphism space in 
the 
graded category $A\gMOD$ of graded left $A$-modules. 
We denote the underlying category 
consisting of the same objects but just the
degree-preserving morphisms by $A\gMod$.
This is the usual Abelian category of graded left $A$-modules. It is equipped with the upward grading shift functor $q$ defined as in the {\em General conventions}, and we have that
\begin{equation}\label{grconv}
\Hom_A(V,W)_d = 
\left(q^{-d} \Hom_A(V,W)\right)_0=
\Hom_A(V, q^{-d} W)_0 = \Hom_A(q^{d} V, W)_0.
\end{equation}
We use the symbol $\cong$ to denote (degree-preserving) isomorphism in $A \gMod$.
We write $\End_A(V)$ for $\Hom_A(V,V)$, which is a graded algebra.

Let $A\pgMod$ be the full subcategory
of $A\gMod$
consisting of finitely generated
projective graded modules.
Also let $K_0(A)$ denote the split Grothendieck group of the additive category $A\pgMod$.
This is a $\Z[q,q^{-1}]$-module with the action of $q$ induced by the grading shift functor.
One could also define $K_0(A)$ equivalently
as the split Grothendieck group of the graded Karoubi envelope of $\cA$, since the latter category is 
contravariantly equivalent to $A\pgMod$
by Yoneda's Lemma. We will not take this point of view here, but note that some care is needed in the identification since contravariant equivalences 
interchange $q$ with $q^{-1}$.

Assume in this paragraph that $A$ is {\em locally finite dimensional and bounded below}, meaning that for every $i,j \in \I$,
the graded vector space $1_i\, A 1_j$ 
is locally finite dimensional, i.e., each of its graded pieces $1_i\, A_d 1_j$ are finite dimensional, and $1_i\, A_d 1_j = 0$
for $d \ll 0$.
Then $K_0(A)$ can be understood in purely combinatorial terms. To explain what we mean, referring to \cite[Sec.~2]{GTB} for more details, we note to start with that the weight spaces of any irreducible graded left $A$-module $L$ are finite dimensional,
and Schur's Lemma holds:
\begin{equation}\label{schurslemma}
\End_A(L) = \kk.
\end{equation}
We say that a graded left $A$-module $V$
is {\em locally finite dimensional} if
$1_i V_d$ is finite dimensional for each $i \in \I$ and $d \in \Z$, and 
{\em bounded below} if 
for each $i \in \I$ we have that 
$1_i V_d = 0$ for $d \ll 0$.
Since the distinguished projective
modules $A 1_i\:(i \in \I)$ are locally finite dimensional and bounded below, it follows that any finitely generated graded left $A$-module also has these properties.
Any graded left $A$-module has an injective hull in $A\gMod$, and any finitely generated graded left $A$-module has a projective cover in $A\gMod$, the latter being a summand of a finite direct sum of degree-shifted copies of the distinguished projective modules $A 1_i\:(i \in \I)$.
Let $L(b)\:(b \in \B)$
be a full set of representatives for the irreducible graded left $A$-modules (up to isomorphism and grading shift), 
and define $P(b)$ to be a projective cover of $L(b)$.
The {\em graded multiplicity}
of $L(b)$ in 
a locally finite-dimensional graded module $V$
is the formal series 
\begin{align*}
[V:L(b)]_q &:=\sum_{d \in \Z} \max 
\bigg(
\left|\{r=1,\dots,n\:|\:V_r / V_{r-1} \cong q^d L(b) \}\right|
\:\bigg|\:
\begin{array}{ll}
\text{for all finite graded filtrations}\\
0 = V_0 \subseteq \cdots \subseteq V_n = V
\end{array}
\bigg)
 q^d.
\end{align*}
Schur's Lemma implies that
\begin{equation}
[V:L(b)]_q = \dim_q \Hom_A(P(b),V).
\end{equation}
Note also that this belongs to $\N\lround q\rround$ when $V$ is finitely generated.
Finally, any finitely generated projective graded left $A$-module $P$
satisfies
\begin{equation}\label{specase}
P \cong \bigoplus_{b \in \B} P(b)^{\oplus 
\overline{\dim_{q} \Hom_A(P, L(b))}}.
\end{equation}
Now it follows that
that $K_0(A)$ is a free $\Z[q,q^{-1}]$-module
with basis $[P(b)]\:(b \in \B)$.

Another basic notion involves induction and restriction. 
For this, we 
start with a pair of small
graded categories, $\cA$ and $\cB$, 
with object sets denoted $\I$ and $\J$, respectively.
Let $A$ and $B$ be their path algebras.
Given a graded functor
$F: \cA \rightarrow \cB$,
there is a graded functor
\begin{equation}
\Res_F:B\gMOD
\rightarrow A\gMOD
\end{equation}
called {\em restriction along $F$}.
This takes a graded left $B$-module $V$
to the graded vector space
$$
1_F V := \bigoplus_{i \in \I} 1_{F i} V
$$
with $\theta \in 1_i A 1_j=\Hom_{\cA}(j,i)$ 
acting as the linear map 
$F\theta:1_{Fj} V \rightarrow 1_{Fi} V$
between the summands indexed by $j$ and $i$,
and as 0 on all other summands.
This notation is for graded left $B$-modules,
but it is readily adapted to a graded right $B$-module $V$, letting
$$
V 1_F := \bigoplus_{i \in \cA} V 1_{Fi}
$$
which is a graded right $A$-module.
The functor $\Res_F$ is isomorphic 
to 
$\bigoplus_{i \in \I} \Hom_B(B 1_{Fi},-)$. Hence, by adjointness of tensor and hom
for locally unital algebras (e.g., see \cite[Lem.~2.7]{BS}), it has a left adjoint \begin{align}\label{Ind}
\Ind_F := B 1_F \otimes_A -&:A\gMOD \rightarrow B\gMOD,
\end{align}
where $B 1_F$ is the graded $(B,A)$-bimodule obtained by restricting the regular $(B,B)$-bimodule
$B$ on the right.
We refer to $\Ind_F$ as {\em induction along $F$}.
If $\alpha:F \Rightarrow G$
is a graded natural transformation between 
graded functors $F,G:\cA \rightarrow \cB$,
we obtain graded bimodule
homomorphisms 
$B 1_G \rightarrow B 1_F$
and
$1_F B \rightarrow 1_G B$
defined
by the linear maps
$1_j B 1_{Gi} \rightarrow 1_{j} B 1_{Fi},
\theta \mapsto \theta \circ \alpha_i$
and 
$1_{Fi} B 1_j \rightarrow 1_{Gi} B 1_j,
\theta \mapsto \alpha_i \circ \theta$,
respectively,
for $i \in \I, j \in \J$.
These bimodule homomorphisms define
graded natural transformations
$\Ind_\alpha:\Ind_G \Rightarrow \Ind_F$
and $\Res_\alpha:\Res_F \Rightarrow \Res_G$.

Suppose finally that the small graded category $\cA$ is
monoidal, with tensor product bifunctor
\begin{equation}\label{starf}
-\star-:\cA \boxtimes \cA \rightarrow \cA,
\end{equation}
where we are using $\boxtimes$ to denote linearized Cartesian product.
Then there is an induced tensor product
bifunctor making $A\gMOD$ into a graded monoidal category in
its own right.
We call this the {\em induction product};
it is also known as {\em Day convolution}.
To define it, observe that
the graded algebra $A \otimes A$ is the path algebra of the graded category
$\cA \boxtimes \cA$. 
The induction product is the graded bifunctor 
\begin{equation}\label{indprod}
-\,\ostar\,-:
A\gMOD \boxtimes A\gMOD \rightarrow A \gMOD
\end{equation}
that is the composition of the usual 
tensor product
$-\otimes -:A\gMOD \boxtimes A\gMOD \rightarrow A\otimes A \gMOD$
followed by the
functor $\Ind_{-\star-}:
\NB\otimes A \gMOD
\rightarrow A\gMOD$
defined by induction
along \cref{starf}.
Note that $-\,\ostar\,-$ is right exact in each argument but it is not necessarily exact.
It is clear from the definition
that
\begin{equation}\label{fitting}
A 1_i \, \ostar \, A 1_j \cong A 1_{i \star j}
\end{equation}
for $i,j \in \I$.
From this, one deduces that 
the restriction of $-\,\ostar\,-$
makes $A\pgMod$ into a monoidal category.
Consequently, $K_0(A)$ is actually a $\Z[q,q^{-1}]$-algebra with
multiplication satisfying
\begin{equation}\label{fitting2}
[A 1_i] [A 1_j] = [A 1_i \,\ostar\, A 1_j] = [A 1_{i\star j}].
\end{equation}

\subsection{Identification of the Grothendieck ring}
Now we apply the general setup just explained to the nil-Brauer category. 
We denote the path algebra of $\cNB_t$
for the fixed value of $t$
simply by $\NB$. Its distinguished idempotents
arising from the identity endomorphisms of
$B^{\star n}\:(n \in \N)$ will be denoted 
by $1_n\:(n \in \N)$. So we have that
$$
\NB = \bigoplus_{m,n \in \N} 1_m \NB 1_n
\qquad\text{where}\qquad
1_m \NB 1_n = \Hom_{\cNB_t}(B^{\star n}, B^{\star m}).
$$
\cref{basisthm} implies that $\NB$ is locally finite dimensional and bounded below,
so that we are in the situation 
discussed in the fifth paragraph of \cref{blah}. 
Note also that $\NB$ is a graded 
$\Gamma$-algebra, with $\beta \in \Gamma$ acting on a morphism by horizontal composition on the right with $\gamma_t(\beta)$ (recall \cref{gammacor}).
Since $\cNB_t$ is monoidal,
we have the induction product
$-\,\ostar\,-:\NB\gMOD \boxtimes \NB\gMOD
\rightarrow \NB\gMOD$ defined as in \cref{indprod}.
It makes $K_0(\NB)$
into a $\Z[q,q^{-1}]$-algebra.
Our goal is to identify this with
the integral form $\UAi$ of the iquantum group.

Recalling the idempotent
$\e_n \in 1_n \NB 1_n$ from \cref{endef},
we define
\begin{equation}\label{pndef}
P(n) := q^{\frac{1}{2}n(n-1)} \NB\, \e_n.
\end{equation}
This is a finitely generated projective graded
left $\NB$-module.
In particular, we have that $P(0) = \NB 1_0$
and $P(1) = \NB 1_1$.
Also let
\begin{equation}\label{pfun}
B:= P(1) \,\ostar\, -:\NB\gMOD \rightarrow \NB\gMOD
\end{equation}
be the endofunctor defined by taking the induction product with the projective module $P(1)$
associated to the generating object $B$ of $\cNB_t$.
From \cref{fitting}, we have that
\begin{equation}\label{siem}
B (\NB 1_n) \cong \NB 1_{n+1}.
\end{equation}
Since it is clearly additive, it follows
that $B$ takes finitely generated projectives
to finitely generated projectives, i.e.,
it restricts to an endofunctor of $\NB \pgMod$.
This is all that we need for now, but we will say more about $B$ viewed 
as an endofunctor of the Abelian category $\NB\gMod$ in \cref{pfunctor} below.

\begin{lemma}\label{missed}
For $n \in \N$, we have that
$$
B P(n) \cong 
\begin{dcases}
P(n+1)^{\oplus [n+1]} \oplus P(n-1)^{\oplus [n]}
&\text{if $n\equiv t\pmod{2}$}\\
P(n+1)^{\oplus [n+1]}
&\text{if $n\not\equiv t\pmod{2}$.}
\end{dcases}
$$
\end{lemma}

\begin{proof}
First consider the case that $n\not\equiv t\pmod{2}$.
By the first part of \cref{id3}(2), we have that $B \star \e_n = \sum_{r=0}^n (\e_{r,n} + \f_{r,n})$ as a sum of mutually orthogonal idempotents. As in 
\cref{fitting}, we deduce that
$$
B P(n) = q^{\frac{1}{2}n(n-1)}
\NB 1_1 \,\ostar\, \NB\, \e_n \cong
\bigoplus_{r=0}^n 
q^{\frac{1}{2}n(n-1)}
\NB (\e_{r,n}+\f_{r,n}).
$$
To complete the proof in this case, we claim that
$q^{\frac{1}{2}n(n-1)}
\NB (\e_{r,n}+\f_{r,n}) \cong q^{2r-n} P(n+1)$
for any $0 \leq r \leq n$.
The second part of
\cref{id3}(2) shows that
right multiplication by $\vv_{r,n}$
defines an invertible $\NB$-module homomorphism
$\NB (\e_{r,n} + \f_{r,n})
\stackrel{\sim}{\rightarrow}
\NB\, \e_{n+1}$
with inverse given by right multiplication by
$\w_{r,n}$. By its definition \cref{uv}, $\vv_{r,n}$ is of degree $-2r$.
Recalling \cref{grconv}, this shows that
$$
q^{\frac{1}{2}n(n-1)}
\NB (\e_{r,n}+\f_{r,n})
\cong q^{\frac{1}{2}n(n-1)+2r} \NB\, \e_{n+1}
\cong q^{-\frac{1}{2}(n+1)n+\frac{1}{2}n(n-1)+2r} P(n+1) = q^{2r-n} P(n+1),
$$
as claimed.

Instead, suppose that $n\equiv t\pmod{2}$.
Then the first part of \cref{id3}(1) gives that
$$
B P(n) = q^{\frac{1}{2}n(n-1)} 
\NB 1_1 \,\ostar\, \NB\, \e_n 
\cong \bigoplus_{r=0}^n 
q^{\frac{1}{2}n(n-1)}
\NB\, \e_{r,n} 
\oplus \bigoplus_{s=1}^n
q^{\frac{1}{2}n(n-1)}
\NB\,\f_{s,n}.
$$
Now
it suffices to show that
$q^{\frac{1}{2}n(n-1)}
\NB\, \e_{r,n} \cong q^{2r-n} P(n+1)$
for $0 \leq r \leq n$
and that
$q^{\frac{1}{2} n(n-1)}
\NB\,\f_{s,n} \cong q^{2s-n-1} P(n-1)$
for $1 \leq s \leq n$.
The first assertion here follows from the second part of \cref{id3}(1)
just like in the previous paragraph (replacing
$\w_{r,n}$ with $\u_{r,n}$).
To prove the second assertion, 
right multiplication by $\y_{s,n}$
defines an invertible $\NB$-module homomorphism
$\NB\,\f_{s,n}
\stackrel{\sim}{\rightarrow}
\NB\, \e_{n-1}$
with inverse given by right multiplication by
$\x_{s,n}$. By its definition \cref{xy}, $\y_{s,n}$ is of degree $2n-2s$, so this shows that
$$
q^{\frac{1}{2}n(n-1)}
\NB\,\f_{s,n}
\cong q^{\frac{1}{2}n(n-1)+2s-2n} \NB\, \e_{n-1}
\cong q^{-\frac{1}{2}(n-1)(n-2)+\frac{1}{2}n(n-1)+2s-2n} P(n-1) = q^{2s-n-1} P(n-1).
$$
\end{proof}

\begin{theorem}\label{therealthing}
The modules
$P(n)\:(n \geq 0)$ 
give a complete set of indecomposable
projective graded left $\NB$-modules
(up to isomorphism and grading shift).
Moreover, there is a unique 
$\Z[q,q^{-1}]$-algebra isomorphism
$$
\kappa_t:K_0(\NB) \stackrel{\sim}{\rightarrow} \UAi
$$
such that
\begin{enumerate}
\item
$\kappa_t([B P]) = b \kappa_t([P])$
for any finitely generated projective graded module $P$.
\end{enumerate}
The following properties also hold
for finitely generated projective graded modules $P, Q$ and $n \geq 0$:
\begin{enumerate}
\item[(2)]
$\kappa_t([\NB 1_n]) = b^n$.
\item[(3)] $\kappa_t([P(n)]) = b^{(n)}$.
\item[(4)]
$\Hom_{\NB} (Q, P) 
\cong \Gamma^{\oplus \overline{(\psi^\imath(\kappa_t([P])),\; \kappa_t([Q]))^\imath}}$ as a graded $\Gamma$-module.
\end{enumerate}
\end{theorem}

\begin{proof}
Let $\lambda_t:\UAi \rightarrow K_0(\NB)$
be the $\Z[q,q^{-1}]$-module homomorphism
taking $b^{(n)}$ to $[P(n)]$ for each $n \geq 0$.
By \cref{idprecurrence} and \cref{missed},
it follows that
$\lambda_t$ intertwines the endomorphism of
$\UAi$ defined by left 
multiplication by $b$
with the endomorphism of $K_0(\NB)$ induced
by the 
functor $B:\NB\pgMod \rightarrow \NB \pgMod$.
Hence, also using \cref{siem},
we have that 
\begin{equation}\label{work1}
\lambda_t(b^n)
=\lambda_t(b^n b^{(0)})
= [B^n P(0)]
= [B^n \NB 1_0]
= [\NB 1_n].
\end{equation}
We also have that
\begin{equation}\label{work2}
\Hom_{\NB}(P(n), P(m))
\cong \Gamma^{\oplus\overline{(b^{(m)},\;b^{(n)})^\imath}}
\end{equation}
for any $m,n \geq 0$. To see this, since
both $\UAi$ and $K_0(\NB)$ are free $\Z[q,q^{-1}]$-modules, it is harmless to extend scalars
from $\Z[q,q^{-1}]$ to $\Q(q)$.
Then $b^{(m)}$ and $b^{(n)}$ are bar-invariant 
$\Q(q)$-linear combinations of 
the elements $b^k\:(k \geq 0)$ (see \cref{idp} for the explicit formula which is not needed here).
Applying $\lambda_t$
gives that $[P(m)]$ and $[P(n)]$ are corresponding linear combinations of $[\NB 1_k]\:(k \geq 0)$.
In this way, the proof of \cref{work2}
is reduced to checking that
\begin{equation}\label{work3}
\Hom_{\NB}(\NB 1_n, \NB 1_m)
\cong \Gamma^{\oplus
\overline{(b^m,\;b^n)^\imath}}
\end{equation}
for all $m,n \geq 0$, which follows from \cref{toinfinityandbeyond}.

By \cref{dontreallyneedthedetailedformula},
we have that $(b^{(n)},\;b^{(n)})^\imath 
\in 1 + q^{-1} \Z\llbracket q^{-1} \rrbracket$.
Hence, by \cref{work2},
we have that $\End_\NB(P(n))_0 \cong \kk$.
This implies the indecomposability of $P(n)$.
Moreover, the isomorphism classes
$[P(n)]\:(n \geq 0)$ are linearly independent
over $\Z[q,q^{-1}]$. This follows because
the matrix $$
\big(\dim_q \Hom_{\NB}(P(n), P(m)\big)_{m,n \geq 0}
$$ 
is invertible by \cref{work2} and
\cref{dontreallyneedthedetailedformula} (or the non-degeneracy of the form
$(\cdot,\cdot)^\imath$).
Hence, for $m \neq n$ the module $P(n)$ is not isomorphic to any grading shift of $P(m)$.
Finally, we observe that
any indecomposable projective graded left $\NB$-module is isomorphic
to $q^d P(n)$ for  unique $d \in \Z, n \in \N$. This is true because
each left ideal $\NB 1_n$ is isomorphic to a
direct sum of grading shifts of the modules
$P(m)$ for $m \geq n$, as follows by induction on $n$ using \cref{siem} and \cref{missed}.

We have now proved the first sentence in the statement of the theorem.
It follows that the isomorphism classes
$[P(n)]\:(n \geq 0)$ give a basis for $K_0(\NB)$ as a free $\Z[q,q^{-1}]$-module.
We deduce immediately that $\lambda_t$ 
is an isomorphism of free $\Z[q,q^{-1}]$-modules. 
Let $\kappa_t := \lambda_t^{-1}$.
This satisfies the property (1). Moreover,
$$
\kappa_t(b^m\cdot b^n) = \kappa_t(b^{m+n})
= [\NB 1_{m+n}] = [\NB 1_m
\,\ostar\,\NB 1_n]
= [\NB 1_m] [\NB 1_n].
$$
It follows that the $\Q(q)$-module isomorphism
$\Q(q)\otimes_{\Z[q,q^{-1}]} \UAi \stackrel{\sim}{\rightarrow} \Q(q) \otimes_{\Z[q,q^{-1}]} K_0(\NB)$
induced by $\kappa_t$ is actually a $\Q(q)$-algebra isomorphism. Hence, $\kappa_t$ itself is a $\Z[q,q^{-1}]$-algebra isomorphism. The uniqueness of an algebra isomorphism $\kappa_t$ satisfying the property (1) is clear.
We also get (2) and (3) since $\lambda_t$
satisfies the appropriate inverse properties by the definition of $\lambda_t$ and \cref{work1}.
Finally, (4) follows from \cref{work2} and sesquilinearity
of the Cartan form.
\end{proof}

\begin{corollary}\label{idempotenttheorem}
The idempotents
$\e_{n}\:(n \geq 0)$ from \cref{endef}
give a complete set of primitive homogeneous idempotents in the nil-Brauer category (up to conjugacy).
\end{corollary}

\begin{proof}
We need to establish the following two assertions:
\begin{itemize}
\item Each $\e_n$
is a primitive homogeneous idempotent in the path algebra $\NB$.
\item Given a primitive homogeneous idempotent $\e \in 1_m \NB 1_m$, 
there is a unique $n \geq 0$ and elements 
$x \in 1_m \NB 1_n, y \in 1_n \NB 1_m$
such that $\e = xy$ and $\e_n = yx$.
\end{itemize}
The first of these is equivalent to
the indecomposability of the projective graded module
$\NB\, \e_n$ established in \cref{therealthing}.
To prove the second assertion,
$\NB\, \e$ is an indecomposable projective graded module, hence, it is isomorphic to $q^d \NB\, \e_n$
for unique $d \in \Z, n \in \N$ by
the definition of $P(n)$ and \cref{therealthing} again.
Let $\theta:\NB\, \e \stackrel{\sim}{\rightarrow} q^d \NB\, \e_n$ be an isomorphism.
Since $\Hom_\NB (\NB\, \e, q^d \NB\, \e_n)_0 =
\Hom_{\NB}(\NB\, \e, \NB\, \e_n)_{-d} \cong \e \NB_{-d} \e_n$, there is a unique $x \in \e \NB_{-d} \e_n$ such that $\theta$ is right multiplication by $x$.
Similarly, there is a unique $y \in \e_n \NB_{d} \e$ such that $\theta^{-1}$
is right multiplication by $y$.
We then have that $xy = \e$ and $yx = \e_n$ as required.
\end{proof}

\begin{corollary}\label{lastgasp}
For $n \geq 0$, we have that
$$
\NB 1_n \cong \bigoplus_{m=0}^{\lfloor \frac{n}{2}\rfloor}
P(n-2m)^{\oplus
\left([n-2m]!\sum_{\alpha \in \Par_t(m\times (n-2m))}  [\alpha_1+1]^2 \cdots [\alpha_m+1]^2\right)}.
$$
\end{corollary}

\begin{proof}
This follows from the theorem together with \cref{texas}.
\end{proof}

Since $K_0(\cNB_t)$ in the introduction is naturally identified with $K_0(\NB)$ here,
\cref{introthm:PIM,introthm:iso} as formulated in the introduction follow from \cref{missed,therealthing}.

%% file: s5-representations.tex
\setcounter{section}{4}

\section{Representation theory}\label{sect5}

In this section, we introduce an explicit graded triangular basis for the path algebra $\NB$ of the nil-Brauer category $\cNB_t$, which fits well with the general machinery developed in \cite{GTB} (extending ideas from \cite{GRS} and \cite[Sec.~5.4]{BS}). This allows us to define standard and proper standard modules, and to classify irreducible graded $\NB$-modules by their lowest weights.
Then, in \cref{Thm:SES}, 
we establish the existence of a certain short exact sequence of functors which can be viewed as a categorification of part of \cref{jy2}.
We use this to describe the effect of the functor $B$ on standard and proper standard modules, thereby proving \cref{introthm:ses} from the introduction.
Finally, we prove character formulae for 
proper standard and irreducible modules,
thereby proving \cref{introthm:decomposition,introthm:char},
and derive further branching rules.

\subsection{Triangular basis}
Recall
that $\RSD(m\times n)$ is a set of representatives for the $\sim$-equivalence classes of 
reduced $m \times n$ string diagrams,
two such diagrams being equivalent if they define the same matchings on their boundaries.
\cref{basisthm} shows moreover 
that $\NB$, the path algebra of $\cNB_t$,
is free as a $\Gamma$-algebra with basis
$\bigcup_{m,n \geq 0} \RSD(m\times n)$.
We now distinguish three special types of 
reduced string diagrams:
\begin{itemize}
\item[(X)]
Reduced string diagrams which only involve cups and non-crossing
propagating strings.
\item[(H)]
Reduced string diagrams with no cups or caps, just propagating strings (which are allowed to cross).
\item[(Y)]
Reduced string diagrams which only involve caps and non-crossing propagating strings.
\end{itemize}
From now on, we actually only need
representatives for the
$\sim$-equivalence classes of 
undotted reduced string diagrams 
of these three types. For types X or Y, we also choose a distinguished point on each
cup or cap.
For type H, we choose a distinguished point on each propagating string.
Then let $\X(a,n) \subset 1_a \NB 1_n$, $\HH(n)
\subset 1_n \NB 1_n$ and $\Y(n,b) \subset 1_n \NB 1_b$ be the sets obtained from the chosen 
$\sim$-equivalence class representatives
of $a \times n$ string diagrams of type X,
of $n \times n$ string diagrams of type H,
and of $n \times b$ string diagrams of type Y, respectively,
obtained by adding closed
dots labeled by non-negative 
multiplicities at each of the distinguished points.
Clearly, $\X(a,n) = \Y(n,b) = \varnothing$
unless $a \geq n \leq b$,
and $\X(n,n) = \{1_n\}=\Y(n,n)$.
Shorthand:
\begin{align*}
\X(n) &:= \bigcup_{a \geq n} \X(a,n),&
\Y(n) &:= \bigcup_{b \geq n} \Y(n,b).
\end{align*}
Also let $\H(n)$ be the set of morphisms
obtained from the ones in $\HH(n)$ by
placing ordered monomials
$\O_1^{m_1} \O_3^{m_3} \O_5^{m_5}\cdots$ in the odd $\O_r$ at the right hand boundary (recall \cref{Or}).
The latter are the images 
of a basis for $\Gamma$
under the isomorphism
$\gamma_t:\Gamma\stackrel{\sim}{\rightarrow}
\End_{\cNB_t}(\one)$ from \cref{gammacor}.

\begin{example}
The following diagram is a 
typical
product $xhy \in 1_{14} \,\NB\, 1_{12}$:
$$
\begin{tikzpicture}[anchorbase,scale=.8]
\draw (-1.5,-3) to [out=90,in=-90] (0.5,3);
\draw (0,-3) to [out=90,in=-90] (-1.5,3);
\draw (1.5,-3) to [out=90,in=-70]
(0,0.5) to [out=110,in=-90] (-0.5,3);
\draw (-2,-3) to [out=90,in=-90] (-2.5,3);
\draw (-3,3) to [out=-90,in=180] (0,1.25) to [out=0,in=-90] (3.5,3);
\draw (-2,3) to [out=-90,in=180] (0.5,2) to [out=0,in=-90] (2,3);
\draw (-1,3) to [out=-90,in=180] (-.5,2.5) to [out=0,in=-90] (0,3);
\draw (1.5,3) to [out=-90,in=180] (2,2.5) to [out=0,in=-90] (2.5,3);
\draw (1,3) to [out=-90,in=180] (2,2) to [out=0,in=-90] (3,3);
\draw (1,-3) to [out=90,in=180] (1.5,-2.5) to [out=0,in=90] (2,-3);
\draw (.5,-3) to [out=90,in=180] (1.75,-2) to [out=0,in=90] (3,-3);
\draw (-2.5,-3) to [out=90,in=180] (.5,-1.5) to [out=0,in=90] (2.5,-3);
\draw (-1,-3) to [out=90,in=180] (-.75,-2.6) to [out=0,in=90] (-0.5,-3);
\draw[dashed,red] (-3.2,1) to (3.5,1);
\draw[dashed,red] (-3.2,-1) to (3.5,-1);
\node at (3.1,.5) {$\O_1$};
\node at (3.1,0) {$\O_3$};
\node at (3.1,-.5) {$\O_3$};
\closeddot{-2.98,2.75};\node at (-3.2,2.75) {$\scriptstyle 2$};
\closeddot{-.94,2.75};\node at (-1.15,2.75) {$\scriptstyle 4$};
\closeddot{2.45,2.75};\node at (2.65,2.75) {$\scriptstyle 1$};
\closeddot{1,2.05};\node at (1.15,1.9) {$\scriptstyle 2$};
\closeddot{-2.25,0};\node at (-2.45,0) {$\scriptstyle 1$};
\closeddot{-.96,.5};\node at (-1.18,.5) {$\scriptstyle 9$};
\closeddot{-.88,-.7};\node at (-1.1,-.7) {$\scriptstyle 3$};
\closeddot{-.97,-2.8};\node at (-1.17,-2.8) {$\scriptstyle 2$};
\closeddot{2.99,-2.8};\node at (3.2,-2.8) {$\scriptstyle 3$};
\closeddot{1.97,-2.8};\node at (2.17,-2.8) {$\scriptstyle 1$};
\node[rectangle,rounded corners,draw,fill=red!15!white,inner sep=4pt] at (-5.5,0) {$h$};
\node[rectangle,rounded corners,draw,fill=red!15!white,inner sep=4.3pt] at (-5.5,2.1) {$x$};
\node[rectangle,rounded corners,draw,fill=red!15!white,inner sep=4pt] at (-5.5,-2.1) {$y$};
\end{tikzpicture}
$$
\end{example}

\begin{example}\label{reader}
Equivalence classes of
undotted reduced string diagrams of type X
with $f$ cups and $n$ propagating strings
are in bijection with the set
of chord diagrams with $f$ free chords and $n$ tethered ones as discussed in \cref{chorddiagrams}.
For example, the chord diagram
\cref{orderfood} corresponds to the 
string diagram
$$
\begin{tikzpicture}
	\draw[-] (0,-1) to (0,1);
	\draw[-] (1,-1) to (1,1);
	\draw[-] (-1,-1) to (-1.5,1);
	\draw[-] (-2.5,1) to [out=-90,in=180] (-.25,-.5) to [out=0,in=-90] (2,1);
	\draw[-] (-2,1) to [out=-90,in=180] (-1.25,.5) to [out=0,in=-90] (-.5,1);
\draw[-] (-1,1) to [out=-90,in=180] (0.25,0) to [out=0,in=-90] (1.5,1);
\draw[-] (-3,1) to [out=-90,in=180] (-1.15,0) to [out=0,in=-90] (.5,1);
\node at (-1,-1.17) {$\scriptstyle 1$};
\node at (0,-1.17) {$\scriptstyle 2$};
\node at (1,-1.17) {$\scriptstyle 3$};
\end{tikzpicture}
$$
We hope the bijection here is apparent;
it is similar to the bijection described in the proof of \cref{toinfinityandbeyond} but now the propagating strings become chords that are 
tethered to the bottom node.
\end{example}

\begin{theorem}\label{tbase}
The products
$x h y$ for $(x,h,y) \in \bigcup_{n \in \N}
\X(n) \times \H(n) \times \Y(n)$
give a graded triangular basis for $\NB$
in the sense of \cite[Def.~1.1]{GTB}
(taking the sets $\mathbf{I}$, $\mathbf{S}$ and $\Lambda$
there all to be equal to $\N$ ordered in the natural way).
\end{theorem}

\begin{proof}
We can choose the set $\RSD(a\times b)$
in \cref{basisthm} so that it consists
of the products $xhy$ for
$(x,h,y) \in \bigcup_{n \in \N}
\X(a,n) \times \HH(n) \times \Y(n,b)$.
These give a basis for $1_a \NB 1_b$
as a free $\Gamma$-module. Since
elements of 
$\H(n)$ are elements of
$\HH(n)$ multiplied by basis elements of $\Gamma$,
it follows that the
products $xhy$ for
$(x,h,y) \in \bigcup_{n \in \N}
\X(a,n) \times \H(n) \times \Y(n,b)$
give a linear basis for $1_a \NB 1_b$.
The remaining axioms of graded triangular 
basis are trivial to check.
\end{proof}

\subsection{Standard modules and BGG reciprocity}
\cref{tbase} is significant because it means we can apply the general theory 
developed in \cite{GTB}.
We recall some of the basic constructions made
there. For $n \in \N$, let $\NB_{\geq n}$
be the quotient of $\NB$ by the
two-sided ideal 
generated by $1_m\:(m \not\geq n)$.
Writing $\bar u$ for the canonical image of
$u \in \NB$ in the quotient $\NB_{\geq n}$,
we let $\NB_n := \bar 1_n \NB_{\geq n} \bar 1_n$. This is a unital graded $\Gamma$-algebra
with basis $\bar h\:\big(h \in \HH(n)\big)$
as a free $\Gamma$-module.
These $\bar h$ are the
usual diagrams for elements of a basis of 
the nil-Hecke algebra
associated to the symmetric group.
In fact, 
$\NB_n$ is precisely this nil-Hecke algebra over the ground ring $\Gamma$.
Put somewhat informally, this follows 
because the following ``local relations'' hold:
\begin{align}
\begin{tikzpicture}[anchorbase,scale=.8]
	\draw[-] (0.28,0) to[out=90,in=-90] (-0.28,.6);
	\draw[-] (-0.28,0) to[out=90,in=-90] (0.28,.6);
	\draw[-] (0.28,-.6) to[out=90,in=-90] (-0.28,0);
	\draw[-] (-0.28,-.6) to[out=90,in=-90] (0.28,0);
\end{tikzpicture}
&=0,
&\begin{tikzpicture}[anchorbase,scale=.8]
	\draw[-] (0.45,.6) to (-0.45,-.6);
	\draw[-] (0.45,-.6) to (-0.45,.6);
        \draw[-] (0,-.6) to[out=90,in=-90] (-.45,0);
        \draw[-] (-0.45,0) to[out=90,in=-90] (0,0.6);
\end{tikzpicture}
&=
\begin{tikzpicture}[anchorbase,scale=.8]
	\draw[-] (0.45,.6) to (-0.45,-.6);
	\draw[-] (0.45,-.6) to (-0.45,.6);
        \draw[-] (0,-.6) to[out=90,in=-90] (.45,0);
        \draw[-] (0.45,0) to[out=90,in=-90] (0,0.6);
\end{tikzpicture}\:,&
\begin{tikzpicture}[anchorbase,scale=.8]
	\draw[-] (0.4,.6) to (-0.4,-.6);
	\draw[-] (0.4,-.6) to (-0.4,.6);
 \closeddot{-0.22,0.3};
\end{tikzpicture}
-
\begin{tikzpicture}[anchorbase,scale=.8]
	\draw[-] (0.4,.6) to (-0.4,-.6);
	\draw[-] (0.4,-.6) to (-0.4,.6);
     \closeddot{.2,-0.3};
\end{tikzpicture}
&=
\begin{tikzpicture}[anchorbase,scale=.8]
 	\draw[-] (.2,-.6) to (.2,.6);
	\draw[-] (-.3,-.6) to (-0.3,.6);
\end{tikzpicture}
=\begin{tikzpicture}[anchorbase,scale=.8]
	\draw[-] (0.4,.6) to (-0.4,-.6);
	\draw[-] (0.4,-.6) to (-0.4,.6);
 \closeddot{-0.2,-0.3};
\end{tikzpicture}
-
\begin{tikzpicture}[anchorbase,scale=.8]
	\draw[-] (0.4,.6) to (-0.4,-.6);
	\draw[-] (0.4,-.6) to (-0.4,.6);
     \closeddot{.2,0.3};
\end{tikzpicture}
\:.
\end{align}
These are derived easily from the defining
relations \cref{rels1,rels4,rels7}, noting that the
final cup/cap terms in \cref{rels4,rels7} become 0 
in the quotient algebra.
Because of this term, the nil-Hecke algebra $\NB_n$
is {\em not} a subalgebra of
 $\NB$---one really
 does need to pass first to the quotient
$\NB_{\geq n}$.
In proper algebraic language, $\NB_n$ is the $\Gamma$-algebra generated
by $x_1,\dots,x_n$ all of degree 2 and $\tau_1,\dots,\tau_{n-1}$ all of degree $-2$,
with $\tau_i$ and $x_i$ denoting the 
crossing of the $i$th and $(i+1)$th strings and the dot on the $i$th string, respectively
(numbering strings by $1,\dots,n$
from left to right).
A complete set of relations is
\begin{align}
x_i x_j &= x_j x_i,\\
\tau_i^2 &= 0,\label{relA}\\
\tau_i \tau_j &= \tau_j \tau_i\quad\text{for $|i-j| > 1$},\label{relB}\\
\tau_i \tau_{i+1}\tau_i &= \tau_{i+1}\tau_i \tau_{i+1},\label{relC}\\
x_{i} \tau_i - \tau_i x_{i+1} &= 1=
\tau_i x_i - x_{i+1} \tau_i.\label{lastrel}
\end{align}
One possible basis for $\NB_n$
as a free graded $\Gamma$-module
is given by 
\begin{equation}\label{heckebasisthm}
x_1^{r_1} \cdots x_n^{r_n}
\tau_w \qquad\left(w \in S_n, r_1,\dots,r_n \geq 0 \right)
\end{equation}
Here, $\tau_w$ is the element of $\NB_n$ defined by multiplying
the generators $\tau_i$ according to some reduced expression of $w$.
Recall also that the {\em center}
of the nil-Hecke algebra $\NB_{n}$
is the algebra 
\begin{equation}\label{thecenter}
\ZZ_n := \Gamma[x_1,\dots,x_{n}]^{S_n}
\subseteq \NB_n
\end{equation}
of symmetric polynomials over $\Gamma$.

The {\em polynomial representation} of $\NB_n$
is the graded $\NB_n$-module
$\Gamma[x_1,\dots,x_n]$, with 
$x_i$ acting in the obvious way by multiplication and
$\tau_i$ acting as the Demazure operator
\begin{equation}
\tau_i f = \frac{f-s_i(f)}{x_i - x_{i+1}},
\end{equation}
using $s_i$ for the basic transposition $(i\:\:i\!+\!1) \in S_n$.
Incorporating also a grading shift,
we obtain the indecomposable projective
graded $\NB_n$-module
$P_n(n) := q^{-\frac{1}{2}n(n-1)} \Gamma[x_1,\dots,x_n]$.
Using \cref{heckebasisthm}, it is easy to see that
$P_n(n)$ is generated by the polynomial $u_n := 1$
(which is of degree $-\frac{1}{2}n(n-1)$ due to the grading shift)
subject just to the relations
that $\tau_i u_n = 0$ for $i=1,\dots,n-1$.

Let $L_n(n) := \head P_n(n)$.
This is an irreducible 
graded $\NB_n$-module,
and every irreducible graded $\NB_n$-module
is isomorphic to $L_n(n)$ up to a grading shift.
Writing $\bar u_n$ for the image of $u_n$ in the quotient $L_n(n)$, the monomials
\begin{equation}\label{LBasis}
x_1^{r_1} \cdots x_n^{r_n} \bar u_n \quad
(0 \leq r_i \leq n-i)
\end{equation}
give a homogeneous linear basis for $L_n(n)$. In particular,
\begin{equation}\label{Ldim}
\dim_q L_n(n) = [n]!.
\end{equation}
It is well known that
\begin{equation}\label{hwvector}
\tau_{w_n} (x_1^{n-1} x_2^{n-2} \cdots x_{n-1})
\bar u_n = \bar u_n.
\end{equation}
Note also that any homogeneous 
element in $\ZZ_n$ of 
positive degree acts as 0 on $\bar u_n$, as does any $\tau_i\:(1 \leq i \leq n-1)$.
This is a full set of relations for $L_n(n)$.

We identify $\NB_{\geq n}\gMOD$
with a subcategory of $\NB\gMOD$ in the obvious way.
Truncation with the idempotent $\bar 1_n$
defines a quotient functor
$\jmath^n:\NB_{\geq n}\gMOD\rightarrow \NB_n\gMOD$.
This has left and right adjoints
called the
{\em standardization}
and {\em costandardization functors}:
\begin{align}
\jmath^n_!:=
\NB_{\geq n} \bar 1_n \otimes_{\NB_n}-&:\NB_n\gMOD \rightarrow \NB\gMOD,\\
\jmath^n_*:=
\bigoplus_{m \geq n}
\Hom_{\NB_n}(\bar 1_n \NB_{\geq n} 1_m,
-)&:
\NB_n\gMOD\rightarrow \NB\gMOD.
\end{align}
The following lemma implies that both of these functors are exact; see also
\cite[Lem.~4.1]{GTB}.

\begin{lemma}\label{exact}
For $n \in \N$,
$\NB_{\geq n} \bar 1_n$
is free as a right $\NB_n$-module with 
basis $\bar x\:(x \in \X(n))$,
and
$\bar 1_n \NB_{\geq n}$
is free as a left $\NB_n$-module with 
basis $\bar y\:(y \in \Y(n))$.
\end{lemma}

\begin{proof}
This is an instance of 
\cite[(4.4)--(4.5)]{GTB}.
\end{proof}

For $n \in \N$, we define
the {\em standard} and {\em proper standard modules} for $\NB$ to be the induced modules
\begin{align}\label{stdpropstd}
\Delta(n) &:= \jmath^n_! P_n(n),&
\bar\Delta(n) &:= \jmath^n_! L_n(n).
\end{align}
These are cyclic graded $\NB$-modules
generated
by the vectors $v_n := 1 \otimes u_n$ and $\bar v_n := 1 \otimes \bar u_n$, respectively.
Since we have in hand a basis for $L_n(n)$,
\cref{exact} implies that the following vectors
give a linear basis for
 $\bar\Delta(n)$:
\begin{align}\label{properstandardbasis}
x (x_1^{r_1}\cdots x_{n}^{r_n}) \bar v_n\qquad
&\left(x \in \X(n)\text{ and }r_1,\dots,r_n\text{ with }0 \leq r_i \leq n-i\text{ for each }i\right).\end{align}
In particular, the lowest weight space
$1_n \bar\Delta(n)$ is naturally identified with
$L_n(n)$.
Vectors in $\bar\Delta(n)$
can be represented diagrammatically by putting
$\bar v_n$
into a labeled node at the bottom, with the left action of $\NB$
being by attaching diagrams to the $n$ strings at the top of that node. For example,
the following is a vector in $1_m \bar\Delta(n)$
for any $u \in 1_m \NB 1_n$:
\begin{equation}\label{piccy}
\begin{tikzpicture}[anchorbase,scale=1.4]
\draw[-] (-0.3,0.15)--(-0.3,0.8);
\draw[-] (0.3,0.15)--(0.3,0.8);
\draw[-] (-0.4,0.6)--(-0.4,1.2);
\draw[-] (0.4,0.6)--(0.4,1.2);
\node at (-0.15,0.37) {$\cdot$};
\node at (0,0.37) {$\cdot$};
\node at (0.15,0.37) {$\cdot$};
\node at (-0.2,1.05) {$\cdot$};
\node at (0,1.05) {$\cdot$};
\node at (0.2,1.05) {$\cdot$};
\node[rectangle,rounded corners,draw,fill=blue!15!white,inner sep=4pt] at (0,0) {$\hspace{3.5mm}\scriptstyle
\bar v_{n}\hspace{3.5mm}$};
\node[rectangle,rounded corners,draw,fill=blue!15!white,inner sep=4pt] at (0,0.7) {$\hspace{5mm}\scriptstyle u\hspace{5mm}$};
\end{tikzpicture}
\end{equation} 
It is clear this vector is 0 if $u$ has some
$\O_r\:(r > 0)$ on its right boundary.
In view of \cref{rels13}, this is also true if $u$ has some 
$\O_r\:(r > 0)$
on its left boundary.

\begin{lemma}\label{twocents}
We have that $\End_{\NB}(\Delta(n))
\cong \ZZ_n$ and $\End_{\NB}(\bar\Delta(n))
\cong \kk$.
\end{lemma}

\begin{proof}
The homomorphism from $\ZZ_n$
to $\End_{\NB}(\Delta(n))$ defined by its action on the lowest weight space
$1_n \Delta(n) \cong P_n(n)$
is an isomorphism because
$$
\End_{\NB}(\Delta(n))
\cong \Hom_{\NB_{\geq n}}(\jmath^n_! P_n(n),
\jmath^n_! P_n(n))
\cong 
\Hom_{\NB_n}(P_n(n), \jmath^n \jmath^n_! P_n(n))
\cong \End_{\NB_n}(P_n(n))
\cong \ZZ_n.
$$
The argument for $\bar\Delta(n)$ 
is similar, reducing to Schur's Lemma \cref{schurslemma}.
\end{proof}

Let $I_n(n)$ be the injective hull of $L_n(n)$ in $\NB_n\gMod$, which exists by the discussion in \cref{blah}.
There are also the costandard and proper costandard modules
\begin{align}
\nabla(n) &:= \jmath^n_* I_n(n),&
\bar\nabla(n) &:= \jmath^n_* L_n(n).
\end{align}
We will not use these so often, but note that they can also be obtained
from $\Delta(n)$ and $\bar\Delta(n)$, respectively,
by applying the
contravariant graded functor
\begin{equation}\label{nameme}
\circledast:\NB\gMOD \rightarrow \NB\gMOD
\end{equation}
which takes a graded module $V
 = \bigoplus_{n \in \N} \bigoplus_{d \in \Z} 1_n V_d$
to the graded dual $V^\circledast = 
\bigoplus_{n \in \N}
\bigoplus_{d \in \Z} (1_n V_{-d})^*$
viewed as a graded $\NB$-module so that
$(a f)(v) := f(\tT(a)v)$ for $a \in \NB, f \in V^\circledast$ and $v \in V$,
where
$\tT: \NB \rightarrow \NB$ is the $\Gamma$-algebra anti-automorphism 
arising from \cref{T}.
The proof of this assertion, i.e.,
\begin{align}\label{dualywualy}
\nabla(n) &\cong \Delta(n)^\circledast,
&
\bar\nabla(n)&\cong\bar\Delta(n)^\circledast,
\end{align}
follows from the general discussion of duality in
\cite[Sec~5]{GTB}, specifically, the formula (5.3) there. One just needs to note that $\tT$ fixes the idempotents $1_n\:(n \in \N)$, hence, it descends to an anti-automorphism $\tT_n:\NB_n \rightarrow \NB_n$
fixing the generators
$x_1,\dots,x_n,\tau_1,\dots,\tau_{n-1}$.
Moreover, the irreducible $\NB_n$-module
$L_n(n)$ is self-dual with respect to the resulting duality $\circledast$
on $\NB_n\gMOD$. This last statement is clear because 
$\dim_q L_n(n)$ is invariant under the bar involution by \cref{Ldim}, and
$L_n(n)$ is the unique irreducible graded left $\NB_n$-module of this 
graded dimension.

For the basic notions of
{\em $\Delta$-flags}, {\em $\bar\Delta$-flags},
{\em $\nabla$-flags} and {\em $\bar\nabla$-flags}, we refer to \cite[Def.~6.3, Def.~6.4]{GTB}.
In particular, a $\Delta$-flag in a graded $\NB$-module $V$
is a (finite) graded filtration
$0 = V_0 \subseteq V_1 \subseteq \cdots \subseteq V_m$
such that $V_i / V_{i-1}\cong
\Delta(n_i)^{\oplus f_i}$
for distinct $n_1,\dots,n_m \in \N$
and $f_i \in \N\lround q^{-1} \rround$.
The notion of 
a $\bar\Delta$-flag is similar, except that the sections of the filtration are {\em $\bar\Delta$-layers}, that is, 
$V_i / V_{i-1} \cong \jmath^{n_i}_! U_i$ for a graded $\NB_{n_i}$-module $U_i$ which is locally finite dimensional and bounded below.
Multiplicities in these four types of filtration are denoted $(V:\Delta(n))_q$,
$(V:\bar \Delta(n))_q$,
$(V:\nabla(n))_q$ and
$(V:\bar \nabla(n))_q$.
For example,
the standard module $\Delta(n)$ has a $\bar\Delta$-flag with the multiplicities
\begin{equation}\label{anothereclair}
(\Delta(n):\bar\Delta(n))_q = [P_n(n):L_n(n)]_q
= \frac{\dim_{q} \Gamma}{(1-q^{2})(1-q^{4}) \cdots (1-q^{2n})}
\end{equation}
and $(\Delta(n):\bar\Delta(m))_q = 0$ for $m \neq n$.
This follows from exactness of $\jmath^n_!$
and the well-known representation theory of $\NB_n$.

Now we can formulate the fundamental theorem
about the structure of $\NB\gMOD$. It follows
by an application of the general theory developed in \cite{GTB},
specifically,
\cite[Th.~4.3, Sec.~5, Cor.~8.4]{GTB},
and is analogous to the basic structural results about Verma and dual Verma modules in Lie theory.

\begin{theorem}\label{wagner}
The following properties hold:
\begin{enumerate}
\item
The standard module $\Delta(n)$ has a unique irreducible
graded quotient $L(n)$.
Also, $L(n)^\circledast \cong L(n)$,
so that $L(n)$ is also the unique irreducible
graded submodule of $\nabla(n)$. 
\item
The 
$\NB$-modules
$L(n)\:(n \in \N)$
give a complete set of irreducible graded $\NB$-modules
up to isomorphism and grading shift.
\item
Let $P(n)$ be the projective cover of $L(n)$
in $\NB\gMod$ and $I(n) \cong P(n)^\circledast$
be its injective hull.
Then $P(n)$ has a $\Delta$-flag
and $I(n)$ has a $\nabla$-flag, with multiplicities satisfying the usual graded
BGG reciprocity formulae
\begin{align*}
(P(n):\Delta(m))_q &= [\bar\Delta(m):L(n)]_q=
[\bar\nabla(m):L(n)]_{q^{-1}}
=(I(n):\nabla(m))_{q^{-1}}\in\N\lround q\rround
\end{align*}
for all $m, n \in \N$.
These multiplicities are $1$ if $m=n$
and 0 unless $m \leq n$.
\end{enumerate}
\end{theorem}

We denote the 
canonical image of
$v_n$ in the irreducible quotient
$L(n)$ of $\Delta(n)$ by $\tilde v_n$.
Vectors in $L(n)$ can be denoted diagrammatically
just like in \cref{piccy}
putting $\tilde v_n$
into the node
at the bottom of the diagram
instead of $\bar v_n$.
Again, the lowest weight space $1_n L(n)$
is naturally identified with the $\NB_n$-module
$L_n(n)$.

\cref{wagner} gives a classification of irreducible graded left $\NB$-modules via their lowest weights. The proof just explained is completely independent of any of the results from \cref{sect4}.
It follows that the modules $P(n)\:(n \geq 0)$
defined in \cref{wagner}(3) give a complete set of pairwise inequivalent indecomposable graded projective left $\NB$-modules.
Such a classification was already established in \cref{therealthing} by a more sophisticated method involving \cref{id3,toinfinityandbeyond}.
The following shows that the two
approaches are consistent with each other:

\begin{lemma}
For $n \geq 0$, the graded module
$P(n)$ defined in \cref{wagner}(3), that is, the projective cover of $L(n)$,
is isomorphic to the graded module
denoted $P(n)$ in the previous section, that is, $q^{\frac{1}{2}n(n-1)} \NB\,\e_n$.
\end{lemma}

\begin{proof}
Since $q^{\frac{1}{2}n(n-1)} \NB\, \e_n$
is an indecomposable projective graded module
by \cref{therealthing},
it suffices to prove that
$$
\Hom_{\NB}\big(q^{\frac{1}{2}n(n-1)} \NB\,\e_n, \;L(n)\big)_0 
\cong \e_n L(n)_{\frac{1}{2}n(n-1)} \neq 0.
$$
This follows because
$(x_1^{n-1} x_2^{n-2}\cdots x_{n-1}) \tilde v_n
\in L(n)$
is a non-zero vector 
of degree $\frac{1}{2}n(n-1)$
such that $\e_n (x_1^{n-1} x_2^{n-2}\cdots x_{n-1}) \tilde v_n = (x_1^{n-1} x_2^{n-2}\cdots x_{n-1}) \tilde v_n$,
as follows from the definition \cref{endef}
of the idempotent $\e_n$
together with \cref{hwvector}.
\end{proof}

\begin{remark}
For convenience, we have worked with
the natural total ordering on $\N$. 
However, the basis in \cref{tbase} is in fact a graded triangular basis with respect to the {\em partial ordering} $\unlhd$ on $\N$
defined by $m \unlhd n
\Leftrightarrow n-m \in 2 \N$;
this is clear since $\X(a,n)$ and $\Y(n,a)$
are empty unless $a \equiv n \pmod{2}$.
Everything established so far
is also true for this order. 
In particular, both $0$ and $1$ are minimal
with respect to $\unlhd$, so by \cref{wagner}(3) we have that
$P(0)=\Delta(0)$ and $P(1)=\Delta(1)$.
\end{remark}

\subsection{The projective functor \texorpdfstring{$B$}{} preserves good filtrations}\label{pfunctor}
Recall the endofunctor $B$ of
$\NB\gMOD$ introduced in \cref{pfun}.
Using the construction \cref{Ind}, 
it can be defined equivalently
as the induction functor
$\Ind_{B \star-}$
where $B \star-:\cNB_t \rightarrow \cNB_t$
is the graded functor defined by tensoring with $B$. This follows easily from the definitions; see \cite[Lem.~2.4]{BV} for details in a similar situation.
In fact, we can go a step further to
make $\NB\gMOD$ into a strict graded
$\cNB_t$-module category, i.e., 
there is a strict 
graded monoidal functor
$\mu$ from $\cNB_t$ to the strict graded monoidal category 
$\gEND(\NB\gMOD)$
consisting of graded endofunctors and graded natural transformations.
This takes the generating object $B$ of $\cNB_t$ to the graded endofunctor 
$\Ind_{B\star-}$ and
the generating morphisms
$\txtdot\,$,
$\txtcrossing\,$,
$\,\txtcap\,$
and $\,\txtcup\,$
to the graded natural transformations
$\Ind_{\!\!\smalltxtdot\star-},
\Ind_{\,\smalltxtcrossing\star-},
\Ind_{\;\smalltxtcup\,\star-}$
and $\Ind_{\;\smalltxtcap\,\star-}$, respectively.
Notice we have switched the cap and the cup here;
this is the usual price for choosing to work with left modules rather than right modules---we are using the contravariant Yoneda embedding.

\begin{lemma}\label{snow}
The functor 
$\Ind_{B\star-}:\NB\gMOD\rightarrow \NB\gMOD$ is isomorphic to the restriction functor
$\Res_{B\star-}:\NB\gMOD\rightarrow \NB\gMOD$.
The isomorphism can be chosen so that it intertwines the endomorphism
$\Ind_{\!\smalltxtdot\star-}:\Ind_{B\star-}\Rightarrow\Ind_{\!B\star-}$ with
$-\Res_{\!\smalltxtdot\star-}:\Res_{B\star-}\Rightarrow \Res_{B\star-}$.
\end{lemma}

\begin{proof}
The functor $\Ind_{B\star-}$
is defined by tensoring with the bimodule
$\NB 1_{B\star-}$
and the functor $\Res_{B\star-}$
is defined by tensoring with the bimodule
$1_{B\star-} \NB$.
The functors are isomorphic because
there is a graded $(\NB,\NB)$-bimodule isomorphism
$\phi:1_{B\star-} \NB
\stackrel{\sim}{\rightarrow}
\NB 1_{B\star-}
$ such that
\begin{align}\label{dualitypic}
\phi\left(
\begin{tikzpicture}[anchorbase,scale=1.4]
\draw[-] (-0.4,-0.15)--(-0.4,-0.5);
\draw[-] (0.4,-0.15)--(0.4,-0.5);
\draw[-] (-0.4,0.15)--(-0.4,0.5);
\draw[-] (-0.3,0.15)--(-0.3,0.5);
\draw[-] (0.4,0.15)--(0.4,0.5);
\node at (-0.1,0.35) {$\cdot$};
\node at (0.05,0.35) {$\cdot$};
\node at (0.2,0.35) {$\cdot$};
\node at (-0.2,-0.4) {$\cdot$};
\node at (-0,-0.4) {$\cdot$};
\node at (0.2,-0.4) {$\cdot$}; 
\node[rectangle,rounded corners,draw,fill=blue!15!white,inner sep=4pt] at (0,0) {$\hspace{5mm}\scriptstyle u\hspace{5mm}$};
\end{tikzpicture}\right)&
= \begin{tikzpicture}[anchorbase,scale=1.4]
\draw[-] (-0.4,-0.15)--(-0.4,-0.5);
\draw[-] (0.4,-0.15)--(0.4,-0.5);
\draw[-] (-0.3,0.15)--(-0.3,0.5);
\draw[-] (0.4,0.15)--(0.4,0.5);
\draw[-] (-0.4,0.15)to[out=up,in=right](-0.6,0.5)to[out=left,in=up](-0.8,0.3)to(-0.8,-0.5);
\node at (-0.1,0.35) {$\cdot$};
\node at (0.05,0.35) {$\cdot$};
\node at (0.2,0.35) {$\cdot$};
\node at (-0.2,-0.4) {$\cdot$};
\node at (-0,-0.4) {$\cdot$};
\node at (0.2,-0.4) {$\cdot$}; 
\node[rectangle,rounded corners,draw,fill=blue!15!white,inner sep=4pt] at (0,0) {$\hspace{5mm}\scriptstyle u\hspace{5mm}$};
\end{tikzpicture}\:,
&
\phi^{-1}\left(\begin{tikzpicture}[anchorbase,scale=1.5]
\draw[-] (-0.4,0.1)--(-0.4,0.5);
\draw[-] (0.4,0.1)--(0.4,0.5);
\draw[-] (-0.4,-0.1)--(-0.4,-0.5);
\draw[-] (-0.3,-0.1)--(-0.3,-0.5);
\draw[-] (0.4,-0.1)--(0.4,-0.5);
\node at (-0.1,-0.4) {$\cdot$};
\node at (0.05,-0.4) {$\cdot$};
\node at (0.2,-0.4) {$\cdot$};
\node at (-0.2,0.35) {$\cdot$};
\node at (-0,0.35) {$\cdot$};
\node at (0.2,0.35) {$\cdot$}; 
\node[rectangle,rounded corners,draw,fill=blue!15!white,inner sep=4pt] at (0,0) {$\hspace{5mm}\scriptstyle v\hspace{5mm}$};
\end{tikzpicture}\right)&= \begin{tikzpicture}[anchorbase,scale=1.5]
\draw[-] (-0.4,-0.1)to[out=down,in=right](-0.6,-0.5)to[out=left,in=down](-0.8,-0.3)to(-0.8,0.5);
\draw[-] (-0.4,0.1)--(-0.4,0.5);
\draw[-] (0.4,0.1)--(0.4,0.5);
\draw[-] (-0.3,-0.1)--(-0.3,-0.5);
\draw[-] (0.4,-0.1)--(0.4,-0.5);
\node at (-0.1,-0.4) {$\cdot$};
\node at (0.05,-0.4) {$\cdot$};
\node at (0.2,-0.4) {$\cdot$};
\node at (-0.2,0.35) {$\cdot$};
\node at (-0,0.35) {$\cdot$};
\node at (0.2,0.35) {$\cdot$}; 
\node[rectangle,rounded corners,draw,fill=blue!15!white,inner sep=4pt] at (0,0) {$\hspace{5mm}\scriptstyle v\hspace{5mm}$};
\end{tikzpicture}\:.
\end{align}
Remembering the sign in the nil-Brauer relations \cref{rels4,rels7},
the resulting isomorphism intertwines 
$\Ind_{\!\smalltxtdot\star-}$ with  $-\Res_{\!\smalltxtdot\star-}$.
\end{proof}

From now on, we denote
the endofunctor $\Ind_{B\star-}$
simply by $B$ (as we did in the previous section).
We often 
use $x$ to denote the endomorphism of $B$
defined by $\Ind_{\!\smalltxtdot\star-}$.
The same letter is used to denote elements of $\X(n)$, but we think it is always clear from context which we mean.

\begin{lemma}\label{exactness}
The endofunctor
$B:\NB\gMOD \rightarrow \NB\gMOD$ 
is self-adjoint. Hence, on the Abelian category
$\NB\gMod$, it is exact, cocontinuous, and preserves finitely generated projectives.
Also $B$ commutes with the duality \cref{nameme}, i.e., we have that
$B \circ \circledast \cong \circledast \circ B$.
\end{lemma}

\begin{proof}
\cref{snow} shows
that $B$ is isomorphic to a right adjoint to $B$. Hence, it is self-adjoint.
The fact that $B$ commutes with duality
follows because $\Res_{|\star-}$ clearly does so.
\end{proof}

\begin{lemma}\label{minpolydegree}
For $n \geq 0$, the degree $\beta(n)$ of the minimal polynomial
of $x_{L(n)}:BL(n) \rightarrow BL(n)$
satisfies $\beta(n) \equiv t\pmod{2}$.
\end{lemma}

\begin{proof}
This is a consequence of \cref{lynch}.
In order for this to be applicable, we need
to check that
$\End_{\NB}(L(n)) = \kk$
and $\dim \End_\NB(BL(n)) < \infty$.
The first property follows from \cref{schurslemma}.
The second follows because
$\End_{\NB}(BL(n)) \cong \Hom_{\NB}(B^2 L(n), L(n))$, and the latter space
is finite dimensional because $B^2 L(n)$ is finitely generated
and $L(b)$ is locally finite dimensional.
\end{proof}

Let $\iota_{1,n}:\NB_n \hookrightarrow \NB_{n+1}$
be the (unital)
graded $\Gamma$-algebra homomorphism
mapping 
$x_i \mapsto x_{i+1}$ and $\tau_j \mapsto \tau_{j+1}$.
We denote the restriction of a graded left
(resp., right) $\NB_{n+1}$-module along the homomorphism $\iota_{1,n}$ by
$\iota^*_{1,n} V$ (resp., $V \iota^*_{1,n}$).
Let $(I_{1,n}, R_{1,n})$ be the resulting adjoint pair of induction and restriction functors
between $\NB_n\gMod$ and $\NB_{n+1}\gMod$.
We have that $I_{1,n} = \NB_{n+1} \iota_{1,n}^*
\otimes_{\NB_n}-$
and $R_{1,n} \cong \iota_{1,n}^*NB_{n+1}
\otimes_{\NB_{n+1}}-$.

\begin{lemma}\label{bundy}
The vectors $x_1^r \tau_1 \cdots \tau_{i-1}\:\left(1 \leq i \leq n+1, r \geq 0\right)$
give a basis for
$\iota_{1,n}^* \NB_{n+1}$ 
as a free graded left $\NB_n$-module.
Similarly, the vectors $\tau_{i-1}\cdots \tau_1 x_1^r\:\left(1 \leq i \leq n+1, 
r\geq 0\right)$
give a basis for
$\NB_{n+1} \iota_{1,n}^*$
as a free graded right $\NB_n$-module.
Hence, the functors $I_{1,n}$ and $R_{1,n}$ are exact.
\end{lemma}

\begin{proof}
This is well known.
The first statement follows easily from  \cref{heckebasisthm}, and the second statement may be deduced from the first by applying an anti-automorphism.
\end{proof}

Recall the isomorphism $\jmath:\Ui \stackrel{\sim}{\rightarrow} \f$
from \cref{jy2}.
Since we are favoring standard modules over costandard modules in our
exposition, we need now the twisted version of this map, which is the
isomorphism
\begin{equation}
\tilde\jmath := \psi \circ \jmath \circ \psi^\imath:\Ui
\stackrel{\sim}{\rightarrow} \f.
\end{equation}
The analog of \cref{jy2}(1) for this is
\begin{equation}
\tilde\jmath(bu) = \theta \tilde\jmath(u) +
\widetilde{R}(\jmath(u))\qquad\text{ for }u \in \Ui,
\end{equation}
where $\widetilde{R} := \psi \circ R \circ \psi:\f\rightarrow \f$.
Equivalently, the inverse isomorphism 
$\tilde\jmath^{-1}: \f \stackrel{\sim}{\rightarrow} \Uinone$ has the property
\begin{align} \label{SESalg}
b \tilde\jmath^{-1}(x) = \tilde\jmath^{-1}(\theta x) + \tilde\jmath^{-1}(\widetilde{R}(x)) \qquad \text{ for } x\in \f.
\end{align}
Our next theorem can be interpreted as a categorification of this identity, with $\jmath^n_!$ corresponding to  $\tilde\jmath^{-1}$,
$I_{1,n}\:(n\geq 0)$ corresponding to multiplication by $\theta$,
and the functors $R_{1,n}\:(n > 0)$ corresponding to the map
$\widetilde{R}$. The fact that the restriction functors
$R_{1,n}$ categorify $\widetilde{R}$ was first pointed out
in \cite{KK}.

\begin{theorem} \label{Thm:SES}
For $n \geq 0$,
there is a short exact sequence of functors\footnote{We mean that one obtains
a short exact sequence in $\NB\gMod$
after evaluating 
on any graded left $\NB_n$-module $V$.}
\begin{equation}\label{cardiff}
0 \longrightarrow \jmath^{n-1}_! \circ R_{1,n-1}
\stackrel{\alpha}{\longrightarrow} B 
\circ \jmath^n_!
\stackrel{\beta}{\longrightarrow}
\jmath^{n+1}_!\circ I_{1,n}
\longrightarrow 0,
\end{equation}
interpreting the first term 
as the zero functor in the case $n=0$.
Moreover, letting
$x':R_{1,n}\Rightarrow R_{1,n}$
and $x'':I_{1,n}\Rightarrow I_{1,n}$
be the degree 2 endomorphisms induced by the endomorphisms of the
bimodules
$\iota^*_{1,n} \NB_{n+1}$ 
and
$\NB_{n+1} \iota^*_{1,n}$
defined by left multiplication by $-x_1$
and by right
multiplication by 
$x_1$, respectively,
we have that 
\begin{align}\label{space}
\alpha \circ \left(\jmath^{n-1}_! x'\right) &= \left(x \jmath^n_!\right) \circ \alpha,&
\beta \circ \left(x \jmath^n_!\right) &= \left(\jmath^{n+1}_! x''\right)\circ \beta.
\end{align}
\end{theorem}

\begin{proof}
All three functors appearing in the short exact sequence are defined by tensoring with
certain graded 
$(\NB, \NB_n)$-bimodules:
$\jmath^{n-1}_! \circ R_{1,n-1}$
is tensoring with the bimodule
$\NB_{\geq (n-1)} \bar 1_{n-1} \otimes_{\NB_{n-1}} 
\iota^*_{1,n-1} \NB_n$,
$B \circ \jmath^n_!$ is tensoring with the bimodule
$1_{B\star-} \NB_{\geq n} \bar 1_n$
(here we have used \cref{snow} to realize $B$ as restriction rather than induction),
and
$\jmath^{n+1}_! \circ I_{1,n}$ is tensoring with
$\NB_{\geq (n+1)} \bar 1_{n+1}
\otimes_{\NB_{n+1}} \NB_{n+1} \iota^*_{1,n}$.
In the next two paragraphs, we construct
a short exact sequence of graded bimodules and degree-preserving bimodule homomorphisms:
$$
0 \longrightarrow 
\NB_{\geq (n-1)} \bar 1_{n-1} \otimes_{\NB_{n-1}} 
\iota^*_{1,n-1} \NB_n\stackrel{a}{\longrightarrow}
1_{B\star-} \NB_{\geq n} \bar 1_n\stackrel{b}{\longrightarrow}
\NB_{\geq (n+1)} \bar 1_{n+1}
\otimes_{\NB_{n+1}} \NB_{n+1} \iota^*_{1,n}
\longrightarrow 0.
$$
\cref{exact,bundy} 
imply that 
$\NB_{\geq (n+1)} \bar 1_{n+1} \otimes_{\NB_{n+1}} 
\NB_{n+1} \iota^*_{1,n} \cong \NB_{\geq (n+1)} \bar 1_{n+1} \iota^*_{1,n}$ is free, hence, projective as a
graded right $\NB_{n}$-module.
So this is a split short exact sequence of right $\NB_{n}$-modules.
It follows that it
remains exact
when we apply the functor $-\otimes_{\NB_n} V$
for any left $\NB_n$-module $V$. Thus, we have
constructed the short exact sequence of functors in the statement of the theorem.

To construct the short exact sequence of bimodules, take $m \geq 0$.
We can assume the set $\X(m+1,n)$ 
is chosen to be 
\begin{equation}\label{chilly}
\X(m+1,n) = 
\Bigg\{
\:\begin{tikzpicture}[anchorbase,scale=1.25]
\draw[-] (-0.45,0.5)--(-0.45,-0.5);
\draw[-] (0.4,-0.15)--(0.4,-0.5);
\draw[-] (-0.2,0.15)--(-0.2,0.5);
\draw[-] (-0.2,-0.15)--(-0.2,-0.5);
\draw[-] (0.4,0.15)--(0.4,0.5);
\node at (-0.05,0.35) {$\cdot$};
\node at (0.1,0.35) {$\cdot$};
\node at (0.25,0.35) {$\cdot$};
\node at (-0.05,-0.35) {$\cdot$};
\node at (0.1,-0.35) {$\cdot$};
\node at (0.25,-0.35) {$\cdot$};
\node[rectangle,rounded corners,draw,fill=blue!15!white,inner sep=4.2pt] at (0.1,0) {$\hspace{3mm}\scriptstyle x\hspace{3mm}$};
\end{tikzpicture}\:\:
\:\Bigg|\: x \in \X(m,n-1)
\Bigg\}
\sqcup
\Bigg\{
\:\begin{tikzpicture}[anchorbase,scale=1.25]
\draw[-] (-0.45,0.5)--(-0.45,-0.2)
to [out=-90,in=180] (-.2,-.5)
to [out=0,in=-90] (0,-.05);
\draw[-] (0.4,-0.05)--(0.4,-0.5);
\draw[-] (-0.2,0.2)--(-0.2,0.5);
\draw[ultra thick] (-0.2,-0.05) to [out=-90,in=90] (0.05,-0.5);
\draw[-] (0.4,0.2)--(0.4,0.5);
\node at (-0.05,0.4) {$\cdot$};
\node at (0.1,0.4) {$\cdot$};
\node at (0.25,0.4) {$\cdot$};
\node at (0.12,-0.27) {$\cdot$};
\node at (0.22,-0.27) {$\cdot$};
\node at (0.32,-0.27) {$\cdot$};
\node at (0.08,-.59) {$\stringlabel{i-1}$};
\closeddot{-.09,-.45};
\node at (-.19,-.35) {$\scriptstyle r$};
\node[rectangle,rounded corners,draw,fill=blue!15!white,inner sep=4.2pt] at (0.1,0.1) {$\hspace{3mm}\scriptstyle x\hspace{3mm}$};
\end{tikzpicture}\:\:
\:\Bigg|\:
\begin{array}{l}
x \in \X(m,n+1)\\1\leq i \leq n+1\\r \geq 0\end{array}
\Bigg\}.
\end{equation}
The first set on the right hand side 
here (which should be interpreted as $\varnothing$ in case $n=0$)
gives
the elements of $\X(m+1,n)$ which have a propagating string at the top left boundary point.
The second set gives all remaining elements of
$\X(m+1,n)$. These have a cup 
at the top left boundary point 
with $(i-1)$ propagating strings
between that and its other boundary point
for some $1 \leq i \leq n+1$; these are represented in the diagram by the single thick string labeled $i-1$.
We can assume the set $\H(n+1)$ is chosen to be
\begin{equation}\label{dip}
\H(n+1) = \Bigg\{
\:\begin{tikzpicture}[anchorbase,scale=1.25]
\draw[-] (-0.45,-0.5) to [out=90,in=-90] (-.45,0) to [out=90,in=-90] (-0.1,0.5);
\draw[-] (0.4,-0.15)--(0.4,-0.5);
\draw[ultra thick] (-0.2,0.1) to [out=90,in=-90] (-.45,0.5);
\draw[-] (-0.2,-0.15)--(-0.2,-0.5);
\draw[-] (0.4,0.1)--(0.4,0.5);
\node at (-0.43,.58) {$\stringlabel{i-1}$};
\closeddot{-.45,-.35};
\node at (-.57,-.33) {$\scriptstyle r$};
\node at (-0.05,0.35) {$\cdot$};
\node at (0.1,0.35) {$\cdot$};
\node at (0.25,0.35) {$\cdot$};
\node at (-0.05,-0.4) {$\cdot$};
\node at (0.1,-0.4) {$\cdot$};
\node at (0.25,-0.4) {$\cdot$};
\node[rectangle,rounded corners,draw,fill=blue!15!white,inner sep=4.2pt] at (0.1,-.07) {$\hspace{3mm}\scriptstyle h\hspace{3mm}$};
\end{tikzpicture}\:\:
\:\Bigg|\:
\begin{array}{l}
h \in \H(n)\\
1 \leq i \leq n+1\\
r \geq 0
\end{array}
\Bigg\}.
\end{equation}
In these diagrams, the propagating string 
with the bottom left boundary point
has $(i-1)$ other strings to the left of its other boundary point for some $1 \leq i \leq n+1$.
The vectors $x \otimes h\:\left(x \in \X(m,n-1),
h \in \H(n)\right)$ give a linear basis
for 
$\bar 1_m \NB_{\geq (n-1)} \bar 1_{n-1}
\otimes_{\NB_{n-1}}\iota^*_{1,n-1}
\NB_n$ by \cref{exact} again.
We define an injective linear map 
$a_m:\bar 1_m \NB_{\geq (n-1)} \bar 1_{n-1}
\otimes_{\NB_{n-1}}\iota^*_{1,n-1}
\NB_n
\hookrightarrow \bar 1_{m+1}\NB_{\geq n}\bar 1_n$
on basis vectors by
\begin{align*}
a_m: \begin{tikzpicture}[anchorbase,scale=1.25]
\draw[-] (0.4,-0.05)--(0.4,-0.4);
\draw[-] (-0.2,0.15)--(-0.2,0.55);
\draw[-] (-0.2,-0.05)--(-0.2,-0.4);
\draw[-] (0.4,0.15)--(0.4,0.55);
\node at (-0.05,0.42) {$\cdot$};
\node at (0.1,0.42) {$\cdot$};
\node at (0.25,0.42) {$\cdot$};
\node at (-0.05,-0.27) {$\cdot$};
\node at (0.1,-0.27) {$\cdot$};
\node at (0.25,-0.27) {$\cdot$};
\node[rectangle,rounded corners,draw,fill=blue!15!white,inner sep=4.2pt] at (0.1,.1) {$\hspace{3mm}\scriptstyle x\hspace{3mm}$};
\end{tikzpicture}
\:\otimes\:
\begin{tikzpicture}[anchorbase,scale=1.25]
\draw[-] (-0.45,-1.1)--(-0.45,-0.7);
\draw[-] (0.45,-1.1)--(0.45,-0.7);
\draw[-] (-0.45,-0.1)--(-0.45,-0.4);
\draw[-] (0.45,-0.1)--(0.45,-0.4);
\node at (-0.25,-0.25) {$\cdot$};
\node at (0,-0.25) {$\cdot$};
\node at (0.25,-0.25) {$\cdot$};
\node at (-0.25,-.96) {$\cdot$};
\node at (0,-.96) {$\cdot$};
\node at (0.25,-.96) {$\cdot$};
\node[rectangle,rounded corners,draw,fill=blue!15!white,inner sep=4.2pt] at (0,-.6) {$\hspace{5mm}\scriptstyle h\hspace{5mm}$};
\end{tikzpicture}
&\mapsto
 \begin{tikzpicture}[anchorbase,scale=1.25]
\draw[-] (-0.45,-1.05)--(-0.45,-0.7);
\draw[-] (0.45,-1.05)--(0.45,-0.7);
\draw[-] (-0.45,0.5)--(-0.45,-0.4);
\draw[-] (0.4,-0.05)--(0.4,-0.4);
\draw[-] (-0.2,0.15)--(-0.2,0.5);
\draw[-] (-0.2,-0.05)--(-0.2,-0.4);
\draw[-] (0.4,0.15)--(0.4,0.5);
\node at (-0.05,0.4) {$\cdot$};
\node at (0.1,0.4) {$\cdot$};
\node at (0.25,0.4) {$\cdot$};
\node at (-0.05,-0.25) {$\cdot$};
\node at (0.1,-0.25) {$\cdot$};
\node at (0.25,-0.25) {$\cdot$};
\node at (-0.25,-0.95) {$\cdot$};
\node at (0,-0.95) {$\cdot$};
\node at (0.25,-0.95) {$\cdot$};
\node[rectangle,rounded corners,draw,fill=blue!15!white,inner sep=4.2pt] at (0.1,.1) {$\hspace{3mm}\scriptstyle x\hspace{3mm}$};
\node[rectangle,rounded corners,draw,fill=blue!15!white,inner sep=4.2pt] at (0,-.6) {$\hspace{5mm}\scriptstyle h\hspace{5mm}$};
\end{tikzpicture}
\end{align*}
for $x \in \X(m,n-1), h \in \H(n)$.
The image of $a_m$ is the subspace 
of $\bar 1_{m+1}\NB_{\geq n} \bar 1_n$
with basis given by the vectors 
$x h\:\left(x \in \X(m+1,n), h \in \H(n)\right)$,
i.e., the basis vectors with 
$x$ in the first set
on the right hand side of \cref{chilly}.
We define $b_m:\bar 1_{m+1} \NB_{\geq n}
\bar 1_n\twoheadrightarrow \bar 1_m\NB_{\geq(n+1)}
\bar 1_{n+1}\otimes_{\NB_{n+1}}
\NB_{n+1}\iota^*_{1,n}$
to be the surjective 
linear map that is 0 on
these basis vectors
and is defined on the remaining basis vectors $x h\:(x \in \X(m+1,n),
h \in \H(n))$ for $x$ in the second set
on the right hand side of \cref{chilly}
by
$$
b_m:
\begin{tikzpicture}[anchorbase,scale=1.25]
\draw[-] (-0.45,0.5)--(-0.45,-0.2)
to [out=-90,in=180] (-.2,-.5)
to [out=0,in=-90] (0,-.05);
\draw[-] (0.4,-0.05) to [out=-90,in=90] (0.6,-0.7);
\draw[-] (-0.2,0.2)--(-0.2,0.5);
\draw[ultra thick] (-0.2,-0.05) to [out=-90,in=90] (0.05,-0.4) to [out=-90,in=90] (0,-.75);
\draw[-] (0.4,0.2)--(0.4,0.5);
\draw[-] (0,-1.05)--(0,-1.35);
\draw[-] (0.6,-1.05)--(0.6,-1.35);
\node at (-0.05,0.4) {$\cdot$};
\node at (0.1,0.4) {$\cdot$};
\node at (0.25,0.4) {$\cdot$};
\node at (.19,-0.5) {$\cdot$};
\node at (0.32,-0.5) {$\cdot$};
\node at (0.45,-0.5) {$\cdot$};
\node at (.14,-1.27) {$\cdot$};
\node at (0.32,-1.27) {$\cdot$};
\node at (0.5,-1.27) {$\cdot$};
\node at (0.25,-.31) {$\stringlabel{i-1}$};
\closeddot{-.09,-.45};
\node at (-.19,-.35) {$\scriptstyle r$};
\node[rectangle,rounded corners,draw,fill=blue!15!white,inner sep=4.2pt] at (0.11,0.1) {$\hspace{3.5mm}\scriptstyle x\hspace{3.5mm}$};
\node[rectangle,rounded corners,draw,fill=blue!15!white,inner sep=4.2pt] at (.3,-.9) {$\hspace{3mm}\scriptstyle h\hspace{3mm}$};
\end{tikzpicture}\:\:
\mapsto
\:\:
\begin{tikzpicture}[anchorbase,scale=1.25]
\draw[-] (0.4,0.7)--(0.4,1.15);
\draw[-] (-0.35,0.7)--(-0.35,1.15);
\draw[-] (0.4,0.2)--(0.4,.5);
\draw[-] (-0.35,0.2)--(-0.35,.5);
\node at (-0.19,0.32) {$\cdot$};
\node at (0.03,0.32) {$\cdot$};
\node at (0.25,0.32) {$\cdot$};
\node at (-0.19,1) {$\cdot$};
\node at (0.03,1) {$\cdot$};
\node at (0.25,1) {$\cdot$};
\node[rectangle,rounded corners,draw,fill=blue!15!white,inner sep=4.2pt] at (0.02,0.67) {$\hspace{3.5mm}\scriptstyle x\hspace{3.5mm}$};
\end{tikzpicture}
\:\:\otimes
\begin{tikzpicture}[baseline=-1mm,scale=1.25]
\draw[-] (-0.45,-0.5) to [out=90,in=-90] (-.45,0) to [out=90,in=-90] (-0.1,0.5);
\draw[-] (0.4,-0.15)--(0.4,-0.5);
\draw[ultra thick] (-0.2,0.1) to [out=90,in=-90] (-.45,0.5);
\draw[-] (-0.2,-0.15)--(-0.2,-0.5);
\draw[-] (0.4,0.1)--(0.4,0.5);
\node at (-0.44,.6) {$\color{blue}\scriptstyle i-1$};
\closeddot{-.45,-.35};
\node at (-.57,-.33) {$\scriptstyle r$};
\node at (-0.05,0.35) {$\cdot$};
\node at (0.1,0.35) {$\cdot$};
\node at (0.25,0.35) {$\cdot$};
\node at (-0.05,-0.4) {$\cdot$};
\node at (0.1,-0.4) {$\cdot$};
\node at (0.25,-0.4) {$\cdot$};
\node[rectangle,rounded corners,draw,fill=blue!15!white,inner sep=4.2pt] at (0.1,-.07) {$\hspace{3mm}\scriptstyle h\hspace{3mm}$};
\end{tikzpicture}
$$
for $x \in \X(m,n+1),
1\leq i \leq n+1$ and $r \geq 0$.
In view of \cref{dip,exact}, 
the image vectors here are a basis for
$\bar 1_m\NB_{\geq(n+1)}
\bar 1_{n+1}\otimes_{\NB_{n+1}}
\NB_{n+1}\iota^*_{1,n}$.
Now we have that $a_m$ is injective,
$b_m$ is surjective and $\im a_m = \ker b_m$.
Then we define $a := \bigoplus_{m \geq 0} a_m$ and $b := \bigoplus_{m \geq 0} b_m$.
This gives the linear maps
in the short exact sequence that we are after, and we have checked the exactness.

Next, we show that
$a$ and $b$ are graded bimodule homomorphisms.
The map $a$ is given equivalently by
multiplication $\NB_{\geq n} \bar 1_{n-1} \otimes_{\NB_{n-1}}\iota_{1,n-1}^* \NB_{n-1} \rightarrow \NB_{\geq n} \bar 1_n \iota_{1,n-1}^*, u\otimes v \mapsto u \iota_{1,n-1}(v) $
for any $u \in \NB_{\geq n} \bar 1_{n-1}, v \in \NB_{n-1}$. 
This is obviously a graded bimodule homomorphism.
For $b$, we show equivalently that the
map $\NB_{\geq(n+1)}\bar 1_{n+1} \otimes_{\NB_{n+1}}
\NB_{n+1}\iota_{1,n}^* \rightarrow \coker a$
that is the inverse of the linear map induced by $b$
is a graded bimodule homomorphism.
This inverse map is defined explicitly by
\begin{align*}
\NB_{\geq(n+1)}\bar 1_{n+1} \otimes_{\NB_{n+1}}
\NB_{n+1}\iota_{1,n}^*&\rightarrow 
1_{B\star-}\NB_{\geq n} \bar 1_{n} / \im a,\\
\begin{tikzpicture}[anchorbase,scale=1.25]
\draw[-] (0.4,0.7)--(0.4,1.15);
\draw[-] (-0.35,0.7)--(-0.35,1.15);
\draw[-] (0.4,0.2)--(0.4,.5);
\draw[-] (-0.35,0.2)--(-0.35,.5);
\node at (-0.19,0.32) {$\cdot$};
\node at (0.03,0.32) {$\cdot$};
\node at (0.25,0.32) {$\cdot$};
\node at (-0.19,1) {$\cdot$};
\node at (0.03,1) {$\cdot$};
\node at (0.25,1) {$\cdot$};
\node[rectangle,rounded corners,draw,fill=blue!15!white,inner sep=4.2pt] at (0.02,0.67) {$\hspace{3.5mm}\scriptstyle u\hspace{3.5mm}$};
\end{tikzpicture}
\:\:\otimes\:\:
\begin{tikzpicture}[anchorbase,scale=1.25]
\draw[-] (0.4,0.7)--(0.4,1.15);
\draw[-] (-0.35,0.7)--(-0.35,1.15);
\draw[-] (0.4,0.2)--(0.4,.5);
\draw[-] (-0.35,0.2)--(-0.35,.5);
\node at (-0.09,0.32) {$\cdot$};
\node at (0.08,0.32) {$\cdot$};
\node at (0.25,0.32) {$\cdot$};
\node at (-0.19,1) {$\cdot$};
\node at (0.03,1) {$\cdot$};
\node at (0.25,1) {$\cdot$};
\draw[-] (-.2,.2)--(-.2,.5);
\node[rectangle,rounded corners,draw,fill=blue!15!white,inner sep=4.2pt] at (0.02,0.67) {$\hspace{3.5mm}\scriptstyle v\hspace{3.5mm}$};
\end{tikzpicture}\:\:
&\mapsto\:\:
\begin{tikzpicture}[anchorbase,scale=1.25]
\draw[-] (-0.5,0.5)--(-0.5,-0.65)
to [out=-90,in=180] (-.35,-.9)
to [out=0,in=-90] (-.2,-.65);
\draw[-] (-0.2,0)--(-0.2,-0.5);
\draw[-] (0.4,0)--(0.4,-0.5);
\draw[-] (-0.2,0.2)--(-0.2,0.5);
\draw[-] (0.4,0.2)--(0.4,0.5);
\draw[-] (0.4,-.65)--(0.4,-.95);
\draw[-] (-.1,-.65)--(-.1,-.95);
\node at (-0.05,0.4) {$\cdot$};
\node at (0.1,0.4) {$\cdot$};
\node at (0.25,0.4) {$\cdot$};
\node at (-.05,-0.2) {$\cdot$};
\node at (0.1,-0.2) {$\cdot$};
\node at (0.25,-0.2) {$\cdot$};
\node at (0.02,-.82) {$\cdot$};
\node at (0.155,-.82) {$\cdot$};
\node at (0.29,-.82) {$\cdot$};
\node[rectangle,rounded corners,draw,fill=blue!15!white,inner sep=4.2pt] at (0.11,0.1) {$\hspace{3.5mm}\scriptstyle u\hspace{3.5mm}$};
\node[rectangle,rounded corners,draw,fill=blue!15!white,inner sep=4.2pt] at (0.11,-0.5) {$\hspace{3.5mm}\scriptstyle v\hspace{3.5mm}$};
\end{tikzpicture}
\:\:+\: \im a
\end{align*}
for any $u \in \NB_{\geq (n+1)} \bar 1_{n+1}, v \in \NB_{n+1}$, which is a graded bimodule homomorphism.

It remains to check \cref{space}. 
Take $m \geq 0$.
By its definition,
$a_m:
\bar 1_m \NB_{\geq (n-1)} \bar 1_{n-1} \otimes_{\NB_{n-1}} 
\iota^*_{1,n-1} \NB_{n}
\rightarrow
\bar 1_{m+1} \NB_{\geq n} \bar 1_n$
intertwines left multiplication by $1 \otimes x_1$ with left
multiplication by $\txtdot \star 1_m$.
This implies the statement about $\alpha$,
noting that a sign appears since
$x:B \Rightarrow B$ corresponds to 
$-\Res_{\!\smalltxtdot\star-}$ 
in \cref{snow}.
Similarly, for $\beta$, one checks from the definition that
$b_m:
\bar 1_{m+1} \NB_{\geq n} \bar 1_n
\rightarrow
\bar 1_m \NB_{\geq (n+1)} \bar 1_{n+1}
\otimes_{\NB_{n+1}} \NB_{n+1} \iota^*_{1,n}
$ intertwines left multiplication by 
$\txtdot\star 1_m$
with right multiplication by
$1 \otimes x_1$.
\end{proof}

\cref{Thm:SES} implies that
the functor $B$ 
preserves modules with $\Delta$-flags and with $\bar\Delta$-flags.
The next two theorems makes this more precise. The combinatorics that emerges here matches \cref{stdrecurrence,witcher2}.

\begin{theorem}\label{burger}
Consider
the short exact sequence
$$
0 \longrightarrow 
K(n) \longrightarrow B \Delta(n)
\longrightarrow Q(n)
\longrightarrow 0
$$
obtained
by applying \cref{Thm:SES} to the $\NB_n$-module $P_n(n)\:(n \geq 0)$.
We denote the endomorphisms 
$\jmath^{n-1}_! x_{\Delta(n)}' :K(n)\rightarrow K(n)$
and $\jmath^{n+1}_! x_{\Delta(n)}'' :Q(n)\rightarrow Q(n)$ from \cref{space}
by $y$ and $z$, respectively.
\begin{enumerate}
\item
Assuming that $n > 0$ so that $K(n) \neq 0$,
we have that $K(n) \cong 
\Delta(n-1)^{\oplus q^{1-n} / (1-q^{2})}$.
More precisely, we have that
$$
K(n) \cong 
q^{1-n}\Gamma[y]  \otimes_\Gamma
\Delta(n-1)$$ 
with the action of $\NB$ being on the second tensor factor. This isomorphism
may be chosen so that
the endomorphism $y$ of $K(n)$
corresponds to multiplication by $y$ on the 
first tensor factor.
\item
We have that $Q(n) \cong \Delta(n+1)^{\oplus [n+1]}$.
More precisely, recalling also \cref{twocents},
$$
Q(n) \cong q^{-n} \ZZ_{n+1}[z] / \big((z-x_1)\cdots (z-x_{n+1})\big) \otimes_{\ZZ_{n+1}}
\Delta(n+1)
$$
with the action of $\NB$ being on the second tensor factor. 
This isomorphism may be chosen
so that the endomorphism $z$ of $Q(n)$
 corresponds to multiplication by $z$ on the first tensor factor.
\end{enumerate}
\end{theorem}

\begin{proof}
(1) 
According to \cref{Thm:SES}, we 
have that $K(n) = \jmath^{n-1}_!(R_{1,n-1} P_n(n))$,
and the endomorphism $y$ of $K(n)$
is obtained by applying the functor $\jmath^{n-1}_!$
to the endomorphism 
we also denote $y:= x'_{P_n(n)}$ of
$R_{1,n-1} P_n(n)$
defined by left multiplication by $-x_1$.
Therefore, by exactness of $\jmath^{n-1}_!$, 
it suffices to prove 
that
$R_{1,n-1} P_n(n) \cong q^{1-n} \Gamma[y] \otimes_\Gamma P_{n-1}(n-1)$
as a graded $\NB_1 \otimes_\kk \NB_{n-1}$-module,
identifying $\NB_1$ with $\Gamma[y]$ so $y=-x_1$.
This follows because
$$
P_n(n) = q^{-\frac{1}{2}n(n-1)}
\Gamma[x_1,x_2,\dots,x_n]
\cong
q^{1-n} \Gamma[y] \otimes_\Gamma q^{-\frac{1}{2}(n-1)(n-2)}
\Gamma[x_2,\dots,x_n].
$$

\vspace{2mm}
\noindent
(2)
By \cref{Thm:SES}, we have that
$Q(n) = \jmath^{n+1}_! (I_{1,n} P_n(n))$,
and the endomorphism $z$ of $Q(n)$ is 
obtained by applying $\jmath^{n+1}_!$ to the endomorphism also denoted $z := x''_{P_n(n)}$
of $I_{1,n} P_n(n)$ defined by right multiplication by $x_1$.
Therefore, it suffices to show that
$$
I_{1,n} P_n(n) \cong q^{-n} \ZZ_{n+1}[z] / \big((z-x_1)\cdots (z-x_{n+1})\big) \otimes_{\ZZ_{n+1}} P_{n+1}(n+1)
$$
as a graded $\NB_{n+1}$-module, where the action 
is on the second tensor factor.
Using \cref{bundy}, it is easy to check that both sides have the same graded dimensions. Hence, it suffices to construct a degree-preserving surjective homomorphism
\begin{equation}\label{want}
\bar \theta:q^{-n} \ZZ_{n+1}[z] / \big((z-x_1)\cdots (z-x_{n+1})\big) \otimes_{\ZZ_{n+1}} P_{n+1}(n+1)
\twoheadrightarrow
\NB_{n+1} \iota_{1,n}^* \otimes_{\NB_n} P_n(n).
\end{equation}
Recall that $P_{n+1}(n+1)$
is generated by $u_{n+1}$ subject to the relations
$\tau_i u_{n+1} = 0$ for $i=1,\dots,n$.
It is easy to see that $\tau_n \cdots \tau_2 \tau_1 
x_1^r 
\otimes u_n$ is annihilated by all $\tau_i$.
Hence, there is a unique graded $\NB_{n+1}$-module homomorphism such that
\begin{align*}
\theta:q^{-n} \ZZ_{n+1}[z] \otimes_{\ZZ_{n+1}} P_{n+1}(n+1)
&\rightarrow \NB_{n+1} \iota_{1,n}^* \otimes_{\NB_n}
P_n(n),
&
z^r \otimes u_{n+1}
&\mapsto \tau_n \cdots \tau_2 \tau_1 x_1^r \otimes
u_{n}
\end{align*}
for any $r \geq 0$.
This takes
$(z-x_1) \cdots (z-x_{n+1}) \otimes u_{n+1}$
to $\tau_n\cdots \tau_2 \tau_1 (x_1-x_1) \cdots (x_1-x_{n+1}) \otimes u_n = 0$.
Hence, we get induced a 
graded $\NB_{n+1}$-module homomorphism
$\bar\theta$ as in \cref{want}.
It remains to show that this is surjective.
The module on the right hand side is cyclic with generator
$1 \otimes u_n$, so we just need to see that it is in the image of $\bar\theta$.
To see this, we show by induction on $m=0,1,\dots,n$
that $1 \otimes u_n$
lies in the submodule generated by
$\tau_m \cdots \tau_2 \tau_1 x_1^r \otimes u_n\:(0 \leq r \leq m)$; the $m=n$ case of this gives what we need.
The base case $m=0$ of the induction is trivial. The induction step follows from the relation
\begin{equation}\label{sierpinski}
\tau_m \cdots \tau_2 \tau_1 x_1^{m} \otimes u_n
=x_{m+1} \tau_m \cdots \tau_2 \tau_1
x_1^{m-1}\otimes u_n
+ \tau_{m-1}\cdots \tau_2 \tau_1 x_1^{m-1}
\otimes u_n,
\end{equation}
which follows using \cref{lastrel}.
\end{proof}

\begin{theorem}\label{burger2}
Consider
the short exact sequence
$$
0 \longrightarrow 
\bar K(n) \longrightarrow B \bar\Delta(n)
\longrightarrow \bar Q(n)
\longrightarrow 0
$$
obtained
by applying \cref{Thm:SES} to the $\NB_n$-module $L_n(n)\:(n \geq 0)$.
We denote the endomorphisms 
$\jmath^{n-1}_! x_{\bar\Delta(n)}' :\bar K(n)\rightarrow \bar K(n)$
and $\jmath^{n+1}_! x_{\bar\Delta(n)}'' :\bar Q(n)\rightarrow \bar Q(n)$ from \cref{space}
by $\bar y$ and $\bar z$, respectively.
\begin{enumerate}
\item
Assuming that $n > 0$ so that $\bar K(n)$ is non-zero, the module $\bar K(n)$ is a $\bar\Delta$-layer
that is equal in the Grothendieck group to
$[n]\left[\bar\Delta(n-1)\right]$. More precisely, letting $\bar K_i(n)$
be the image of $\bar y^i:\bar K(n)\rightarrow \bar K(n)$ defines
a graded filtration
$$
\bar K(n)
= \bar K_0(n) > \bar K_1(n)>\cdots > \bar K_{n}(n) = 0
$$ 
such that $\bar K_{i-1}(n) / \bar K_i(n)
\cong q^{2i-n-1} \bar\Delta(n-1)$
for $i =1,\dots,n$.
Also
\begin{equation}\label{hd1}
\dim_q \Hom_{\NB}(\bar K(n), \bar L(n-1)) = q^{n-1}.
\end{equation}
\item
The module $\bar Q(n)$ is a $\bar\Delta$-layer
equal in the Grothendieck group to $q^{-n}\left[\bar\Delta(n+1)\right] / (1-q^{2})$.
More precisely, 
letting $\bar Q_i(n)$
be the image of $\bar z^i:\bar Q(n)\rightarrow \bar Q(n)$ defines a 
graded filtration
$$
\bar Q(n) = \bar Q_0(n) > \bar Q_1(n) > \bar Q_2(n)>\cdots
$$
such that $\bar Q_{i-1}(n) / \bar Q_i(n)
\cong q^{2i-n-2} \bar\Delta(n+1)$ for $i \geq 1$.
Also
\begin{equation}\label{hd2}
\dim_q \Hom_{\NB}(\bar Q(n), \bar L(n+1)) = q^{n}.
\end{equation}
\end{enumerate}
\end{theorem}

\begin{proof}
(1)
Let $V := R_{1,n-1}L_n(n)$ and 
$\bar y:V \rightarrow V$ be the endomorphism defined by multiplication by $-x_1$.
Let $V_i := \im \bar y^i$.
Like in the proof of the previous theorem, the proof of the first assertion in (1)
reduces to showing that
$V_{i-1} / V_{i} \cong q^{n+1-2i} L_{n-1}(n-1)$
as a graded $\NB_{n-1}$-module
for $i=1,\dots,n$, and that $V_n = 0$.
We have that 
$$
\sum_{r=0}^n (-1)^r x_1^{n-r} e_{r,n}
=(x_1-x_1) (x_1-x_2) \cdots (x_1-x_n) = 0,
$$
where
$e_{r,n}$ is the $r$th elementary symmetric polynomial in $x_1,\dots,x_n$. 
Also let $e_{r,n}'$ be the $r$th elementary
symmetric polynomial in $x_2,\dots,x_n$.
Since $e_{r,n}$ acts as 0 on $L_n(n)$
for $r \geq 1$,
it follows that $x_1^n$ acts as 0 too.
This implies that $V_n = 0$.
Now take $1 \leq i \leq n$. 
We claim that there is a graded $\NB_{n-1}$-module
homomorphism
\begin{align*}
\theta_i: q^{2i-n-1} 
L_{n-1}(n-1)
&\rightarrow V_{i-1} / V_{i},&
\bar u_{n-1} &\mapsto 
x_1^{i-1} \bar u_n + V_{i}.
\end{align*}
This follows 
using the generators and relations for
$L_{n-1}(n-1)$ discussed earlier
since $\tau_2,\dots,\tau_{n-1}$
annihilate $x_1^{i-1} \bar u_n$,
as does $e_{r,n}'$
for each $r \geq 1$.
To see the latter assertion,
we have that \begin{equation}\label{resist}
e_r'
= e_{r,n}
- x_1 e'_{r-1}.
\end{equation}
The first term on the right-hand side of \cref{resist} is 0 on $x_1^{i-1} \bar u_n$,
and the second term maps it to $V_{i}$.
This proves the claim.
Finally, each $\theta_i$ is actually an isomorphism. This follows by considering
the explicit bases for $L_n(n)$
and $L_{n-1}(n-1)$ from \cref{LBasis}.

It remains to prove \cref{hd1}.
We have that
\begin{align*}
\Hom_{\NB}(\bar K(n),
L(n-1)) &=
\Hom_{\NB_{\geq (n-1)}}
\left(\jmath^{n-1}_! (R_{1,n} L_n(n)),
L(n-1)\right)\\
&\cong
\Hom_{\NB_{n-1}}
\left(R_{1,n} L_n(n),
\jmath^{n-1} L(n-1)\right)\\&\cong
\Hom_{\NB_{n-1}}(
R_{1,n} L_n(n),
L_{n-1}(n-1)).
\end{align*}
Let $f:R_{1,n} L_n(n)\rightarrow L_{n-1}(n-1)$
be an $\NB_{n-1}$-module homomorphism.
Since $x_1 = (x_1+\cdots+x_n) - (x_2+\cdots+x_n)$ and $x_1+\cdots+x_n$
annihilates $L_n(n)$
as it is central of positive degree,
we see that $$
f(x_1^i \bar u_n) 
= (-1)^i f\left((x_2+\cdots+x_n)^i\bar u_n\right)
=(-1)^i (x_1+\cdots+x_{n-1})^i f(\bar u_n).
$$
This is 0 for $i \geq 1$.
It follows that $f$ sends the submodule $V_1$
defined in the previous paragraph to 0.
Thus, it factors through the quotient
$V_0 / V_1 \cong q^{1-n} L_{n-1}(n-1)$.
Using Schur's Lemma, we deduce that
\begin{multline}
\qquad\dim_q 
\Hom_{\NB_{n-1}}(R_{1,n} L_n(n), L_{n-1}(n-1))
=\\
\dim_q 
\Hom_{\NB_{n-1}}\left(q^{n-1} L_{n-1}(n-1), L_{n-1}(n-1)\right)
= q^{n-1}.\qquad\label{store}
\end{multline}

\noindent
(2)
Let $W := I_{1,n} L_n(n) = \NB_{n+1}\iota_{1,n}^* \otimes_{\NB_n} L_n(n)$ and $\bar z:W \rightarrow W$
be the endomorphism defined by 
right multiplying the bimodule $\NB_{n+1} \iota_{1,n}^*$
by $x_1$.
Let $W_i := \im \bar z^i$.
For the first assertion, 
we need to show that $W_{i-1} / W_i
\cong q^{2i-n-2} L_{n+1}(n+1)$ for each $i \geq 1$.
The argument using \cref{sierpinski}
explained at the end of the proof of \cref{burger}
shows that $W$ is generated as an
$\NB_{n+1}$-module
by the vectors $\tau_n \cdots \tau_2 \tau_1 x_1^j\otimes \bar u_n$
for all $j\geq 0$ (actually, one just needs them for
$0 \leq j \leq n$).
It follows that $W_i$ is generated
by the vectors
 $\tau_n \cdots \tau_2 \tau_1 x_1^j \otimes \bar u_n$
for all $j\geq i$, and 
$W_{i-1} / W_i$ is a cyclic $\NB_{n+1}$-module 
generated
by 
$\tau_n \cdots \tau_2 \tau_1 x_1^{i-1} 
\otimes \bar u_n + W_i$.
For any $i \geq 1$, we claim that there is a surjective
graded $\NB_{n+1}$-module homomorphism
\begin{align*}
\theta_i: q^{2i-n-2} L_{n+1}(n+1)
&\twoheadrightarrow W_{i-1} / W_i,
&
\bar u_{n+1} &\mapsto
\tau_n\cdots\tau_2\tau_1 x_1^{i-1}
\otimes \bar u_n + W_i.
\end{align*}
To see this, it just remains to check the relations:
each of $\tau_1,\dots,\tau_n$ annihilates
$\tau_n\cdots\tau_2\tau_1 x_1^{i-1}
\otimes \bar u_n + W_i$ by some easy commutation relations
using \cref{relA,relB,relC},
and $e_{r,n+1}$ does
too for $r \geq 1$, 
as may be deduced using \cref{resist}.
Finally, one checks graded dimensions
using \cref{Ldim,exact} to see that each $\theta_i$ must actually be an isomorphism.

Now consider \cref{hd2}. This reduces like before to showing that 
$\dim_q \Hom_{\NB_{n+1}}(I_{1,n} L_n(n),
L_{n+1}(n+1)) = q^{n}.
$
For this, we note using adjointness 
and duality 
that
\begin{align*}
\Hom_{\NB_{n+1}}(I_{1,n} L_n(n),
L_{n+1}(n+1)) &\cong
\Hom_{\NB_{n+1}}(L_n(n),
R_{1,n} L_{n+1}(n+1))\\&\cong
\Hom_{\NB_{n+1}}(R_{1,n} L_{n+1}(n+1),
L_{n}(n)).
\end{align*}
This is of graded dimension $q^{n}$ by \cref{store}.
\end{proof}

\subsection{Character formulae}\label{grch}
The {\em graded character} of a
locally finite-dimensional graded left $\NB$-module $V$ is defined by
\begin{equation}
\ch V := \sum_{n \geq 0}
\dim_q(1_n V) \chi^n.
\end{equation}
In general, this is a power series
in the formal variable $\chi$ with coefficients that are themselves
formal series of the form $\sum_{n \in \Z} a_n q^n$
for $a_n \in \N$. 
The graded character of any finitely generated graded module (or any
module that is locally finite dimensional and bounded below) belongs to
$\Z\lround q \rround \llbracket \chi \rrbracket$.
This is an integral form for the
completion $\Q\lround q \rround \llbracket \chi\rrbracket$ of the character ring from \cref{charring}.

We obviously have that
\begin{equation}
\ch (V^\circledast) = \psi^\imath(\ch V)
\end{equation}
where $\psi^\imath$ on the right-hand side is the bar involution on the character ring from \cref{barinv}. Also
\begin{equation}\label{reminder0}
\ch (B V) = b (\ch V)
\end{equation}
where the action of $b$ on $\Z\lround q\rround\llbracket \chi \rrbracket$ on the right-hand side is defined as in \cref{dualaction}.
This identity is easy to see if one views $B$
as the functor $\Res_{|\star-}$ as explained in \cref{snow}.

The irreducible module $L(n)$
has (globally) finite-dimensional weight
spaces by general theory, so its graded character actually lies in
$\Z[q,q^{-1}] \llbracket \chi \rrbracket$,
as does the formal character of any 
graded module of finite length.
By lowest weight theory, we clearly have that
\begin{equation}
\ch L(n) \equiv [n]! \chi^n \pmod{\chi^{n+1} \Z[q,q^{-1}]\llbracket \chi\rrbracket},
\end{equation}
which implies that the irreducible characters
are linearly independent. They are also invariant under $\psi^\imath$ since $L(n)$
is self-dual.
Now recall the following expressions 
defined/computed in
\cref{reaching,oops}: 
\begin{align}\label{reminders1}
\bar\delta_n &=
[n]! 
\sum_{f \geq 0} \frac{T_{f,n}(q^{-2}) }{(1-q^{2})^f} \chi^{n+2f},\\
\el_n &=
[n]! \sum_{m \geq 0}
\left(\sum_{\alpha \in \Par_t(m\times n)}
[\alpha_1+1]^2\cdots [\alpha_m+1]^2\right) \chi^{n+2m}.\label{reminders2}
\end{align}
These are the graded characters
of proper standard and irreducible modules:

\begin{theorem}
\label{daughters}
For any $n \in \N$, we have that 
$\ch \bar\Delta(n) = \bar\delta_n$
and $\ch L(n) = \el_n$.
\end{theorem}

\begin{proof}
The equality $\ch\bar\Delta(n)=\bar\delta_n$ follows on computing the graded character of $\bar\Delta(n)$ by counting
vectors of each degree
in the basis \cref{properstandardbasis},
using also the combinatorics discussed
in \cref{reader}.
To prove that $\ch L(n) = \el_n$,
\cref{lastgasp} implies that 
\begin{align*}
\dim_q 1_{n} L(n-2m) &= 
\dim_q \Hom_{\NB}(\NB 1_{n}, L(n-2m))\\
&=
[n-2m]! 
\sum_{\alpha \in \Par_t(m\times (n-2m))}
[\alpha_1+1]^2\cdots [\alpha_m+1]^2.
\end{align*}
Replacing $n$ by $n+2m$ throughout, 
this shows that the $\chi^{n+2m}$-coefficient
of $\ch L(n)$ is the same as this coefficient in the formula
\cref{reminders2} for $\el_n$.
\end{proof}

Using also the identity \cref{oopsformula},
\cref{daughters} proves \cref{introthm:char} from the introduction, and \cref{introthm:decomposition} follows from \cref{frombefore}.

\subsection{Branching rules}
We end by describing the effect of the projective functor $B$ on the 
irreducible module $L(n)$.
In view of \cref{daughters,reminder0}, we can reinterpret \cref{witcher1} as
\begin{equation}\label{sundaymorning}
\ch \left(B L(n)\right) = [n] \ch L(n-1)
+ \delta_{n\not\equiv t} [n+1] \ch L(n+1).
\end{equation}
Since the irreducible characters are linearly independent, this provides
complete information about the composition factors of $B L(n)$.
In particular, we see that
\begin{equation}\label{somethingeasy}
B L(0) 
\cong 
\begin{dcases}
L(1)&\text{if $t=1$}\\
0&\text{if $t=0$}.
\end{dcases}
\end{equation}
 
\begin{lemma}\label{branch1}
Interpreting $L(-1)$ as $0$, 
the following hold for all $n \geq 0$:
\begin{enumerate}
\item
$\head B \bar\Delta(n)
\cong 
\begin{dcases}
q^{-n} L(n+1) \oplus q^{1-n} L(n-1)\hspace{4mm}&\text{if $n \equiv t \pmod{2}$}\\
q^{-n} L(n+1)&\text{if $n \not\equiv t \pmod{2}$.}
\end{dcases}$
\item
$\soc B \bar\nabla(n)
\cong \begin{dcases}
q^{n} L(n+1) \oplus q^{n-1} L(n-1)\hspace{5mm}&\text{if $n \equiv t \pmod{2}$}\\
q^{n} L(n+1)&\text{if $n \not\equiv t \pmod{2}$.}
\end{dcases}$
\item
$\head B L(n) \cong
\begin{dcases}
q^{1-n} L(n-1)\hspace{28.7mm}&\text{if $n \equiv t \pmod{2}$}\\
q^{-n} L(n+1)&\text{if $n \not\equiv t \pmod{2}$.}
\end{dcases}$
\item
$\soc B L(n) \cong
\begin{dcases}
q^{n-1} L(n-1)\hspace{27.4mm}&\text{if $n \equiv t \pmod{2}$}\\
q^{n} L(n+1)&\text{if $n \not\equiv t \pmod{2}$.}
\end{dcases}$
\end{enumerate}
\end{lemma}

\begin{proof}
We first treat the case $n=0$.
Parts (3) and (4) are immediate from \cref{somethingeasy}.
For (1), \cref{burger2}(2)
shows that  $B \bar\Delta(0) \cong \bar Q(0)$, and this module
has irreducible head $L(1)$.
Then (2) follows (1) by duality.
Assume for the rest of the proof that $n \geq 1$.

By duality, (1) and (2) are equivalent, as are (3) and (4).
By \cref{burger2}, especially
\cref{hd1,hd2}, it is clear
that $\head B \bar\Delta(n)$
is isomorphic {\em either} to $q^{-n} L(n+1) \oplus q^{1-n} L(n-1)$ {\em or} to $q^{-n} L(n+1)$.
The following claim completes the proof of (1) and (2) when $n\not\equiv t\pmod{2}$.

\vspace{2mm}
\noindent
{\bf Claim.} {\em
If $n \not\equiv t\pmod{2}$
then
$\Hom_{\NB}(B\bar\Delta(n), L(n-1)) = 0$.}

\vspace{1mm}
\noindent
To prove this, 
we let $V := \Res_{|\star-} \bar\Delta(n)$,
this being isomorphic to $B \bar \Delta(n)$
by \cref{snow}.
In this incarnation, 
the submodule $\bar K(n)$ from \cref{burger2}(1) is identified with the submodule $K$ of $V$
generated by the vectors
$x_1^{i-1} \bar v_n$ for $1 \leq i \leq n$.
This is apparent from the proofs of \cref{Thm:SES} and 
\cref{burger2}(1).
Any non-zero homomorphism $f:K
\rightarrow L(n-1)$ resulting 
from \cref{hd1}
is necessarily homogeneous of degree
$n-1$, and must 
take $\bar v_n$ to a {\em non-zero} vector of the minimal degree $-\frac{1}{2}(n-1)(n-2)$ 
in $1_{n-1}L(n-1)$.
We are trying to show that
$f$ does not extend to a homogeneous homomorphism
$\hat f: V\rightarrow L(n-1)$.
Suppose for a contradiction that there is such an extension.
Consider the vectors
\begin{align*}
v &:=\begin{tikzpicture}[anchorbase,scale=1.4]
\draw[-] (-0.4,0.15)--(-0.4,0.8);
\draw[-] (0.4,0.15)--(0.4,0.8);
\draw[-] (0.2,0.15)--(0.2,0.8);
\draw[-] (-.65,.8) to[out=-90,in=180] (0,.5) to[out=0,in=-90] (.65,.8);
\node at (-0.25,0.35) {$\cdot$};
\node at (-.1,0.35) {$\cdot$};
\node at (0.05,0.35) {$\cdot$};
\node at (-0.25,0.65) {$\cdot$};
\node at (-.1,0.65) {$\cdot$};
\node at (0.05,0.65) {$\cdot$};
\node[rectangle,rounded corners,draw,fill=blue!15!white,inner sep=4pt] at (0,0) {$\hspace{5mm}\scriptstyle \bar v_n\hspace{5mm}$};
\end{tikzpicture}
\:\:&
w &:=\begin{tikzpicture}[anchorbase,scale=1.4]
\draw[-] (-.6,.9) to (-.6,.75);
\draw[-] (-0.4,0.15)--(-0.4,0.9);
\draw[-] (0.2,0.15)--(0.2,0.9);
\draw[-] (0.4,0.15)--(0.4,0.75);
\draw[-] (-.6,.8) to[out=-90,in=180] (0,.5) to[out=0,in=-90] (.6,.75) to [out=90,in=0] (.5,.85) to [out=180,in=90] (.4,.75);
\node at (-0.25,0.35) {$\cdot$};
\node at (-.1,0.35) {$\cdot$};
\node at (0.05,0.35) {$\cdot$};
\node at (-0.25,0.7) {$\cdot$};
\node at (-.1,0.7) {$\cdot$};
\node at (0.05,0.7) {$\cdot$};
\closeddot{.6,.75};\node at (.75,.75) {$\scriptstyle n$};
\node[rectangle,rounded corners,draw,fill=blue!15!white,inner sep=4pt] at (0,0) {$\hspace{5mm}\scriptstyle\bar v_n\hspace{5mm}$};
\end{tikzpicture}
\end{align*}
The vector $v$ 
is of degree $-\frac{1}{2}n(n-1)-2n$,
so 
$\hat f(v)$ is of degree $-\frac{1}{2}(n-1)(n-2)-2n$,
which is smaller than the degree of any non-zero vector in 
$1_{n+1} \bar\Delta(n-1)$, hence, in $1_{n+1} L(n-1)$.
So $\hat f(v) = 0$.
Since $w$ is obtained from $v$ by acting 
with some element of $\NB$, we deduce
that $\hat f(w) = 0$ too.
Now we calculate using \cref{curlstuff,rels10}
and the defining relations of $L_n(n)$
to see that 
$$
w = \begin{tikzpicture}[anchorbase,scale=1.4]
\draw[-] (-.6,.9) to (-.6,.75);
\draw[-] (-0.4,0.15)--(-0.4,0.9);
\draw[-] (0.2,0.15)--(0.2,0.9);
\draw[-] (0.4,0.15)--(0.4,0.75);
\draw[-] (-.6,.8) to[out=-90,in=180] (0,.5) to[out=0,in=-90] (.6,.75) to [out=90,in=0] (.5,.85) to [out=180,in=90] (.4,.75);
\node at (-0.25,0.35) {$\cdot$};
\node at (-.1,0.35) {$\cdot$};
\node at (0.05,0.35) {$\cdot$};
\node at (-0.25,0.7) {$\cdot$};
\node at (-.1,0.7) {$\cdot$};
\node at (0.05,0.7) {$\cdot$};
\closeddot{.6,.75};\node at (.75,.75) {$\scriptstyle n$};
\node[rectangle,rounded corners,draw,fill=blue!15!white,inner sep=4pt] at (0,0) {$\hspace{5mm}\scriptstyle\bar v_n\hspace{5mm}$};
\end{tikzpicture}
=- \begin{tikzpicture}[anchorbase,scale=1.4]
\draw[-] (-.6,.9) to (-.6,.75);
\draw[-] (-0.4,0.15)--(-0.4,0.9);
\draw[-] (0.2,0.15)--(0.2,0.9);
\draw[-] (-.6,.8) to[out=-90,in=170] (-0.1,.5) to[out=-10,in=90] (.4,.15);
\node at (-0.25,0.35) {$\cdot$};
\node at (-.1,0.35) {$\cdot$};
\node at (0.05,0.35) {$\cdot$};
\node at (-0.25,0.7) {$\cdot$};
\node at (-.1,0.7) {$\cdot$};
\node at (0.05,0.7) {$\cdot$};
\closeddot{.34,.32};\node at (.63,.38) {$\scriptstyle n-1$};
\node[rectangle,rounded corners,draw,fill=blue!15!white,inner sep=4pt] at (0,0) {$\hspace{5mm}\scriptstyle\bar v_n\hspace{5mm}$};
\end{tikzpicture}
=
(-1)^n \begin{tikzpicture}[anchorbase,scale=1.4]
\draw[-] (-0.4,0.15)--(-0.4,0.9);
\draw[-] (-0.2,0.15)--(-0.2,0.9);
\draw[-] (0.4,0.15)--(0.4,0.9);
\node at (-0.1,0.5) {$\cdot$};
\node at (0.05,0.5) {$\cdot$};
\node at (0.2,0.5) {$\cdot$};
\node[rectangle,rounded corners,draw,fill=blue!15!white,inner sep=4pt] at (0,0) {$\hspace{5mm}\scriptstyle\bar v_n\hspace{5mm}$};
\end{tikzpicture}
= (-1)^n \bar v_n.
$$
The first equality here
requires $n \not\equiv t \pmod{2}$---otherwise, it would be 0.
Now we have that $\hat f(w) = (-1)^n \hat f(\bar v_n) = 0$ but
$\hat f(\bar v_n) \neq 0$.
This contradiction proves the claim.

Next, consider $\head B L(n)$.
For $m \geq 0$, 
$\Hom_{\NB}(B L(n), L(m))$
embeds naturally into both of the spaces
$\Hom_{\NB}(B \bar\Delta(n), L(m))$
and $\Hom_{\NB}(B L(n), \bar\nabla(m))
\cong \Hom_{\NB}(L(n), B \bar\nabla(m))$.
So the parts of (1)--(2) proved so far
imply:
\begin{itemize}
\item $\dim_q \Hom_{\NB}(B L(n), L(m))=0$ 
if $m \neq n\pm 1$.
\item 
$\dim_q \Hom_{\NB}(B L(n), L(n+1)) = 0$ or $q^{n}$.
\item 
$\dim_q \Hom_{\NB}(B L(n), L(n-1)) = 0$ or $q^{n-1}$.
\end{itemize}
If $n \not\equiv t\pmod{2}$
then $\Hom_{\NB}(B L(n), L(n-1)) = 0$
as $\Hom_{\NB}(B \bar\Delta(n), L(n-1))=0$.
Since $B L(n) \neq 0$ by \cref{sundaymorning},
we must therefore have that
$\Hom_{\NB}(B L(n), L(n+1)) \neq 0$, so 
its graded dimension is $q^{n}$.
Hence,
$\head B L(n) \cong q^{-n} L(n+1)$ in this situation.
Instead, if $n \equiv t \pmod{2}$
then we have that $\Hom_{\NB}(B L(n), L(n+1))=0$
as $\Hom_{\NB}(L(n), B \bar\nabla(n+1)) = 0$.
Since $B L(n) \neq 0$, 
we must therefore have that $\Hom_{\NB}(B L(n), L(n-1)) \neq 0$.
So it has graded dimension $q^{n-1}$,
and we have proved that $\head B L(n) \cong q^{1-n} L(n-1)$.
Now (3) and (4) are proved.

Finally, we complete the proof of (1) and (2) in the remaining case that $n\equiv t\pmod{2}$. We need to show that 
$\Hom_{\NB}(B\bar\Delta(n),L(n-1))$
and $\Hom_{\NB}(L(n-1), B\bar\nabla(n))$
are non-zero. This follows because
$\Hom_{\NB}(B L(n), L(n-1))$
and
$\Hom_{\NB}(L(n-1), BL(n))$
are non-zero by (3)--(4).
\end{proof}

\begin{theorem}\label{lowerbound}
For $n \geq 0$, the module $V := B L(n)$ is
uniserial.
To describe its unique composition series, let
$x:V \rightarrow V$ denote the nilpotent endomorphism $x_{L(n)}$,
$V_i := \im x^i$ and $V^i := \ker x^i$.
\begin{enumerate}
\item If $n \equiv t \pmod{2}$ then the unique composition series is
$$
V = V_0 = V^n > V_1 = V^{n-1} > V_2=V^{n-2} > \cdots >  V^1 > V_n = V^0 = 0
$$
with  $V_{i-1} / V_i = V^{n+1-i} / V^{n-i} \cong q^{2i-n-1} L(n-1)$
for each $i=1,\dots,n$.
\item
If $n \not \equiv t \pmod{2}$ then the unique composition series is
$$
V = V_0 > V^n > V_1 > V^{n-1} > V_2 > V^{n-2} >\cdots > V^1 > V_n >
V^0 = 0
$$
with 
$V_{i-1} / V^{n+1-i} \cong q^{2i-n-2}L(n+1)$ for $i=1,\dots,n+1$ and
$V^{n+1-i} / V_{i} \cong q^{2i-n-1} L(n-1)$ for $i=1,\dots,n$.
\end{enumerate}
Moreover,
$\End_{\NB}(V) = \kk[x] / \left(x^{\beta(n)}\right)$
with $\beta(n) = n$ if $n\equiv t\pmod{2}$
or $n+1$ if $n \not\equiv t\pmod{2}$.
\end{theorem}

\begin{proof}
Since
$V$ is a quotient of 
$B \bar\Delta(n)$, \cref{burger2} implies that
there is a short exact sequence
$$
0 \longrightarrow K \longrightarrow V
\longrightarrow Q \longrightarrow 0
$$
where $K$ is a quotient of $\bar K(n)$
and $Q$ is a quotient of $\bar Q(n)$.
The filtrations of $\bar K(n)$ 
and $\bar Q(n)$ described in \cref{burger2}
induce filtrations $K = K_0 \geq K_1 \geq \cdots \geq K_n = 0$ and $Q = Q_0 \geq Q_1 \geq \cdots$ with $K_{i-1} / K_i$ being
a (possibly zero) quotient of $q^{2i-n-1} \bar\Delta(n-1)$ 
for $i=1,\dots,n$, 
and $Q_{i-1} / Q_i$ being a (possibly
zero) quotient of $q^{2i-n-2} \bar\Delta(n+1)$
for $i \geq 1$.
By \cref{sundaymorning}, we know that 
$[V:L(n-1)]_q = [n]$.
Since $[Q:L(n-1)]_q = 0$,
these composition factors can only come
from the heads of $K_{i-1} / K_{i}$ for
$i=1,\dots,n$. So we must have
that $K_0 > K_1 > \cdots > K_n = 0$.
Since $K_i = x^i K$ by definition, 
this shows that $x^{n-1} \neq 0$.

Now suppose that $n \equiv t \pmod{2}$.
Then all composition factors of $V$ are isomorphic (up to degree shift) to $L(n-1)$
by \cref{sundaymorning} again.
We deduce that $V=K$, $V_i = K_i$
and $V_{i-1} / V_i \cong q^{2i-n-1}L(n-1)$
for each $i$. 
Thus, we have constructed the filtration described in (1).
We also  know from \cref{branch1}(3) 
that $\head V \cong q^{1-n} L(n-1)$
so that $\dim \End_{\NB}(V) \leq [V:L(n-1)] = n$.
As $x^{n-1} \neq 0$, the endomorphisms $1,x,\dots,x^{n-1}$ are linearly independent. So we have that
$\End_{\NB}(V) = \kk[x] / (x^n)$ as at the end of the statement of the theorem.
Moreover, $V$ is uniserial because $V$,
hence, each $V_i = x^i V$ has irreducible head,
i.e., $V_i$ is the unique maximal submodule $\rad V_{i-1}$ of $V_{i-1}$ for $i=1,\dots,n$.

It remains to treat the case 
$n \not\equiv t\pmod{2}$.
Since $\head V \cong q^{-n} L(n+1)$
and $[V:L(n+1)]_q = [n+1]$, we have that
$\dim \End_{\NB}(V) \leq [V:L(n+1)] = n+1$.
We know already that $x^{n-1} \neq 0$.
We cannot have $x^n = 0$ as this would contradict \cref{minpolydegree}.
So the nilpotency degree of $x$ is exactly
$n+1$, and $\End_{\NB}(V) = \kk[x] / (x^{n+1})$
as required for the final statement of the theorem.
It follows that
$V=V_0 > V_1 > \cdots > V_n > V_{n+1} = 0$.
Since $\head V \cong q^{-n} L(n+1)$, 
each $V_i$ 
has irreducible head
$q^{2i-n} L(n+1)$.
Since $\soc V \cong q^{n} L(n+1)$
we have that $V_n = \im x^n = \soc V$.
This is also the image of the restriction of 
$x^{n+1-i}$
to $V_{i-1}$, and $x^{n+1-i} V_{i} = 0$,
so $x^{n+1-i}$ induces a homomorphism
$V_{i-1} / V_{i} \twoheadrightarrow q^{2i-n-2}L(n+1)$. It follows that $V^{n+1-i}=\rad V_{i-1}$. 
We have now shown that
$$
V = V_0 > V^n \geq V_1 > V^{n-1} \geq V_2 > \cdots > V^1 \geq V_n > V^0 = 0
$$
with $V_{i-1} / V^{n+1-i} \cong q^{2i-n-2}L(n+1)$ for $i=1,\dots,n+1$.
We claim that $V^{n+1-i} / V_i$
has $q^{2i-n-1} L(n-1)$ 
as a composition factor.
This follows because $\head K_{i-1} \cong q^{2i-n-1} L(n-1)$,
$x^{n+1-i} K_{i-1} = 0$ and $x^{n-i} K_{i-1} \neq 0$, so
$V^{n+1-i} / V^{n-i}$ has $q^{2i-n-1} L(n-1)$ as a composition factor.
Combined with the information from 
\cref{sundaymorning}, the claim implies that
$V^{n+1-i} / V_i\cong q^{2i-n-1} L(n-1)$,
and we have constructed the filtration in (2).
Finally, we observe that $V$ is uniserial
because $V_{i-1}$ has irreducible
head $q^{2i-n-2} L(n+1)$
for $i=1,\dots,n+1$, hence,
$V_{i-1}/V_i$ is uniserial of length 2
for $i=1,\dots,n$ or length 1 for $i=n+1$.
\end{proof}